\documentclass{article}
\usepackage{amssymb}
\usepackage{amsmath}
\usepackage{sw20elba}

\newtheorem{theorem}{Theorem}

\newtheorem{condition}[theorem]{Condition}

\newtheorem{corollary}[theorem]{Corollary}

\newtheorem{definition}[theorem]{Definition}
\newtheorem{example}[theorem]{Example}

\newtheorem{lemma}[theorem]{Lemma}

\newtheorem{proposition}[theorem]{Proposition}
\newtheorem{remark}[theorem]{Remark}

\newenvironment{proof}[1][Proof]{\textbf{#1.} }{\ \rule{0.5em}{0.5em}}
\input{tcilatex}

\begin{document}

\title{Wavepacket preservation under nonlinear evolution }
\author{A. Babin and A. Figotin}
\maketitle

\begin{abstract}
We study nonlinear systems of hyperbolic (in a wider sense) PDE's in $%
\mathbb{R}^{d}$ describing wave propagation in dispersive nonlinear media
such as, for example, electromagnetic waves in nonlinear photonic crystals.
The initial data is assumed to be a finite sum of wavepackets referred to as
a \emph{multi-wavepacket}. The wavepackets and the medium nonlinearity are
characterized by two principal small parameters $\beta $ and $\varrho $
where: (i) $\frac{1}{\beta }$ is a factor describing spatial extension of
involved wavepackets; (ii) $\frac{1}{\varrho }$ is a factor describing the
relative magnitude of the linear part of the evolution equation compared to
its nonlinearity. A key element in our approach is a proper definition of a
wavepacket. Remarkably, the introduced definition has a flexibility
sufficient for a wavepacket to preserve its defining properties under a
general nonlinear evolution for long times. In particular, the corresponding
wave vectors and the band numbers of involved wavepackets are "conserved
quantities". We also prove that the evolution of a multi-wavepacket is
described with high accuracy by a properly constructed system of envelope
equations with a universal nonlinearity. The universal nonlinearity is
obtained by a time averaging applied to the original nonlinearity, in
simpler cases the averaged system turns into a system of Nonlinear
Schrodinger equations.
\end{abstract}

\section{Introduction}

The underlying physical subject of this work is propagation  of a \emph{%
multi-wavepacket} (a finite system of wavepackets) in a spatially dispersive
and nonlinear medium, and we are particularly interested in electromagnetic
waves propagation in nonlinear photonic crystals, see \cite{SipeBhat}, \cite%
{Slusher}, \cite{VolkSipe}, \cite{BF4}-\cite{BF7} and references therein,
with the nonlinear optics constitutive relations, \cite{Bo}, \cite[Sections
1,2]{ButCot}, \cite{Mil}, \cite{Sauter}. The mathematical subject of
interest is the following general nonlinear evolutionary system%
\begin{equation}
\partial _{\tau }\mathbf{U}=-\frac{\mathrm{i}}{\varrho }\mathbf{L}\left( -%
\mathrm{i}\nabla \right) \mathbf{U}+\mathbf{F}\left( \mathbf{U}\right) ,\
\left. \mathbf{U}\left( \mathbf{r},\tau \right) \right\vert _{\tau =0}=%
\mathbf{h}\left( \mathbf{r}\right) ,\ \mathbf{r}\in \mathbb{R}^{d},
\label{difeqintr}
\end{equation}%
where (i) $\mathbf{U}=\mathbf{U}\left( \mathbf{r},\tau \right) $, $\mathbf{r}%
\in \mathbb{R}^{d}$, $\mathbf{U}\in \mathbb{C}^{2J}$ is a $2J$ dimensional
vector; (ii) $\mathbf{L}\left( -\mathrm{i}\nabla \right) $ is a linear
self-adjoint differential (pseudodifferential) operator with constant
coefficients with the symbol $\mathbf{L}\left( \mathbf{k}\right) $, which is
a Hermitian $2J\times 2J$ matrix; (iii) $\ \mathbf{F}$ is a polynomial
nonlinearity such that $\mathbf{F}\left( \mathbf{0}\right) =\mathbf{0}$, $%
\mathbf{F}^{\prime }\left( \mathbf{0}\right) =\mathbf{0}$ and $\mathbf{F}%
\left( \mathbf{U}\right) $ is translation-invariant, i.e. if $T_{\mathbf{a}}%
\mathbf{U}\left( \mathbf{r}\right) =\mathbf{U}\left( \mathbf{r}+\mathbf{a}%
\right) $ for $\mathbf{a}\in \mathbb{R}^{d}$ then $\mathbf{F}\left( T_{%
\mathbf{a}}\mathbf{U}\right) =T_{\mathbf{a}}\mathbf{F}\left( \mathbf{U}%
\right) $; (iv) $\mathbf{h}=\mathbf{h}\left( \mathbf{r}\right) $ is assumed
to be the sum of a finite number of wavepackets $\mathbf{h}_{l}$, $%
l=1,\ldots ,N$; (v) $\varrho >0$ is a \emph{small parameter}. In the case of
nonlinear photonic crystals the components of the vector field $\mathbf{U}%
\left( \mathbf{r}\right) $ are the modal amplitudes of the electromagnetic
field and the nonlinearity $\mathbf{F}\left( \mathbf{U}\right) $ is
constructed from the nonlinear\ medium polarization in the adiabatic
approximation, \cite[Sections 2.4.2]{ButCot}. The systems of the form (\ref%
{difeqintr})  also describe as a particular case well-known equations,
namely:  \ complexification of the Nonlinear Schrodinger equation; coupled
envelope equations which arise in nonlinear birefringent optical media, \cite%
[Section 2i]{MolNew}; nonlinear Klein-Gordon and Sine-Gordon equations \cite[%
Section 14.1]{W}, \cite[Section 5.8.3]{Nayfeh}, \cite[Section 9.6]{OstrPot}.
Such equations appear in a number of physical problems: elementary
particles; dislocations in crystals; propagation of Bloch's domain walls in
the theory of ferromagnetism; self-induced transparency in nonlinear optics;
the propagation of magnetic flux quanta in long Josephson transmission
lines. Significance and importance of wavepacket solutions from the both
physical and mathematical points of view is discussed in \cite{BF4}-\cite%
{BF7}, \cite[Section 2]{MolNew}, \cite{SipeBhat}, \cite{VolkSipe}.

There are numerous problems involving small parameters only in the initial
data which can be reduced to the form (\ref{difeqintr}), for instance, 
problems with high frequency initial data or small initial data with
consequent evolution on long time intervals (see Section 3 for details).

We study the nonlinear evolution equation (\ref{difeqintr}) on a finite time
interval 
\begin{equation}
0\leq \tau \leq \tau _{\ast },\text{ where }\tau _{\ast }>0\text{ is a fixed
number.}  \label{tstarrho}
\end{equation}%
The time $\tau _{\ast }$ may depend on the $L^{\infty }$ norm of the initial
data $\mathbf{h}$ but, importantly, $\tau _{\ast }$\emph{\ does not depend
on }$\varrho $. \emph{We consider classes of initial data such that wave
evolution governed by (\ref{difeqintr}) is significantly nonlinear on time
interval }$\left[ 0,\tau _{\ast }\right] $\emph{\ and the effect of the
nonlinearity }$F\left( \mathbf{U}\right) $\emph{\ does not vanish as }$%
\varrho \rightarrow 0$\emph{.}

Since the both linear operator $\mathbf{L}\left( -\mathrm{i}\nabla \right) $
and the nonlinearity $\mathbf{F}\left( \mathbf{U}\right) $ are translation
invariant, it is natural and convenient to recast the evolution equation (%
\ref{difeqintr}) by applying to it the Fourier transform with respect to the
space variables $\mathbf{r}$, namely 
\begin{equation}
\partial _{\tau }\mathbf{\hat{U}}\left( \mathbf{k}\right) =-\frac{\mathrm{i}%
}{\varrho }\mathbf{L}\left( \mathbf{k}\right) \mathbf{\hat{U}}\left( \mathbf{%
k}\right) +\mathbf{\hat{F}}\left( \mathbf{\hat{U}}\right) \left( \mathbf{k}%
\right) ,\ \left. \mathbf{\hat{U}}\left( \mathbf{k}\right) \right\vert
_{\tau =0}=\mathbf{\hat{h}}\left( \mathbf{k}\right) ,  \label{difeqfou}
\end{equation}%
where $\mathbf{\hat{U}}\left( \mathbf{k}\right) $ is the Fourier transform
of $\mathbf{U}\left( \mathbf{r}\right) $, i.e.%
\begin{equation}
\mathbf{\hat{U}}\left( \mathbf{k}\right) =\int_{\mathbb{R}^{d}}\mathbf{U}%
\left( \mathbf{r}\right) \mathrm{e}^{-\mathrm{i}\mathbf{r}\cdot \mathbf{k}}%
\mathrm{d}\mathbf{r},\text{ }\mathbf{U}\left( \mathbf{r}\right) =\left( 2\pi
\right) ^{-d}\int_{\mathbb{R}^{d}}\mathbf{\hat{U}}\left( \mathbf{k}\right) 
\mathrm{e}^{\mathrm{i}\mathbf{r}\cdot \mathbf{k}}\mathrm{d}\mathbf{r},\text{%
\ where }\mathbf{r},\mathbf{k}\in \mathbb{R}^{d},  \label{Ftrans}
\end{equation}%
and $\mathbf{\hat{F}}$ is the Fourier form of the nonlinear operator $%
\mathbf{F}\left( \mathbf{U}\right) $ involving convolutions.

The nonlinear evolution equations (\ref{difeqintr}), (\ref{difeqfou}) are
commonly interpreted as describing wave propagation in a nonlinear medium.
We assume that the linear part $\mathbf{L}\left( \mathbf{\mathbf{k}}\right) $
is a $2J\times 2J$ Hermitian matrix with eigenvalues $\omega _{n,\zeta
}\left( \mathbf{k}\right) $ and eigenvectors $\mathbf{g}_{n,\zeta }\left( 
\mathbf{k}\right) $ satisfying%
\begin{equation}
\mathbf{L}\left( \mathbf{\mathbf{k}}\right) \mathbf{g}_{n,\zeta }\left( 
\mathbf{k}\right) =\omega _{n,\zeta }\left( \mathbf{k}\right) \mathbf{g}%
_{n,\zeta }\left( \mathbf{k}\right) ,\ \zeta =\pm ,\ \omega _{n,+}\left( 
\mathbf{k}\right) \geq 0,\ \omega _{n,-}\left( \mathbf{k}\right) \leq 0,\
n=1,\ldots ,J,  \label{OmomL}
\end{equation}%
where $\omega _{n,\zeta }\left( \mathbf{k}\right) $ are real-valued,
continuous for all non-singular $\mathbf{k}$ \ functions, and vectors $%
\mathbf{g}_{n,\zeta }\left( \mathbf{k}\right) \in \mathbb{C}^{2J}$ have unit
length in the standard Euclidean norm. The functions $\omega _{n,\zeta
}\left( \mathbf{k}\right) $, $n=1,\ldots ,J$, are called \emph{dispersion
relations} between the frequency $\omega $ and the \emph{wavevector} $%
\mathbf{k}$ with $n$ being the \emph{band number}. We assume that the
eigenvalues are naturally ordered by 
\begin{equation}
\omega _{J,+}\left( \mathbf{k}\right) \geq \ldots \geq \omega _{1,+}\left( 
\mathbf{k}\right) \geq 0\geq \omega _{1,-}\left( \mathbf{k}\right) \geq
\ldots \geq \omega _{J,-}\left( \mathbf{k}\right) ,  \label{omgr0}
\end{equation}%
and for \emph{almost every} $\mathbf{k}$ (with respect to the standard
Lebesgue measure) the eigenvalues are distinct and, consequently, the above
inequalities become strict. Importantly, we also assume the following \emph{%
diagonal symmetry} condition%
\begin{equation}
\omega _{n,-\zeta }\left( -\mathbf{k}\right) =-\omega _{n,\zeta }\left( 
\mathbf{k}\right) ,\ \zeta =\pm ,\ n=1,\ldots ,J,  \label{invsym}
\end{equation}%
which is naturally present in many physical problems (see also Remark \ref%
{Remark symmetry} below), and is a fundamental condition imposed on the
matrix $\mathbf{L}\left( \mathbf{\mathbf{k}}\right) $. In addition to that
in many examples we also have%
\begin{equation}
\mathbf{g}_{n,\zeta }\left( \mathbf{k}\right) =\mathbf{g}_{n,-\zeta }^{\ast
}\left( -\mathbf{k}\right) ,\text{ where }z^{\ast }\text{ is complex
conjugate to }z.  \label{omgr0a}
\end{equation}%
Very often we will use the following abbreviation 
\begin{equation}
\omega _{n,+}\left( \mathbf{k}\right) =\omega _{n}\left( \mathbf{k}\right) .
\label{ompl}
\end{equation}%
From (\ref{invsym}) \ we obtain%
\begin{equation}
\omega _{n,-}\left( \mathbf{k}\right) =-\omega _{n}\left( -\mathbf{k}\right)
,\;\omega _{n,\zeta }\left( \mathbf{k}\right) =\zeta \omega _{n}\left( \zeta 
\mathbf{k}\right) ,\;\zeta =\pm .  \label{omz}
\end{equation}%
We also will often use the orthogonal projection $\Pi _{n,\zeta }\left( 
\mathbf{\mathbf{k}}\right) $ in $\mathbb{C}^{2J}$ onto the complex line
defined by the eigenvector $\mathbf{g}_{n,\zeta }\left( \mathbf{k}\right) $,
namely 
\begin{equation}
\Pi _{n,\zeta }\left( \mathbf{\mathbf{k}}\right) \mathbf{\hat{u}}\left( 
\mathbf{k}\right) =\tilde{u}_{n,\zeta }\left( \mathbf{k}\right) \mathbf{g}%
_{n,\zeta }\left( \mathbf{k}\right) =\mathbf{\hat{u}}_{n,\zeta }\left( 
\mathbf{k}\right) ,\ n=1,\ldots ,J,\ \zeta =\pm .  \label{Pin}
\end{equation}

As it is indicated by the title of this paper we study the nonlinear problem
(\ref{difeqintr}) for initial data $\mathbf{\hat{h}}$ in the form of a
properly defined \emph{wavepacket} or, more generally, a sum of wavepackets
which we refer to as \emph{multi-wavepacket}. The simplest example of a
wavepacket $\mathbf{w}$ is provided by the following formula 
\begin{equation}
\mathbf{w}\left( \mathbf{r},\beta \right) =\Phi _{+}\left( \beta \mathbf{r}%
\right) \mathrm{e}^{\mathrm{i}\mathbf{k}_{\ast }\cdot \mathbf{r}}\mathbf{g}%
_{n,+}\left( \mathbf{k}_{\ast }\right) ,\ \mathbf{r}\in \mathbb{R}^{d},
\label{wpint}
\end{equation}%
where $\mathbf{k}_{\ast }\in \mathbb{R}^{d}$ is a \emph{wavepacket wave
vector}, $n$ is \emph{band number}, and $\beta >0$ is a small parameter. We
refer to the pair $\left( n,\mathbf{k}_{\ast }\right) $ in (\ref{wpint}) as 
\emph{wavepacket }$nk$\emph{-pair}. Observe that the space extension of the
wavepacket $\mathbf{w}\left( \mathbf{r},\beta \right) $ is proportional to $%
\beta ^{-1}$ and it is large for small $\beta $. Notice also that if $\beta
\rightarrow 0$ the wavepacket $\mathbf{w}\left( \mathbf{r},\beta \right) $
as in (\ref{wpint}) tends, up to a constant factor, to the elementary
eigenmode $\mathrm{e}^{\mathrm{i}\mathbf{k}_{\ast }\cdot \mathbf{r}}\mathbf{g%
}_{n,\zeta }\left( \mathbf{k}_{\ast }\right) $ of the operator $\mathbf{L}%
\left( -\mathrm{i}\nabla \right) $ with the corresponding eigenvalue $\omega
_{n,\zeta }\left( \mathbf{k}_{\ast }\right) $. We refer to wavepackets \ of
the simple form (\ref{wpint}) as \emph{simple wavepackets} to underline the
very special way the parameter $\beta $ enters its representation. The
function $\Phi _{\zeta }\left( \mathbf{r}\right) $, which we call \emph{%
wavepacket envelope}, describes its shape and it can be any scalar
complex-valued regular enough function, for example a function from Schwartz
space. Importantly, as $\beta \rightarrow 0$ the $L^{\infty }$ \emph{norm of
a wavepacket \ (\ref{wpint}) remains constant, and, hence, nonlinear effects
in (\ref{difeqintr}) remain strong}.

Evolution of wavepackets in problems which can be reduced to the form (\ref%
{difeqintr}) were studied for a variety of equations in numerous physical
and mathematical papers, mostly by asymptotic expansions with respect to a
single small parameter similar to $\beta $, see \cite{BenYoussefLannes02}, 
\cite{BonaCL05}, \cite{ColinLannes}, \cite{CraigSulemS92}, \cite{GiaMielke}, 
\cite{JolyMR98}, \cite{KalyakinUMN}, \cite{Maslov83}, \cite{PW}, \cite%
{Schneider98a}, \cite{Schneider05} and references therein. We are interested
in general properties of evolutionary systems of the form (\ref{difeqintr})
with wavepacket initial data which hold for a wide class of nonlinearities
and all values of the space dimensions $d$ of the number $2J$ of the system
components. Our approach is not based on asymptotic expansions but involves
the two small parameters $\beta $ and $\varrho $ with mild constraints on
their relative smallness. The constraints can be expressed either in the
form of certain inequalities or equalities, and a possible simple form of
such a constraint can be a power law 
\begin{equation}
\beta =C\varrho ^{\varkappa }\text{ where }C>0\text{ and }\varkappa >0\text{
are arbitrary constants.}  \label{powerk}
\end{equation}%
Of course, general features of wavepacket evolution are independent of
particular values of\ the constant $C$. In addition to that, some
fundamental properties such as wavepacket invariance, are also totally
independent on particular choice of the values of $\varkappa $ in (\ref%
{powerk}), whereas other properties are independent of $\varkappa $ as it
varies in certain intervals. For for instance, dispersion effects are
dominant for $\varkappa <1/2$, whereas the\ wavepacket superposition
principle of \cite{BF7} holds for $\varkappa <1$.

The qualitative picture of wavepacket evolution dependence on small $\beta $
and $\varrho $ is as follows. The parameter $\beta $ enters the problem (\ref%
{difeqintr}) through the multi-wavepacket initial data $\mathbf{h}\left( 
\mathbf{r},\beta \right) $ whereas $\varrho $ enters it through the factor $%
\frac{1}{\varrho }$ before the linear part. Evidently the factor $\frac{1}{%
\varrho }$ determines the relative magnitude of the linear part compared to
the nonlinearity and since $\frac{1}{\varrho }$ is large, one expects the
linear part to provide an important input into solutions properties. This
input includes, in particular, \emph{key role of eigenmodes and
eigenfrequencies (dispersion relations)} in expressing the nonlinear
evolution. Importantly, in many cases of interest though $\frac{1}{\varrho }$
is large, \emph{nonlinear phenomena are significant} and this is the case
when $\beta \leq C\varrho ^{1/2}$. More precisely, if $\beta \leq C\varrho
^{1/2}$ then, as in the case of finite-dimensional nonlinear ODE
evolutionary systems, the large values of $\frac{1}{\varrho }$ lead to a
well defined solution factorization into the fast (high frequency) and the
slow (low frequency) components. The interplay between the fast and slow
components is also similar to the ODE case, namely, the nonlinear evolution
is associated primarily with the slow component governed by a nonlinear
equation obtained from the original one by a certain canonical time
averaging procedure. Our further analysis of the above mentioned interplay
shows the following. Firstly, the \emph{linear superposition principle holds}%
, \cite{BF7}, that is if $\varkappa <1$ is as in (\ref{powerk}) and the
initial data is a sum of generic wavepackets then the solution is the sum of
the solutions for single involved wavepackets with precision $\frac{\varrho 
}{\beta ^{1+\epsilon }}$ with arbitrary small $\epsilon $. Secondly,
properly defined \emph{wavepackets and their linear combinations are
preserved} under the nonlinear evolution (\ref{difeqintr}), which is a
subject of this paper.

In the light of the above discussion we introduce the \emph{slow variable} $%
\mathbf{\hat{u}}\left( \mathbf{k},\tau \right) $ by the formula 
\begin{equation}
\mathbf{\hat{U}}\left( \mathbf{k},\tau \right) =\mathrm{e}^{-\frac{\mathrm{i}%
\tau }{\varrho }\mathbf{L}\left( \mathbf{k}\right) }\mathbf{\hat{u}}\left( 
\mathbf{k},\tau \right) ,  \label{Uu0}
\end{equation}%
and recast the equation (\ref{difeqfou}) for it as follows 
\begin{equation}
\ \partial _{\tau }\mathbf{\hat{u}}=\mathrm{e}^{\frac{\mathrm{i}\tau }{%
\varrho }\mathbf{L}}\mathbf{\hat{F}}\left( \mathrm{e}^{\frac{-\mathrm{i}\tau 
}{\varrho }\mathbf{L}}\mathbf{\hat{u}}\right) ,\ \left. \mathbf{\hat{u}}%
\right\vert _{\tau =0}=\mathbf{\hat{h}}.  \label{dfsF}
\end{equation}%
Then we obtain an integral form of (\ref{dfsF}) by integrating it with
respect to $\tau $: 
\begin{equation}
\mathbf{\hat{u}}=\mathcal{F}\left( \mathbf{\hat{u}}\right) +\mathbf{\hat{h}}%
,\ \mathcal{F}\left( \mathbf{\hat{u}}\right) =\mathcal{F}\left( \varrho
\right) \left( \mathbf{\hat{u}}\right) =\int_{0}^{\tau }\mathrm{e}^{\frac{%
\mathrm{i}\tau ^{\prime }}{\varrho }\mathbf{L}}\mathbf{\hat{F}}\left( 
\mathrm{e}^{\frac{-\mathrm{i}\tau ^{\prime }}{\varrho }\mathbf{L}}\mathbf{%
\hat{u}}\left( \tau ^{\prime }\right) \right) \,\mathrm{d}\tau ^{\prime }
\label{ubaseq}
\end{equation}%
with explicitly defined nonlinear polynomial integral operator $\mathcal{F}%
\left( \varrho \right) $, which depends on the parameter $\varrho $. This
operator is bounded uniformly with respect to $\varrho $ in the Banach space 
$E=C\left( \left[ 0,\tau _{\ast }\right] ,L^{1}\right) $ of functions $%
\mathbf{\hat{v}}\left( \mathbf{k},\tau \right) $, $0\leq \tau \leq \tau
_{\ast }$, with the norm 
\begin{equation}
\left\Vert \mathbf{\hat{v}}\left( \mathbf{k},\tau \right) \right\Vert
_{E}=\left\Vert \mathbf{\hat{v}}\left( \mathbf{k},\tau \right) \right\Vert
_{C\left( \left[ 0,\tau _{\ast }\right] ,L^{1}\right) }=\sup_{0\leq \tau
\leq \tau _{\ast }}\int_{\mathbb{R}^{d}}\left\vert \mathbf{\hat{v}}\left( 
\mathbf{k},\tau \right) \right\vert \,\mathrm{d}\mathbf{k},  \label{Elat}
\end{equation}%
where $L^{1}$ is the Lebesgue space of functions $\mathbf{\hat{v}}\left( 
\mathbf{k}\right) $ with the standard norm 
\begin{equation}
\left\Vert \mathbf{\hat{v}}\left( \mathbf{\cdot }\right) \right\Vert
_{L^{1}}=\int_{\mathbb{R}^{d}}\left\vert \mathbf{\hat{v}}\left( \mathbf{k}%
\right) \right\vert \,\mathrm{d}\mathbf{k}.  \label{L1}
\end{equation}%
Sometimes we use more general weighted spaces $L^{1,a}$ with the norm 
\begin{equation}
\left\Vert \mathbf{\hat{v}}\right\Vert _{L^{1,a}}=\int_{\mathbb{R}%
^{d}}\left( 1+\left\vert \mathbf{k}\right\vert \right) ^{a}\left\vert 
\mathbf{\hat{v}}\left( \mathbf{k}\right) \right\vert \,\mathrm{d}\mathbf{k}%
,\ a\geq 0.  \label{L1a}
\end{equation}%
A rather elementary existence and uniqueness theorem (Theorem \ref{Theorem
Existence1}) implies that for a small and, importantly, independent of $%
\varrho $ constant $\tau _{\ast }>0$ this equation has a unique solution 
\begin{equation}
\mathbf{\hat{u}}\left( \tau \right) =\mathcal{G}\left( \mathcal{F}\left(
\varrho \right) ,\mathbf{\hat{h}}\right) \left( \tau \right) ,\ \tau \in %
\left[ 0,\tau _{\ast }\right] ,\;\mathbf{\hat{u}}\in C^{1}\left( \left[
0,\tau _{\ast }\right] ,L^{1}\right) ,  \label{ugh1}
\end{equation}%
where $\mathcal{G}$ denotes the \emph{solution operator} for the equation (%
\ref{ubaseq}), the operator depends on operator $\mathcal{F}\left( \varrho
\right) $, which itself depends on the parameter $\varrho $. If $\mathbf{%
\hat{u}}\left( \mathbf{k},\tau \right) $ is a solution to the equation (\ref%
{ubaseq}) we call the function $\mathbf{U}\left( \mathbf{r},\tau \right) $
defined by (\ref{Uu0}), (\ref{Ftrans}) an $\emph{F}$\emph{-solution} to the
equation (\ref{difeqintr}). We denote by $\hat{L}^{1}$ the space of
functions $\mathbf{V}\left( \mathbf{r}\right) $ such that their Fourier
transform $\mathbf{\hat{V}}\left( \mathbf{k}\right) $ belongs to $L^{1}$,
and define $\left\Vert \mathbf{V}\right\Vert _{\hat{L}^{1}}$ $=\left\Vert 
\mathbf{\hat{V}}\right\Vert _{L^{1}}$. Since 
\begin{equation}
\left\Vert \mathbf{V}\right\Vert _{L^{\infty }}\leq \left( 2\pi \right)
^{-d}\left\Vert \mathbf{\hat{V}}\right\Vert _{L^{1}}\text{ and }\hat{L}%
^{1}\subset L^{\infty }  \label{Linf}
\end{equation}%
$F$-solutions to (\ref{difeqintr}) belong to $C^{1}\left( \left[ 0,\tau
_{\ast }\right] ,\hat{L}^{1}\right) \subset C^{1}\left( \left[ 0,\tau _{\ast
}\right] ,L^{\infty }\right) $.

We would like to define wavepackets in a form which explicitly allows them
to be real valued. This is accomplished based on the symmetry (\ref{invsym})
of the dispersion relations by introduction of a \emph{doublet wavepacket} 
\begin{equation}
\mathbf{w}\left( \mathbf{r},\beta \right) =\Phi _{+}\left( \beta \mathbf{r}%
\right) \mathrm{e}^{\mathrm{i}\mathbf{k}_{\ast }\cdot \mathbf{r}}\mathbf{g}%
_{n,+}\left( \mathbf{k}_{\ast }\right) +\Phi _{-}\left( \beta \mathbf{r}%
\right) \mathrm{e}^{-\mathrm{i}\mathbf{k}_{\ast }\cdot \mathbf{r}}\mathbf{g}%
_{n,-}\left( -\mathbf{k}_{\ast }\right) .  \label{wpn1}
\end{equation}%
Such a wavepacket is real if $\Phi _{-}\left( \mathbf{r}\right) ,\mathbf{g}%
_{n,-}\left( -\mathbf{k}_{\ast }\right) $ is complex conjugate to $\Phi
_{+}\left( \mathbf{r}\right) $, $\mathbf{g}_{n,+}\left( \mathbf{k}_{\ast
}\right) $, i.e. if%
\begin{equation}
\Phi _{-}\left( \mathbf{r}\right) =\Phi _{+}^{\ast }\left( \mathbf{r}\right)
,\ \mathbf{g}_{n,+}\left( \mathbf{k}_{\ast }\right) =\mathbf{g}_{n,-}\left( -%
\mathbf{k}_{\ast }\right) ^{\ast }.  \label{wpn1a}
\end{equation}%
Considering wavepackets with $\ nk$\emph{-}pair $\left( n,\mathbf{k}_{\ast
}\right) $ we usually mean doublet ones as in (\ref{wpn1}), but sometimes $%
\Phi _{+}$ or $\Phi _{-}$ may be zero producing (\ref{wpint}).

To identify characteristic properties of a wavepacket suitable for our
needs, let us look at the Fourier transform $\mathbf{\hat{w}}\left( \mathbf{k%
},\beta \right) $ of an elementary wavepacket $\mathbf{w}\left( \mathbf{r}%
,\beta \right) $ defined by (\ref{wpint}), that is 
\begin{equation}
\mathbf{\hat{w}}\left( \mathbf{k},\beta \right) =\beta ^{-d}\hat{\Phi}\left(
\beta ^{-1}\left( \mathbf{k}-\mathbf{k}_{\ast }\right) \right) \mathbf{g}%
_{n,\zeta }\left( \mathbf{k}_{\ast }\right) .  \label{wpint1}
\end{equation}%
We call such $\mathbf{\hat{w}}\left( \mathbf{k},\beta \right) $ wavepacket
too, obviously it possesses the following properties: (i) its $L^{1}$ norm
is bounded (in fact, constant), uniformly in $\beta \rightarrow 0$; (ii) for
every $\epsilon >0$ the value $\mathbf{\hat{w}}\left( \mathbf{k},\beta
\right) \rightarrow 0$ for every $\mathbf{k}$ outside a $\beta ^{1-\epsilon
} $-neighborhood of $\mathbf{k}_{\ast }$, and the convergence is faster than
any power of $\beta $ if $\Phi $ is a Schwartz function. To explicitly
interpret the last property we introduce a \emph{cutoff function} $\Psi
\left( \mathbf{\eta }\right) $%
\begin{equation}
\Psi \left( \mathbf{\eta }\right) =1\text{ for }\left\vert \mathbf{\eta }%
\right\vert \leq 1,\ \Psi \left( \mathbf{\eta }\right) =0\text{ for }%
\left\vert \mathbf{\eta }\right\vert >1,  \label{j0}
\end{equation}%
together with its shifted/rescaled modification 
\begin{equation}
\Psi \left( \mathbf{k};\mathbf{k}_{\ast }\right) =\Psi \left( \mathbf{k};%
\mathbf{k}_{\ast },\beta ^{1-\epsilon }\right) =\Psi \left( \beta ^{-\left(
1-\epsilon \right) }\left( \mathbf{k}-\mathbf{k}_{\ast }\right) \right) .
\label{Psik}
\end{equation}%
If in an elementary wavepacket $\mathbf{w}\left( \mathbf{r},\beta \right) $
defined by (\ref{wpint1}) $\Phi _{\zeta }\left( \mathbf{r}\right) $ is a
Schwartz function then%
\begin{equation*}
\left\Vert \left( 1-\Psi \left( \mathbf{\cdot },\mathbf{k}_{\ast },\beta
^{1-\epsilon }\right) \right) \mathbf{\hat{w}}\left( \mathbf{\cdot },\beta
\right) \right\Vert \leq C_{\epsilon ,s}\beta ^{s},\ 0<\beta \leq 1,
\end{equation*}%
which holds for arbitrarily small $\epsilon >0$ and arbitrarily large $s>0$.
Based on the above discussion we give the following definition of a
wavepacket which is a minor variation of \cite[Definiton 8]{BF7}.

\begin{definition}[single-band wavepacket]
\label{dwavepack} Let $0<\epsilon <1$ \ be a fixed number. For a given band
number $n\in \left\{ 1,\ldots ,J\right\} $ and a wavevector $\mathbf{k}%
_{\ast }\in \mathbb{R}^{d}$ a function $\mathbf{\hat{h}}\left( \beta ,%
\mathbf{k}\right) $ is called a \emph{\ wavepacket with }$nk$\emph{-pair } $%
\left( n,\mathbf{k}_{\ast }\right) $ and the degree of regularity $s>0$ if
there exists such $\beta _{0}>0$ \ that for $\beta <\beta _{0}$ the
following conditions are satisfied: (i) $\mathbf{\hat{h}}\left( \beta ,%
\mathbf{k}\right) $ is $L^{1}$-bounded uniformly in $\beta $, i.e.%
\begin{equation}
\left\Vert \mathbf{\hat{h}}\left( \beta ,\mathbf{\cdot }\right) \right\Vert
_{L^{1}}\leq C,\ 0<\beta <\beta _{0}\text{ for some }C>0;  \label{L1b}
\end{equation}%
(ii) $\mathbf{\hat{h}}\left( \beta ,\mathbf{k}\right) $\ has the following
structure:\ 
\begin{equation}
\mathbf{\hat{h}}\left( \beta ,\mathbf{k}\right) =\mathbf{\hat{h}}_{-}\left(
\beta ,\mathbf{k}\right) +\mathbf{\hat{h}}_{+}\left( \beta ,\mathbf{k}%
\right) +\hat{D}_{h},\ 0<\beta <\beta _{0},\text{ where}  \label{hbold}
\end{equation}%
\begin{equation}
\mathbf{\hat{h}}_{\zeta }\left( \beta ,\mathbf{k}\right) =\Psi \left( 
\mathbf{k},\zeta \mathbf{k}_{\ast },\beta ^{1-\epsilon }\right) \Pi
_{n,\zeta }\left( \mathbf{k}\right) \mathbf{\hat{h}}_{\zeta }\left( \beta ,%
\mathbf{k}\right) ,\ \zeta =\pm ,  \label{halfwave}
\end{equation}%
with $\Psi \left( \mathbf{\cdot },\zeta \mathbf{k}_{\ast },\beta
^{1-\epsilon }\right) $ defined by (\ref{Psik}) and $\hat{D}_{h}$ satisfying
the following \emph{tail estimate}: 
\begin{equation}
\left\Vert \hat{D}_{h}\right\Vert _{L^{1}}\leq C^{\prime }\beta ^{s},\
0<\beta <\beta _{0}\text{ for some }C^{\prime }>0.  \label{sourloc}
\end{equation}%
The inverse Fourier transform $\mathbf{h}\left( \beta ,\mathbf{r}\right) $
of a wavepacket $\mathbf{\hat{h}}\left( \beta ,\mathbf{k}\right) $ is also
called a wavepacket.
\end{definition}

Point (ii) of the above Definition means that the wavepacket $\mathbf{\hat{h}%
}\left( \beta ,\mathbf{k}\right) $ is composed of two functions $\mathbf{%
\hat{h}}_{\zeta }\left( \beta ,\mathbf{k}\right) $, $\zeta =\pm $, which
take values in $n$-th band eigenspace of $\mathbf{L}\left( \mathbf{k}\right) 
$ and are localized near $\zeta \mathbf{k}_{\ast }$, where $\left( n,\mathbf{%
k}_{\ast }\right) $ is the $nk$\emph{-}pair \ of the wavepacket. The number $%
\beta _{0}$ usually is small and may depend on a wavepacket.

Evidently, if a wavepacket has the degree of regularity $s$, it also has a
smaller degree of regularity $s^{\prime }\leq s$ with the same $\epsilon $.
Observe that the degree of regularity $s$ is related to the smoothness of $%
\Phi _{\zeta }\left( \mathbf{r}\right) $ in (\ref{wpint}) so that the higher
is the smoothness the higher $\frac{s}{\epsilon }$ can be taken. Namely, if $%
\hat{\Phi}_{\zeta }\in L^{1,a}$ then one can take any $\frac{s}{\epsilon }<a$%
, see Lemma \ref{Lemma Lpsiwp} below. For example, if in the elementary
wavepacket $\mathbf{w}\left( \mathbf{r},\beta \right) $ defined by (\ref%
{wpint}) $\Phi _{\zeta }\left( \mathbf{r}\right) $ is a Schwartz function
then it has arbitrarily large degree of regularity.

It turns out remarkably that wavepackets satisfying Definition \ref%
{dwavepack} preserve their defining properties under nonlinear evolution. It
is remarkable, in particular, since it is well-known that determination of
classes of solutions which preserve their form under generic nonlinear
evolution usually leads to infinite expansions, such as multi-scale
expansions, power expansions, modal expansions etc with serious difficulties
in establishing the convergence. Such expansions often are formally
invariant, but they involve infinitely many rather complex terms and
establishing the convergence is a very hard problem indeed if there is any
convergence at all. Our Definition \ref{dwavepack} of a wavepacket involves
only a finite number of terms and its invariance is provided by the flexible 
\emph{tail term} $\hat{D}_{h}$. We also find remarkable the very simplicity
of the definition which nevetherless allows for a sufficiently detailed
analysis of the dynamics, including, in particular, rigorously justified
NLS-type approximations of wavepacket dynamics presented in the following
sections.

Our special interest is in waves that are finite sums of wavepackets and we
refer to them as \emph{multi-wavepackets}.

\begin{definition}[multi-wavepacket]
\label{dmwavepack} Let $S$ be a set of $nk$-pairs: 
\begin{equation}
S=\left\{ \left( n_{l},\mathbf{k}_{\ast l}\right) ,\ l=1,\ldots ,N\right\}
\subset \Sigma =\left\{ 1,\ldots ,J\right\} \times \mathbb{R}^{d},\ \left(
n_{l},\mathbf{k}_{\ast l}\right) \neq \left( n_{l^{\prime }},\mathbf{k}%
_{\ast l^{\prime }}\right) \text{ for }l\neq l^{\prime },  \label{P0}
\end{equation}%
and $N=\left\vert S\right\vert $ be their number. Let $K_{S}$ be a set
consisting of all different wavevectors $\mathbf{k}_{\ast l}$ involved in $S$
with $\left\vert K_{S}\right\vert \leq N$ being the number of its elements. $%
K_{S}$ is called \emph{wavepacket }$k$\emph{-spectrum} and without loss of
genericity we assume the indexing of elements in $S$ to be such that 
\begin{equation}
K_{S}=\left\{ \mathbf{k}_{\ast i},i=1,\ldots ,\left\vert K_{S}\right\vert
\right\} ,\text{ i.e. }l_{i}=i\text{ for }1\leq i\leq \left\vert
K_{S}\right\vert \text{.}  \label{K0}
\end{equation}%
A function $\mathbf{\hat{h}}\left( \beta \right) =\mathbf{\hat{h}}\left(
\beta ,\mathbf{k}\right) $ is called a \emph{multi-wavepacket} with $nk$%
\emph{-spectrum} $S$ if it is a finite sum of wavepackets, namely 
\begin{equation}
\mathbf{\hat{h}}\left( \beta ,\mathbf{k}\right) =\sum_{l=1}^{N}\mathbf{\hat{h%
}}_{l}\left( \beta ,\mathbf{k}\right) ,\ 0<\beta <\beta _{0}\text{ for some }%
\beta _{0}>0,  \label{JJ1}
\end{equation}%
where $\mathbf{\hat{h}}_{l}$, $l=1,\ldots ,N$, is a wavepacket with $nk$%
-pair $\left( \mathbf{k}_{\ast l},n_{l}\right) \in S$ as in Definition \ref%
{dwavepack}.
\end{definition}

Note that if $\mathbf{\hat{h}}\left( \beta ,\mathbf{k}\right) $ is a
wavepacket then $\mathbf{\hat{h}}\left( \beta ,\mathbf{k}\right) +O\left(
\beta ^{s}\right) $ is a wavepacket as well\ with the same $nk$-spectrum,
and the same is true for multi-wavepackets. Hence, we can introduce
multi-wavepackets equivalence relation "$\simeq $" of the degree $s$ by 
\begin{equation}
\mathbf{\hat{h}}_{1}\left( \beta ,\mathbf{k}\right) \simeq \mathbf{\hat{h}}%
_{2}\left( \beta ,\mathbf{k}\right) \text{ if }\left\Vert \mathbf{\hat{h}}%
_{1}\left( \beta ,\mathbf{k}\right) -\mathbf{\hat{h}}_{2}\left( \beta ,%
\mathbf{k}\right) \right\Vert _{L^{1}}\leq C\beta ^{s}\text{ for some
constant }C>0\text{.}  \label{equiv}
\end{equation}

Let us turn now to the abstract nonlinear problem (\ref{ubaseq}) where (i) $%
\mathcal{F}=\mathcal{F}\left( \varrho \right) $ depends on $\varrho $ and
(ii) the initial data $\mathbf{\hat{h}}=\mathbf{\hat{h}}\left( \beta \right) 
$\ is a multi-wavepacket depending on $\beta $. We would like to state our
first theorem on multi-wavepacket preservation under the evolution (\ref%
{ubaseq}) for $\beta ,\varrho \rightarrow 0$, which holds, as it turns out,
provided its $nk$\emph{-}spectrum $S$ satisfies certain natural condition
called \emph{resonance invariance}. This condition is intimately related to
the so-called \emph{phase and frequency matching} conditions for stronger
nonlinear interactions, and its concise formulation is as follows. We define
for given dispersion relations $\left\{ \omega _{n}\left( \mathbf{k}\right)
\right\} $ and any finite set $S\subset \left\{ 1,\ldots ,J\right\} \times 
\mathbb{R}^{d}$ another finite set $\mathcal{R}\left( S\right) \subset
\left\{ 1,\ldots ,J\right\} \times \mathbb{R}^{d}$ where $\mathcal{R}$ is a
certain algebraic operation described in Definition \ref{Definition omclos}
below. It turns out that for any $S$ always $S\subseteq \mathcal{R}\left(
S\right) $ but \emph{if, in fact, }$\mathcal{R}\left( S\right) =S$\emph{\ we
call }$S$\emph{\ resonance invariant}. The condition of resonance invariance
is instrumental for the multi-wavepacket preservation, and there are
examples showing that if it fails, i.e. $\mathcal{R}\left( S\right) \neq S$,
the wavepacket preservation does not hold. \emph{\ Importantly, the
resonance invariance }$\mathcal{R}\left( S\right) =S$\emph{\ allows
resonances inside the multi-wavepacket, that includes, in particular,
resonances associated with the second and the third harmonic generations,
resonant four-wave interaction etc. }

\begin{theorem}[multi-wavepacket preservation]
\label{Theorem invarwave}Suppose that the nonlinear evolution is governed by
(\ref{ubaseq}) and the initial data $\mathbf{\hat{h}}=\mathbf{\hat{h}}\left(
\beta ,\mathbf{k}\right) $ is a multi-wavepacket with $nk$-spectrum $S$ and
the regularity degree $s$, and assume $S$ to be resonance invariant (see
Definition \ref{Definition omclos} below). Let dependence between parametrs $%
\varrho $ and $\beta \ $be given be any function $\varrho =\rho \left( \beta
\right) $ satisfying 
\begin{equation}
0<\rho \left( \beta \right) \leq C\beta ^{s},\text{ for some constant }C>0,
\label{rbr}
\end{equation}%
and let us set $\varrho =\rho \left( \beta \right) $. Then the solution $%
\mathbf{\hat{u}}\left( \tau ,\beta \right) =\mathcal{G}\left( \mathcal{F}%
\left( \rho \left( \beta \right) \right) ,\mathbf{\hat{h}}\left( \beta
\right) \right) \left( \tau \right) $ to (\ref{ubaseq}) for any $\tau \in %
\left[ 0,\tau _{\ast }\right] $ is a multi-wavepacket with $nk$-spectrum $S$
and the regularity degree $s$, i.e.%
\begin{equation}
\mathbf{\hat{u}}\left( \tau ,\beta ;\mathbf{k}\right)
=\dsum\nolimits_{l=1}^{N}\mathbf{\hat{u}}_{l}\left( \tau ,\beta ;\mathbf{k}%
\right) ,\text{ where }\mathbf{\hat{u}}_{l}\text{ is wavepacket with }nk%
\text{-pair }\left( n_{l},\mathbf{k}_{\ast l}\right) \in S.  \label{usumul}
\end{equation}%
The time interval length $\tau _{\ast }>0$ depends only on $L^{1}$-norms of $%
\mathbf{\hat{h}}_{l}\left( \beta ,\mathbf{k}\right) $ and $N$. The
presentation (\ref{usumul}) is unique up to the equivalence (\ref{equiv}).
\end{theorem}

The above statement can be interpreted as follows. Modes in $nk$-spectrum $S$
are\ always resonance coupled with modes in $\mathcal{R}\left( S\right) $
through the nonlinear interactions, but if $\mathcal{R}\left( S\right) =S$
then (i) all resonance interactions occur inside $S$ and (ii) only small
vicinity of $S$ is involved in nonlinear interactions leading to the
multi-wavepacket preservation.

Many nonlinear evolution problems with small initial data can be readily
reduced by elementary rescaling to the system (\ref{difeqintr}) with a large
parameter $\frac{1}{\varrho }$ before its linear part. For example, suppose
that $\mathbf{F}\left( \mathbf{V}\right) $ is a homogeneous nonlinearity of
degree $m$ ($m=3$ for cubic one) and\ that the nonlinear evolution is
governed by 
\begin{equation}
\partial _{t}\mathbf{V}=-\mathrm{i}\mathbf{L}\left( -\mathrm{i}\nabla
\right) \mathbf{V}+\mathbf{F}\left( \mathbf{V}\right) ,\ \left. \mathbf{V}%
\left( \mathbf{r},t\right) \right\vert _{t=0}=\varrho ^{1/\left( m-1\right) }%
\mathbf{h}\left( \mathbf{r}\right) ,\ \mathbf{r}\in \mathbb{R}^{d},
\label{difeqintrv}
\end{equation}%
considered for small $\varrho $ on the large time interval $0\leq t\leq 
\frac{\tau _{\ast }}{\varrho }$ with a fixed $\tau _{\ast }>0$. Then the
following simple change of variables 
\begin{equation}
\mathbf{V}\left( t\right) =\varrho ^{1/\left( m-1\right) }\mathbf{U}\left(
\tau \right) ,\ \tau =t\varrho  \label{VU}
\end{equation}%
transforms the problem (\ref{difeqintrv}) into equivalent problem (\ref%
{difeqintr}). In this case the inequality (\ref{rbr}) describes a constraint
between spatial extension $\frac{1}{\beta }$ and the amplitude factor $%
\varrho ^{1/\left( m-1\right) }=\rho \left( \beta \right) ^{1/\left(
m-1\right) }$ of the initial data. Observe that the equation (\ref%
{difeqintrv}) does not have any small parameters and the both small
parameters $\varrho $ and $\beta $ enter the problem through its initial
data. Theorem \ref{Theorem invarwave} \ can be restated for the problem (\ref%
{difeqintrv}) as follows.

\begin{corollary}[multi-wavepacket preservation]
\label{cinvar} Let $\mathbf{V}\left( \mathbf{r},t\right) $ be a solution to
the nonlinear system (\ref{difeqintrv}), $\rho \left( \beta \right) $ is as
in (\ref{rbr}) and we set $\varrho =\rho \left( \beta \right) $. Then if the
initial data is such that $\varrho ^{-1/\left( m-1\right) }\mathbf{\hat{V}}%
\left( \mathbf{k},0\right) =\mathbf{\hat{h}}\left( \mathbf{k}\right) $ is a
multi-wavepacket then $\varrho ^{-1/\left( m-1\right) }\mathbf{\hat{V}}%
\left( \mathbf{k},t\right) $ remains to be a multi-wavepacket with the same $%
nk$-spectrum and the degree of regularity for all times $t\in \left[ 0,\frac{%
\tau _{\ast }}{\varrho }\right] $.
\end{corollary}

The statements of Theorems \ref{Theorem invarwave} and Corollary \ref{cinvar}
directly follow from the following general theorem which makes no
assumptions on the relations between $\beta ,\varrho \rightarrow 0$.

\begin{theorem}[multi-wavepacket approximation]
\label{Theorem sumwave} Let the initial data $\mathbf{\hat{h}}$ in the
integral equation (\ref{ubaseq}) be a multi-wavepacket $\mathbf{\hat{h}}%
\left( \beta ,\mathbf{k}\right) $ with $nk$-spectrum $S$ as in (\ref{P0}),
the regularity degree $s$ and with the parameter $\epsilon >0$ as in
Definition \ref{dwavepack}. \ Assume that $S$ is resonance invariant in the
sense of Definition \ref{Definition omclos} below. Let the cutoff function $%
\Psi \left( \mathbf{k},\mathbf{k}_{\ast }\right) $ and the eigenvector
projectors $\Pi _{n,\pm }\left( \mathbf{\mathbf{k}}\right) $ be defined by (%
\ref{Psik}) and (\ref{Pin}) respectively. For a solution $\mathbf{\hat{u}}$
of (\ref{ubaseq}) \ we set%
\begin{equation}
\mathbf{\hat{u}}_{l}\left( \tau ,\beta ;\mathbf{k}\right) =\left[
\dsum\nolimits_{\zeta =\pm }\Psi \left( \mathbf{k},\zeta \mathbf{k}_{\ast
l}\right) \Pi _{n_{l},\zeta }\left( \mathbf{\mathbf{k}}\right) \right] 
\mathbf{\hat{u}}\left( \tau ,\beta ;\mathbf{k}\right) ,\ l=1,\ldots ,N.
\label{ups1}
\end{equation}%
Then every such $\mathbf{\hat{u}}_{l}\left( \mathbf{k};\tau ,\beta \right) $
is a wavepacket and%
\begin{equation}
\sup_{0\leq \tau \leq \tau _{\ast }}\left\Vert \mathbf{\hat{u}}\left( \tau
,\beta ;\mathbf{k}\right) -\dsum\nolimits_{l=1}^{N}\mathbf{\hat{u}}%
_{l}\left( \tau ,\beta ;\mathbf{k}\right) \right\Vert _{L^{1}}\leq
C_{1}\varrho +C_{2}\beta ^{s}  \label{uui}
\end{equation}%
where the constant $C_{1}$ does not depend on $\epsilon ,s$ and $\beta $,
and the constant $C_{2}$ does not depend on $\beta .$
\end{theorem}

It is interesting to note that the statement of the Theorem \ref{Theorem
sumwave} can be extended to the special limit case $\beta =0$, $\mathbf{k}%
_{\ast l}=0$. In this case the initial data of (\ref{difeqintr}) are
constants in $\mathbf{r}$ and we can consider solutions $\mathbf{U}$ (\ref%
{difeqintr}) which do not depend on $\mathbf{r}$. Then $\nabla \mathbf{U=0}$%
, the linear operator $\mathbf{L}\left( -\mathrm{i}\nabla \right) $ reduces
to the multiplication by a matrix $\mathbf{L}_{0}=\mathbf{L}\left( 0\right) $
and the system (\ref{difeqintr}) turns into a system of ordinary
differential equations (ODE). Notice that (i) the structure of the
eigenvalues (\ref{invsym}) implies that the linear part is time-reversible;
(ii) the nonlinear part can be an arbitrary polynomial. The extension of
Theorem \ref{Theorem sumwave} to this case (see Theorem \ref{Theorem
sumwave1}) reads that in a generic, non-resonant situation if initial data
are bounded and a set of eigenmodes of the matrix $\mathbf{L}_{0}$ is
excited at $\tau =0$ then in the course of evolution on a time interval $%
\left[ 0,\tau _{\ast }\right] $ where $\tau _{\ast }$ depends on magnitude
of initial data (i) all remaining modes remain unexcited with accuracy
proportional to $\varrho $, and (ii) only the originally excited modes can
significantly evolve with this level of accuracy. For finite-dimensional
systems governed by ODE's such a statement can be derived from the classical
time-averaging principle and the time-averaged equations remain nonlinear.
For infinitely-dimensional systems governed by PDE and with the linear
operator having a continuous spectrum, as in Theorem \ref{Theorem sumwave},
the analysis is more complex but the time-averaging still plays important
role yielding an accurate approximation governed by a certain universal
nonlinear PDE.

We would like to point out also that though Theorem \ref{Theorem invarwave}
is a simple corollary of more general Theorem \ref{Theorem sumwave}, \emph{%
it is important that the statement (\ref{uui}) can be formulated as
multi-wavepacket invariance}. That, in particular, allows to take values $%
\mathbf{\hat{u}}\left( \tau _{\ast }\right) $ as new wavepacket initial data
for (\ref{difeqintr}) and extend the wavepacket \ invariance of a solution
to the next time interval $\tau _{\ast }\leq \tau \leq \tau _{\ast 1}$. This
observation allows to extend the wavepacket \ invariance to larger values of 
$\tau $ (up to blow-up time or infinity) if some additional information
about solutions with wavepacket initial data is available. In particular,
the following theorem holds.

\begin{theorem}
\label{Theorem globalinv}Assume that all conditions of Theorem \ref{Theorem
invarwave} are satisfied and, in addition to that, solutions $\mathbf{\hat{u}%
}\left( \tau \right) $ of (\ref{ubaseq}) with the multi-wavepacket initial
data $\mathbf{\hat{h}}\left( \beta \right) $ exist on a maximal interval $%
0\leq \tau <\tau _{0}$ such that $\left\Vert \mathbf{\hat{u}}\right\Vert
_{C\left( \left[ 0,\tau _{1}\right] ,L^{1}\right) }\leq R\left( \tau
_{1}\right) $ for any $\tau _{1}<\tau _{0}\leq \infty $ where $R\left( \tau
_{1}\right) $ does not depend on $\beta ,\varrho $. Then the solution $%
\mathbf{\hat{u}}\left( \tau ,\beta \right) =\mathcal{G}\left( \mathcal{F}%
\left( \rho \left( \beta \right) \right) ,\mathbf{\hat{h}}\left( \beta
\right) \right) \left( \tau \right) $ to (\ref{ubaseq}) for any $\tau <\tau
_{0}$ is a multi-wavepacket with $nk$-spectrum $S$ and the regularity degree 
$s$, that is (\ref{usumul}) holds.
\end{theorem}

Note that the wavepacket form of solutions can be used to obtain long-time
estimates of solutions. Namely, very often behavior of every single
wavepacket is well approximated by its own nonlinear Schrodinger equation
(NLS), see \cite{Colin}, \cite{KSM}, \cite{ColinLannes}, \cite{GiaMielke}, 
\cite{KalyakinUMN}, \cite{Kalyakin2}, \cite{PW}, \cite{Schneider98a}, \cite%
{Schneider05}, \cite{SU} and references therein, see also Section 6. Many
features of the dynamics governed by NLS-type equations are well-understood,
see \cite{Bourgain}, \cite{Caz}, \cite{SchlagK}, \cite{Schlag}, \cite{Sulem}%
, \cite{Weinstein} and references therein. These results can be used to
obtain long-time estimates for every single wavepacket (as, for example, in 
\cite{Kalyakin2}) and, with the help of the superposition principle, for the
multiwavepacket solution.

The wavepacket representation (\ref{usumul}) from Theorem \ref{Theorem
invarwave} can be used for more detailed analysis of dynamics of wavepackets 
$\mathbf{\hat{u}}_{l}\left( \tau ,\beta \right) $ and interaction between
them. The following theorem illustrates that by describing wavepacket
interaction based on a system with a weakly universal nonlinearity similar
to so-called coupled modes systems or NLS.

\begin{theorem}[NLS-type approximation]
\label{Theorem sumcanwave} Let the conditions of Theorem \ref{Theorem
sumwave} hold and, in addition to that, the initial data $\mathbf{\hat{h}}%
_{l}\left( \mathbf{k}\right) $ are of the form $\ \mathbf{\hat{h}}_{l}=%
\mathbf{\hat{h}}_{l,+}+\mathbf{\hat{h}}_{l,-}+\hat{D}_{l}$ where 
\begin{equation*}
\mathbf{\hat{h}}_{l,\zeta }\left( \mathbf{k}\right) =\beta ^{-d}\hat{H}%
_{l,\zeta }\left( \beta ^{-1}\left( \mathbf{k}-\zeta \mathbf{k}_{\ast
l}\right) \right) \mathbf{g}_{n_{l},\zeta }\left( \mathbf{k}\right) \ \text{
for }\left\vert \mathbf{k}-\mathbf{k}_{\ast l}\right\vert \leq \beta
^{1-\epsilon },\;\zeta =\pm ,
\end{equation*}%
$\hat{D}_{l}$ satisfies (\ref{sourloc}), and every function $\hat{H}%
_{l,\zeta }\left( \mathbf{\eta }\right) $, which may depend on $\beta $, is
defined for all $\mathbf{\eta }$ \ and is bounded in $L^{1,a}$ with $a>\frac{%
s}{\epsilon }$ uniformly in $\beta $. Then one can write a nonlinear system
of differential equations for $2N$ scalar envelope functions $z_{l,\zeta
}\left( \tau ,\mathbf{r}\right) $ \ with the initial data $H_{l,\zeta }$, a
linear part of the system has order $\mu \leq 3$ \ and the nonlinearity is
weakly universal as in (\ref{canphys}) and has order $\nu \leq 1$. Let $\hat{%
z}_{l,\zeta }\left( \tau ,\mathbf{k}\right) $, $l=1,...,N$, be the Fourier
transform of a solution to this system Then there exist $\beta _{0}>0$ and a
constant $C$ which does not depend on $\beta ,\varrho $ such that for $\beta
\leq \beta _{0}$ the solution $\mathbf{\hat{u}}$ of (\ref{ubaseq}) with
initial data $\mathbf{\hat{h}}$ can be approximated as follows: 
\begin{equation}
\dsum\limits_{l=1}^{N}\left\Vert \mathbf{\hat{u}}_{l}\left( \tau ,\beta
\right) -\beta ^{-d}\hat{z}_{l,\zeta }\left( \tau ,\beta ^{-1}\left( \mathbf{%
\cdot }-\mathbf{k}_{\ast l}\right) \right) \mathbf{g}_{n_{l},\zeta
}\right\Vert _{E}\leq C\left[ \varrho +\frac{\beta ^{\left( \mu +1\right)
\left( 1-\epsilon \right) }}{\varrho }+\beta ^{\left( \nu +1\right) \left(
1-\epsilon \right) }+\beta ^{s}\right] .  \label{wpint2}
\end{equation}
\end{theorem}

The above-mentioned system with a weakly universal nonlinearity is
constructed based on the equation (\ref{difeqintr}) and $nk$-spectrum $S$
with the help of time averaging (\ref{Gav}) described below. Note that in
the simplest case when $\mu =2$, $\nu =0$, $N=1$ (and $J$ is arbitrary) the
resulting system with a universal nonlinearity is equivalent to classical
Nonlinear Schrodinger equation (NLS). If $N=2$ and $\mathbf{k}_{\ast 1}=-%
\mathbf{k}_{\ast 2}$ we obtain well-known coupled modes system for
counterpropagating waves. This theorem applied to particular systems implies
approximation theorems similar to results of (i) \cite{KalyakinUMN}, \cite%
{SU}, \cite{BF6}, \cite{GiaMielke} on NLS approximation; (ii) \cite{BF6}, 
\cite{GoodmanWH01}, \cite{PW}, \cite{SU01} on coupled mode approximation;
(iii) \cite{SchneiderWayne03} on three-wave approximations. Note also that (%
\ref{wpint2}) implies that if $\varrho =\beta ^{\varkappa ^{\prime }}$ with $%
1<\varkappa ^{\prime }<2$, then the both the first order hyperbolic
equations ($\mu =1$, $\nu =0$) and the second-order NLS ($\mu =2$, $\nu =0$)
provide an approximation for a solution $\mathbf{\hat{u}}$ of (\ref{ubaseq}%
), but NLS provides a better approximation $O\left( \beta ^{\left(
1-\epsilon \right) }\right) $ compared with $O\left( \beta ^{2\left(
1-\epsilon \right) -\varkappa ^{\prime }}\right) $ for first order
hyperbolic equations.

Observe that in the form (\ref{wpn1}) for a simple wavepacket we require $%
\mathbf{g}_{n,\pm }\left( \mathbf{k}_{\ast }\right) $ to be an eigenvector
of the Hermitian matrix $\mathbf{L}\left( \mathbf{k}_{\ast }\right) $, and
one can wonder if $\mathbf{g}_{n,\pm }\left( \mathbf{k}_{\ast }\right) $ can
be replaced with an arbitrary pair of vectors $\mathbf{g}_{\pm }$ in the
case $J>1$. The answer is affirmative, since one can always expand any $%
\mathbf{g}$ with respect to the basis $\mathbf{g}_{n,\pm }\left( \mathbf{k}%
\right) $ using $\Pi _{n,\pm }\left( \mathbf{k}\right) $, but the result
will be a multi-wavepacket with up to $2J$ components rather than a single
wavepacket.

The rest of the paper is organized as follows. In the next section we
illustrate important points of parameter dependence and wavepacket
preservation based on examples. In Section 3 we formulate conditions of
wavepacket preservation including the key resonance invariance condition. In
sections 4 we provide examples of different forms of equations and systems
which involve small or large parameters and can be written in the form of (%
\ref{difeqintr}) after a rescaling. In section 5 we introduce and discuss
integrated modal forms of the evolution equation. In Section 6 \ we
introduce and study the wavepacket interaction system in its relation to the
original system. In Section 7 we approximate the wavepacket interaction
system\ by a certain minimal wavepacket interaction system \ which in
simplest cases turns into the NLS or the coupled modes system.

\section{Preliminary discussion \ and examples}

Observe that the multi-wavepacket preservation as described in Theorems \ref%
{Theorem invarwave}-\ref{Theorem sumcanwave} states in different forms that
(i) its modal composition is essentially preserved; (ii) its $nk$-spectrum
(the set of $nk$-pairs $\left\{ \mathbf{k}_{\ast l},n_{l}\right\} $) remains
the same at all times; (iii) no new modes are excited with a good accuracy
as a result of the nonlinear evolution. The preservation of
multi-wavepackets as they evolve shows also that only the nonlinear
interactions between small neighborhoods of points $\left( \mathbf{k}_{\ast
l},n_{l}\right) \mathbf{\ }$are essential and contribute constructively to
the nonlinear dynamics, whereas the amplitudes of modes with wavevectors $%
\mathbf{k}$ outside those neighborhoods is vanishingly small as $\beta
,\varrho \rightarrow 0$. The later is quite remarkable since the coupling
term $\mathbf{\hat{F}}\left( \mathbf{\hat{U}}\right) \left( \mathbf{k}%
\right) $ in (\ref{difeqfou}) for such $\mathbf{k}$ is not small. A
qualitative explanation to that, confirmed by rigorous analysis, is based on
a fact that the contribution of this term to the solution is a time integral
involving highly oscillatory functions that becomes vanishingly small as $%
\beta ,\varrho \rightarrow 0$. This mechanism is similar to the classical
averaging mechanism for systems of ordinary differential equations
described, for instance, in \cite{BM}; the relevance of the averaging
mechanism for long-wave asymptotics for hyperbolic systems of PDE is
well-known, see \cite{KalyakinUMN}.

We would like to relate now the multi-wavepacket preservation property to
the linear superposition for wavepackets established in \cite{BF7}.
According to that principle if the initial state $\mathbf{h}=\dsum \mathbf{h}%
_{l}$, with $\mathbf{h}_{l}$, $l=1,\ldots ,N$ being "generic" wavepackets,
then the solution $\mathbf{\hat{u}}\left( \tau \right) =\mathcal{G}\left( 
\mathbf{h}\right) \left( \tau \right) $ to the evolution equation (\ref{dfsF}%
) equals with high accuracy to the sum of individual solutions $\mathbf{u}%
_{l}$ of $N$ equations with respective initial data $\mathbf{h}_{l}$.
Namely, if $\beta ,\varrho >0$ satisfy the following relation%
\begin{equation}
\beta ,\varrho \rightarrow 0,\ \beta \geq C_{1}\varrho \text{ with some }%
C_{1}>0,  \label{rbb1}
\end{equation}%
then for all times\emph{\ }$0\leq \tau \leq \tau _{\ast }$ we have%
\begin{gather}
\mathcal{G}\left( \dsum\nolimits_{l=1}^{N}\mathbf{w}_{l}\right) \left( \tau
\right) =\dsum\nolimits_{l=1}^{N}\mathcal{G}\left( \mathbf{w}_{l}\right)
\left( \tau \right) +\mathbf{D}\left( \tau \right) ,  \label{Gsum} \\
\left\Vert \mathbf{D}\left( \tau \right) \right\Vert _{E}=\sup\limits_{0\leq
\tau \leq \tau _{\ast }}\left\Vert \mathbf{D}\left( \tau \right) \right\Vert
_{L^{\infty }}\leq C_{\epsilon }\frac{\varrho }{\beta ^{1+\epsilon }}+C\beta 
\text{ for any }\epsilon >0.  \label{Dbetkap}
\end{gather}%
The \emph{linear superposition principle} is formulated in \cite{BF7} for $%
\beta =C_{2}\varrho ^{1/2}$, but, in fact, the provided proofs of (\ref{Gsum}%
), (\ref{Dbetkap}) remain valid as long as (\ref{rbb1}) holds. Obviously,
the bound $\beta \geq C_{1}\varrho $ in (\ref{rbb1}) determines when (\ref%
{Dbetkap}) becomes trivial. This bound is sharp and examples below show that
when $\beta \sim \varrho $ the remainder $\mathbf{D}\left( \tau \right) $ in
(\ref{Gsum}) does not tend to zero when $\beta \rightarrow 0$.

The both the multi-wavepacket preservation and\ the linear superposition
apply to sums of generic wavepackets. It is important no notice though that
the multi-wavepacket preservation holds for any dependence between $\varrho $
and $\beta $ which satisfy (\ref{rbr}), that is $\varrho \left( \beta
\right) \leq C\beta ^{q}$ with arbitrary small $q$ \ whereas the linear
superposition holds if $\varrho \left( \beta \right) \leq C\beta $. Thus,
the bounds (\ref{rbb1}) on $\beta $\ determine the range of its values for
which the both multi-wavepacket preservation and the linear superposition
hold simultaneously (provided some genericity conditions are satisfied). \
In this range wavepacket preservation provides additional information on
behavior of solutions with single wavepacket initial data, namely that the
solution remains a single wavepacket. Obviously, linear superposition
principle does not follow from multi-wavepacket invariance. \ Below we use
simple examples and models to discuss different ranges of parameters $%
\varrho $ and $\beta $ where wavepacket preservation is valid but the
solutions of equations exhibit different behavior.

\subsection{An exactly solvable model \ and the effect of large group
velocity}

Here we introduce a simple exactly solvable model for our general system (%
\ref{difeqintr}) which makes explicit that in the limit $\varrho \rightarrow
0$ nonlinear effects do not vanish, in particular the blow-up time does not
tend to infinity. This example also shows that on the time scale where $\tau 
$ is of order $1$ solutions undergo significant nonlinear evolution. The
influence of $\varrho $ on solutions through the group velocity in this
example can be seen explicitly. The model is the following system of two
coupled nonlinear first order hyperbolic equations for variables $%
u_{1}\left( x,\tau \right) $, $u_{2}\left( x,\tau \right) $ with
one-dimensional spatial variable $x$: 
\begin{gather}
\partial _{\tau }u_{1}=-\frac{c_{1}}{\varrho }\partial _{x}u_{1}+F_{1}\left(
u_{1},u_{2}\right) ,  \label{toy1} \\
\partial _{\tau }u_{2}=-\frac{c_{2}}{\varrho }\partial _{x}u_{2}+F_{2}\left(
u_{1},u_{2}\right) ,\ c_{1}\neq c_{2},\ \left. u_{1}\right\vert _{\tau
=0}=h_{1}\left( x\right) ,\ \left. u_{2}\right\vert _{\tau =0}=h_{2}\left(
x\right) ,  \label{toy2}
\end{gather}%
where the initial data $h_{1},h_{2}$ in (\ref{toy2}) are of wavepacket form: 
\begin{equation}
h_{1}\left( x\right) =\Phi _{1}\left( \beta x\right) \cos k_{1\ast }x,\
h_{2}\left( x\right) =\Phi _{2}\left( \beta x\right) \cos k_{2\ast }x,\
\left\vert k_{1\ast }\right\vert \neq \left\vert k_{2\ast }\right\vert .
\label{toyinit}
\end{equation}%
We take the nonlinearity to be quadratic and of the following simple form%
\begin{equation}
F_{1}\left( u_{1},u_{2}\right) =u_{1}^{2}+a_{1}u_{1}u_{2},\ F_{2}\left(
u_{1},u_{2}\right) =u_{2}^{2}+a_{2}u_{1}u_{2}.  \label{Fquad}
\end{equation}%
The  system (\ref{toy1})-(\ref{toyinit}) allows for explicit form of
solutions with one-wavepacket initial data, describing a wave propagating
with a constant speed controlled by the linear part and with a shape
evolution controlled by the nonlinearity. This simplest case is compared
then with the case of two-wavepacket initial data, for which explicit
solution is not available.

In the case when $h_{2}=0$ \ the second equation has trivial solution $%
u_{2}=0$ and the system (\ref{toy1})-(\ref{toy2}) reduces to a single
equation (\ref{toy1}). \ The solution to this equation has the form of a
traveling wave $v_{1}\left( x-\frac{c_{1}}{\varrho }\tau ,\tau \right) $
where $v_{1}\left( y,\tau \right) $ is a solution of the ordinary
differential equation 
\begin{equation}
\partial _{\tau }v_{1}=F_{1}\left( v_{1},0\right) ,\ v_{1}\left( y,0\right)
=h_{1}\left( y\right) .  \label{od1}
\end{equation}%
The explicit formula in the case (\ref{od1}) yields 
\begin{equation}
v_{1}\left( x,\tau \right) =\frac{h_{1}\left( x-\frac{c_{1}\tau }{\varrho }%
\right) }{1-\tau h_{1}\left( x-\frac{c_{1}\tau }{\varrho }\right) }=\frac{%
\Phi _{1}\left( \beta \left( x-\frac{c_{1}\tau }{\varrho }\right) \right)
\cos k_{1\ast }\beta \left( x-\frac{c_{1}\tau }{\varrho }\right) }{1-\tau
\Phi _{1}\left( \beta \left( x-\frac{c_{1}\tau }{\varrho }\right) \right)
\cos k_{1\ast }\beta \left( x-\frac{c_{1}\tau }{\varrho }\right) }
\label{toysol}
\end{equation}%
for a time interval $0\leq \tau <\tau _{0}$ where $\tau _{0}=\frac{1}{%
\sup_{y}\left\vert h_{1}\left( y\right) \right\vert }$ is the \emph{blow-up
time}. Obviously, the blow-up time does not depend on $\varrho $.
Consequently, the wave propagates with the velocity $\frac{c_{1}}{\varrho }$
with its shape evolution being controlled by the nonlinearity. \ Similarly,
when $h_{1}=0$ the first equation has the trivial solution $u_{1}=0$ and the
system (\ref{toy1})-(\ref{toy2}) reduces to a\ single equation (\ref{toy2})
which has a solution in the form of a traveling wave $v_{2}\left( x-\frac{%
c_{2}}{\varrho }\tau ,\tau \right) $ propagating with the velocity $\frac{%
c_{2}}{\varrho }$. Observe that for the simple model (\ref{toy1})-(\ref%
{toyinit}) the group velocity coincides with the velocity of a traveling
wave.

The above model is not exactly solvable if the both initial conditions $h_{1}
$ and $h_{2}$ do not vanish. But one can still see the way $\varrho $
influences the nonlinear dynamics quite explicitly by applying the
superposition principle from \cite{BF6}. Indeed, let us assume that $h_{1}$
and $h_{2}$ are two nonzero initial wavepackets. Then the approximate
superposition principle is applicable (in order to put the system in the
framework of \cite{BF6} we use 4-component extension (\ref{toyext}) and set $%
\varrho =\beta ^{\varkappa ^{\prime }}$, $\varkappa ^{\prime }>1$).
According to the principle the exact solution $\left( u_{1},u_{2}\right) $
is approximated by $\left( v_{1}\left( x-\frac{c_{1}}{\varrho }\tau ,\tau
\right) ,v_{2}\left( x-\frac{c_{2}}{\varrho }\tau ,\tau \right) \right) $,
which is explicitly given by (\ref{toysol}) \ with the accuracy $O\left( 
\frac{\varrho }{\beta ^{1+\epsilon }}\right) =O\left( \beta ^{\varkappa
^{\prime }-1-\epsilon }\right) $ with arbitrary small $\epsilon $ if $%
c_{1}\neq c_{2}$. As it as shown in \cite{BF6} the validity of such
approximate presentation is due to the large difference $\frac{c_{1}-c_{2}}{%
\varrho }$ of the group velocities of two wavepackets.

\subsection{Dispersive effects and nonlinearity}

Based on an elementary example of the Nonlinear Schrodinger equation (NLS)%
\begin{equation}
\partial _{\tau }u=-\frac{\mathrm{i}}{\varrho }\left[ \gamma _{0}u+\mathrm{i}%
\gamma _{1}\partial _{x}u+\gamma _{2}\partial _{x}^{2}u\right]
+b_{1}\left\vert u\right\vert ^{2}u,\ u=u\left( x,\tau \right) ,\ x\in 
\mathbb{R}  \label{NLS0}
\end{equation}%
with the initial data in the form of a wavepacket $\left. u\right\vert
_{\tau =0}=\Phi \left( \beta x\right) \mathrm{e}^{\mathrm{i}k_{\ast }x}$ we
would like to explain here why we are interested mostly in the case 
\begin{equation}
\frac{\varrho }{\beta ^{2}}\geq C>0,  \label{wdisp}
\end{equation}%
when the dispersion \emph{is not dominant}. To make the dependence of $u$ on 
$\beta $ and $\varrho $ explicit we change the variables%
\begin{equation}
u\left( x\right) =v\left( \beta x\right) \mathrm{e}^{\mathrm{i}k_{\ast }x},\
\beta x=z,  \label{uve}
\end{equation}%
and obtain equation%
\begin{equation}
\partial _{\tau }v=-\frac{\mathrm{i}}{\varrho }\left[ \gamma _{0}^{\prime }v+%
\mathrm{i}\beta \gamma _{1}^{\prime }\partial _{z}v+\gamma _{2}\beta
^{2}\partial _{z}^{2}v_{1}\right] +b\left\vert v\right\vert ^{2}v,\ \left.
v\right\vert _{\tau =0}=\Phi \left( z\right) ,  \label{gev1}
\end{equation}%
where $\gamma _{1}^{\prime }=\gamma _{1}/\beta +2\gamma _{2}k_{\ast }.$
Changing variables once more 
\begin{equation}
v\left( z,\tau \right) =\mathrm{e}^{-\frac{\mathrm{i}\tau }{\varrho }\gamma
_{0}^{\prime }}w\left( z+\frac{\beta }{\varrho }\gamma _{1}^{\prime }\tau
,\tau \right) ,\ z+\frac{\beta }{\varrho }\gamma _{1}^{\prime }\tau =y,
\label{moving}
\end{equation}%
we obtain for the envelope $w$ the following standard NLS equation 
\begin{equation}
\partial _{\tau }w=-\frac{\mathrm{i}\beta ^{2}}{\varrho }\gamma _{2}\partial
_{y}^{2}w+b\left\vert w\right\vert ^{2}w,\ \left. w\right\vert _{\tau
=0}=\Phi \left( y\right) ,\ 0\leq \tau \leq \tau _{\ast },  \label{stNLS}
\end{equation}%
with initial data independent of the parameters $\beta ,\varrho $. The
behavior of the solution $w$ to the equation (\ref{stNLS}) on the time
interval $0\leq \tau \leq \tau _{\ast }$ is determined by the dispersion
parameter $\frac{\beta ^{2}}{\varrho }$, and evidently linear dispersive
effects become significant when $\frac{\varrho }{\beta ^{2}}\ $is not too
large. If $\frac{\beta ^{2}}{\varrho }\rightarrow \infty $ and $\beta
\rightarrow 0$ \emph{the solution tends to zero at every fixed }$\tau =\tau
_{0}>0$. Indeed, if we take $\varrho =\beta ^{\varkappa ^{\prime }}$, $%
\varkappa ^{\prime }>2$, \ and make another change of variables $\tau
=t\beta ^{\varkappa ^{\prime }-2}$, $w=\beta ^{1-\varkappa ^{\prime }/2}W$,
\ the equation (\ref{stNLS}) reduces to the following problem with small
initial data 
\begin{equation}
\partial _{t}W=-\mathrm{i}\gamma _{2}\partial _{y}^{2}W+b\left\vert
W\right\vert ^{2}W,\ \left. W\right\vert _{t=0}=\beta ^{\varkappa ^{\prime
}/2-1}\Phi \left( y\right) .
\end{equation}%
For small enough $\beta $ the solution $W$ to this problem exists for all $t$
and $W\left( t\right) \rightarrow 0$ as $t\rightarrow \infty $ (see \cite%
{Caz}). In particular, for $t=\tau _{0}\beta ^{2-\varkappa ^{\prime }}$ we
have $w\left( \tau _{0}\right) \rightarrow 0$ when $\beta \rightarrow 0$. \ 

In the general case, the solution dependence on small $\beta ,\varrho $ is
as follows. The dependence on large $\frac{1}{\varrho }$ in (\ref{NLS0}) is
completely described by the change of variables (\ref{moving}), yielding a
wave which (i) moves as a whole with a large group velocity $\frac{-\gamma
_{1}^{\prime }}{\varrho }$; (ii) has a slowly evolving shape as described by 
$v$ and $w$ in (\ref{uve}), (\ref{moving}), (\ref{stNLS}).

The above observations show that for small $\frac{\varrho }{\beta ^{2}}$ the
dispersive effects dominate and control the nonlinear ones. Keeping that in
mind and being interested in stronger nonlinear effects we focus primarily
on the case (\ref{wdisp}), i.e. $\frac{\varrho }{\beta ^{2}}\geq C>0$, for
which there are two scenarios of the nonlinear evolution. In the first
scenario, when $\frac{\beta ^{2}}{\varrho }\rightarrow 0$, the linear
dispersion produces only a small correction to the solution of the equation $%
\partial _{\tau }w=b\left\vert w\right\vert ^{2}w$ with that nonlinear
equation governing the nonlinear dynamics of the envelope $w$ for $\tau
_{\ast }$ being smaller than the blow-up time. In the second scenario, when $%
\beta ^{2}\sim \varrho $, the equation (\ref{stNLS}) becomes independent of $%
\beta ,\varrho $ and describes the evolution of the envelope $w$ governed by
an interplay between the dispersion and the nonlinearity. The case $\beta
^{2}\sim \varrho $ can be also characterized as one where dispersive effects
do occur but they don't dominate nonlinear effects, and, as it is well
known, the dispersion can exactly balance the nonlinearity yielding solitons.

\subsection{A coupled modes system}

Here we illustrate statements of the general theorem on the wavepacket
preservation and the approximate superposition principle by a simple but
still nontrivial example. Let us consider a system of two coupled NLS type
equations for variables $u_{1}\left( x,\tau \right) $, $u_{2}\left( x,\tau
\right) $ with one-dimensional spatial variable $x$%
\begin{gather}
\partial _{\tau }u_{1}=-\frac{\mathrm{i}}{\varrho }\left[ \gamma _{01}+%
\mathrm{i}\gamma _{11}\partial _{x}+\gamma _{21}\partial _{x}^{2}\right]
u_{1}+\left( b_{11}\left\vert u_{1}\right\vert ^{2}+b_{12}\left\vert
u_{2}\right\vert ^{2}\right) u_{1}+c_{12}\left\vert u_{2}\right\vert
^{2}u_{2},  \label{comode1} \\
\partial _{\tau }u_{2}=-\frac{\mathrm{i}}{\varrho }\left[ \gamma _{02}+%
\mathrm{i}\gamma _{12}\partial _{x}+\gamma _{22}\partial _{x}^{2}\right]
u_{2}+\left( b_{21}\left\vert u_{1}\right\vert ^{2}+b_{22}\left\vert
u_{2}\right\vert ^{2}\right) u_{2}+c_{22}\left\vert u_{1}\right\vert
^{2}u_{1},  \label{comode2} \\
\left. u_{1}\right\vert _{\tau =0}=h_{1}\left( x\right) =\Phi _{1}\left(
\beta x\right) \mathrm{e}^{\mathrm{i}k_{\ast 1}x},\ \left. u_{2}\right\vert
_{\tau =0}=h_{2}\left( x\right) =\Phi _{2}\left( \beta x\right) \mathrm{e}^{%
\mathrm{i}k_{\ast 2}x},  \label{comode3}
\end{gather}%
where $\gamma _{ij}$ are real and $b_{ij}$ are complex coefficients and the
initial data in (\ref{comode3}) are in the form of wavepackets with $\Phi
_{j}\left( y\right) $ being Schwartz functions. Notice that if in the
coupled modes system (\ref{comode1})-(\ref{comode3}) $h_{2}=0$ and $%
c_{12}=c_{22}=0$ then it has trivial solution $u_{2}=0$, and reduces to a
single NLS equation of the form (\ref{NLS0}). The dependence of the solution 
$\left\{ u_{1},u_{2}\right\} $ on the large $\frac{1}{\varrho }$ is captured
by the change of variables (\ref{moving}). Namely, $u_{1}$ is a wave with a
slowly varying envelope described by $v_{1}$ which moves with large velocity 
$\frac{-\gamma _{11}^{\prime }}{\varrho }$. The dependence on $\beta $ is of
the form $v_{1}\left( y,\tau \right) =w_{1}\left( \beta y,\tau \right) $
(see following subsection for details).\ Similarly we can consider the case
when $h_{1}=0$ for which the first equation has\ trivial solution $u_{1}=0$,
so the system (\ref{comode1})-(\ref{comode2}) reduces to a single equation (%
\ref{comode2}) with the solution represented by a wave having large spacial
extension proportional to $\frac{1}{\beta }$ and moving with the large
velocity $\frac{-\gamma _{12}^{\prime }}{\varrho }$.

\subsubsection{The superposition principle}

Let us assume here that $h_{1}\neq 0$, $h_{2}\neq 0$, $c_{12}\neq 0$, $%
c_{22}\neq 0$ and $\beta =\varrho ^{\varkappa }$, \ $0<\varkappa <1$.
Applying the superposition principle we obtain for generic $k_{\ast 1}$, $%
k_{\ast 2}$ the following representation of the exact solution 
\begin{equation*}
u_{1}\left( x,\tau \right) =v_{1}\left( x,\tau \right) \mathrm{e}^{\mathrm{i}%
k_{\ast 1}x}+D_{1},\ u_{2}\left( x,\tau \right) =v_{2}\left( x,\tau \right) 
\mathrm{e}^{\mathrm{i}k_{\ast 2}x}+D_{2}
\end{equation*}%
where $v_{1}\left( x,\tau \right) $ is a solution of the NLS equation (\ref%
{comode1}) with $b_{12}=c_{12}=0$, with $v_{2}\left( x,\tau \right) $ being
a solution to a similar decoupled NLS equation for $b_{22}=c_{22}=0$, and $%
D_{1}$ and $D_{2}$ are small terms satisfying 
\begin{equation}
\sup\nolimits_{0\leq \tau \leq \tau _{\ast }}\left\Vert D_{1}\left( \cdot
,\tau \right) \right\Vert _{L^{\infty }}+\sup\nolimits_{0\leq \tau \leq \tau
_{\ast }}\left\Vert D_{2}\left( \cdot ,\tau \right) \right\Vert _{L^{\infty
}}\leq C\beta ^{\varkappa ^{\prime }-1-\epsilon }+C\beta ,\ \varkappa
^{\prime }=\varkappa ^{-1}.  \label{exDest}
\end{equation}%
\emph{We would like to emphasize here that the coupling terms }$%
b_{12}\left\vert u_{2}\right\vert ^{2}u_{1}+c_{12}\left\vert
u_{2}\right\vert ^{2}u_{2}$\emph{\ and }$b_{21}\left\vert u_{1}\right\vert
^{2}u_{2}+c_{22}\left\vert u_{2}\right\vert ^{2}u_{2}$\emph{\ in the
equations (\ref{comode1})-(\ref{comode2}) are not small whereas their
ultimate contributions to the solutions are small}. One can
explain/interpret that phenomenon as being due to the destructive wave
interference and mismatch of group velocities.

\subsubsection{Wavepacket preservation}

Here we assume that $h_{1}\neq 0$, $h_{2}=0$, $c_{12}\neq 0$, $c_{22}\neq 0$
and $\varrho =\beta ^{\varkappa ^{\prime }}$,\ $\ 0<\varkappa ^{\prime }\leq
2$. According to the wavepacket preservation we have%
\begin{equation*}
u_{1}\left( x,\tau \right) =v_{1}\left( x,\tau \right) \mathrm{e}^{\mathrm{i}%
k_{\ast 1}x}+D_{1},\ u_{2}\left( x,\tau \right) =D_{1},
\end{equation*}%
where $v_{1}\left( x,\tau \right) $ is a solution of (\ref{comode1}) with $%
b_{12}=0$ ,\ $c_{12}=0$, and $D_{1}$ and $D_{2}$ are small terms satisfying 
\begin{equation*}
\sup\nolimits_{0\leq \tau \leq \tau _{\ast }}\left\Vert D_{1}\left( \cdot
,\tau \right) \right\Vert _{L^{\infty }}+\sup\nolimits_{0\leq \tau \leq \tau
_{\ast }}\left\Vert D_{2}\left( \cdot ,\tau \right) \right\Vert _{L^{\infty
}}\leq C\varrho
\end{equation*}%
Notice once more (see the above section) an interesting phenomenon: the
equation (\ref{comode2}) for $u_{2}\left( x,\tau \right) $ has a coupling
term $b_{21}\left\vert u_{1}\right\vert ^{2}u_{2}+c_{22}\left\vert
u_{1}\right\vert ^{2}u_{1}$ which does not become small as $\beta ,\varrho
\rightarrow 0$, but, remarkably, its ultimate contribution to the solution
is small.

\subsubsection{Limitations of the superposition principle}

Now we provide an example based on the system (\ref{comode1})-(\ref{comode3}%
) with $c_{12}=c_{22}=0$ showing that the above estimate (\ref{exDest}) in
the superposition principle is sharp in the sense that $\beta ^{\varkappa
^{\prime }-1-\epsilon }$ cannot be replaced by $\beta ^{\varkappa ^{\prime
}-1+\epsilon }$ with $\varkappa ^{\prime }\geq 1$. We set here $\varkappa
^{\prime }=1$ and $\varrho =\beta $. After the change of variables (\ref{uve}%
) for $u_{1},u_{2}$ followed by yet another change of variables $\beta x=z$, 
$v_{1}=\mathrm{e}^{-\mathrm{i}\tau \frac{\gamma _{01}^{\prime }}{\beta }%
}w_{1}$, $v_{2}=\mathrm{e}^{-\mathrm{i}\tau \frac{\gamma _{01}^{\prime }}{%
\beta }}w_{2}$ we obtain from (\ref{comode1})-(\ref{comode3}) the following
system: 
\begin{gather*}
\partial _{\tau }w_{1}=-\mathrm{i}\left[ \mathrm{i}\gamma _{11}^{\prime
}\partial _{z}w_{1}+\beta \gamma _{21}\partial _{z}^{2}w_{1}\right] +\left(
b_{11}\left\vert w_{1}\right\vert ^{2}+b_{12}\left\vert w_{2}\right\vert
^{2}\right) w_{1}, \\
\partial _{\tau }w_{2}=-\mathrm{i}\left[ \mathrm{i}\gamma _{12}^{\prime
}\partial _{z}w_{2}+\beta \gamma _{22}\partial _{z}^{2}w_{2}\right] +\left(
b_{21}\left\vert w_{1}\right\vert ^{2}+b_{22}\left\vert w_{2}\right\vert
^{2}\right) w_{2}, \\
\left. w_{1}\right\vert _{\tau =0}=\Phi _{1}\left( z\right) ,\left.
w_{2}\right\vert _{\tau =0}=\Phi _{2}\left( z\right) .
\end{gather*}%
This system has a regular dependence on $\beta $ as $\beta \rightarrow 0$
with the solution converging in $L^{\infty }$ to the solution of the system
with $\beta =0$. If we set now in the last system $b_{12}=b_{21}=0$ it turns
into a system of two decoupled equations. Notice then that the difference
between the solutions of the decoupled system and the original one does not
tend to zero as $\beta \rightarrow 0$, implying that the superposition
principle does not hold when $\varrho =\beta $.

\subsection{Wavepacket interaction system with a universal nonlinearity}

We will prove in the following sections that the dynamics of a
multi-wavepacket with a universally resonance invariant $nk$-spectrum for a
general system can be approximated with the accuracy $O\left( \varrho
\right) $ by substituting the nonlinearity with a properly constructed \emph{%
universal} or weakly universal one. Here we provide an example of a system,
called \emph{wavepacket interaction system}, with a universal nonlinearity
and show that its dynamics preserves simple wavepackets as in (\ref{wpint}).
It is shown later that universal nonlinearities are related to universally
invariant multi-wavepackets in the sense of Definition \ref{Definition
omclos}.

Wavepacket interaction system with universal nonlinearity has the form
similar to NLS, namely 
\begin{gather}
\partial _{\tau }u_{j,\zeta }=\frac{1}{\varrho }\left[ -\mathrm{i}\zeta
\gamma _{0,j}+\gamma _{1,j}\cdot \nabla _{\mathbf{r}}u_{j,\zeta }-\mathrm{i}%
\zeta \nabla _{\mathbf{r}}\cdot \gamma _{2,j}\nabla _{\mathbf{r}}u_{j,\zeta }%
\right] +F_{j,\zeta }\left( \vec{u}\right) ,\ \mathbf{r}\in \mathbb{R}^{d},
\label{cansystext} \\
\vec{u}=\left( u_{1+},u_{1-},\ldots ,u_{N+},u_{N-}\right) ,\ \ j=1,\ldots
,N,\ \zeta =\pm ,  \label{uar} \\
\left. u_{j,\zeta }\right\vert _{\tau =0}=h_{j,\zeta },h_{j,\zeta }\left( 
\mathbf{r}\right) =\Phi _{j}\left( \beta \mathbf{r}\right) \mathrm{e}^{%
\mathrm{i}\zeta \mathbf{k}_{\ast j}\cdot \mathbf{r}},
\end{gather}%
where for every $j$ coefficient$\ \gamma _{0,j}\in \mathbb{R},$ $\gamma
_{1,j}\in \mathbb{R}^{d}$ is a vector, $\gamma _{2,j}$ \ is a symmetric $%
d\times d$ matrix, $\gamma _{1,j}\cdot \nabla _{\mathbf{r}}$ is a first
order scalar differential operator, $\nabla _{\mathbf{r}}\cdot \gamma
_{2,j}\nabla _{\mathbf{r}}$ is the second order scalar differential
operator,and the \emph{universal polynomial nonlinearities} $F_{j,\zeta }$
have the following form: 
\begin{gather}
F_{j,\zeta }\left( \vec{u}\right) =\sum\nolimits_{\nu =1}^{\nu
_{F}}\sum\nolimits_{\left\vert \vec{\nu}\right\vert =\nu }b_{\vec{\nu}%
,j,\zeta }\dprod\nolimits_{l=1}^{N}\left( u_{l,+}u_{l,-}\right) ^{\nu
_{l}}u_{j,\zeta },  \label{Ppm} \\
\text{where }\vec{\nu}=\left( \nu _{1},\ldots ,\nu _{N}\right) ,\ \
j=1,\ldots ,N,\;\zeta =\pm .  \notag
\end{gather}

\begin{remark}
Notice that if we set $h_{j,-}=h_{j,+}^{\ast }$, $b_{\vec{\nu},j,+}=b_{\vec{%
\nu},j,-}^{\ast }=b_{\vec{\nu},j}$ and $u_{j,+}=u_{j,-}^{\ast }=u_{j}$ then $%
u_{l,+}u_{l,-}=\left\vert u_{l,+}\right\vert ^{2}$ and $F_{j,+}\left( \vec{u}%
\right) $ turns into 
\begin{equation}
F_{j}\left( u_{1},\ldots ,u_{N}\right) =\sum\nolimits_{\nu =1}^{\nu
_{F}}\sum\nolimits_{\left\vert \vec{\nu}\right\vert =\nu }b_{\vec{\nu}%
,j}\dprod\nolimits_{l=1}^{N}\left\vert u_{l}\right\vert ^{2\nu _{l}}u_{j},\
\   \label{Fjcan}
\end{equation}%
and equations of (\ref{cansystext}) with $\zeta =+$ \ turn into 
\begin{gather}
\partial _{\tau }u_{j}=\frac{1}{\varrho }\left[ -\mathrm{i}\gamma
_{0,j}+\gamma _{1j}\cdot \nabla _{\mathbf{r}}u_{j}-\mathrm{i}\nabla _{%
\mathbf{r}}\cdot \gamma _{2,j}\nabla _{\mathbf{r}}u_{j}\right] +F_{j}\left(
u_{1},\ldots ,u_{N}\right) ,  \label{cansyst1} \\
\left. u_{j}\right\vert _{\tau =0}=h_{j,+},\ j=1,\ldots ,N,\;\zeta =\pm . 
\notag
\end{gather}%
Obviously, a solution of (\ref{cansyst1}) defines a solution $u_{j,+}=u_{j}$%
, $u_{j,-}=u_{j}^{\ast }$ of the system (\ref{cansystext}). In the simplest
case $N=1$, $d=1$ (\ref{cansyst1}) takes the form of classical NLS: $%
\partial _{\tau }u=\frac{\gamma _{1}}{\varrho }\partial _{x}u-\mathrm{i}%
\frac{\gamma _{2}}{\varrho }\partial _{x}^{2}u+b\left\vert u\right\vert
^{2}u.$
\end{remark}

Note that the \emph{universal nonlinearity} $F_{j,\zeta }$ has a
characteristic property 
\begin{equation}
F_{j,\zeta }\left( \mathrm{e}^{\mathrm{i}\phi _{1}t}u_{1,+},\mathrm{e}^{-%
\mathrm{i}\phi _{1}t}u_{1,-},\ldots ,\mathrm{e}^{\mathrm{i}\phi
_{N}t}u_{N,+},\mathrm{e}^{-\mathrm{i}\phi _{N}t}u_{N,-}\right) =\mathrm{e}^{%
\mathrm{i}\zeta \phi _{j}t}F_{j,\zeta }\left( u_{1+},u_{1-},\ldots
,u_{N+},u_{N-}\right) .  \label{invfi1}
\end{equation}%
holding for arbitrary \ set values $\phi _{i}$. We also consider more
general nonlinearities $F$ for which (\ref{invfi1}) holds for a fixed set of
frequencies $\phi _{l}=\omega _{n_{l}}\left( \mathbf{k}_{\ast l}\right) $,
and call them \emph{weakly universal. }We introduce now the \emph{averaging
operator} $A_{T}$ acting on polynomial functions $F:\left( \mathbb{C}%
^{2}\right) ^{N}\rightarrow \left( \mathbb{C}^{2}\right) ^{N}$ by 
\begin{gather}
\left( A_{T}F\right) _{j,\zeta }=\left( A_{T,\vec{\phi}}F\right) _{j,\zeta }=
\label{cana1} \\
\frac{1}{T}\int_{0}^{T}\mathrm{e}^{-\mathrm{i}\zeta \phi _{j}t}F_{j,\zeta
}\left( \mathrm{e}^{\mathrm{i}\phi _{1}t}u_{1,+},\mathrm{e}^{-\mathrm{i}\phi
_{1}t}u_{1,-},\ldots ,\mathrm{e}^{\mathrm{i}\phi _{N}t}u_{N,+},\mathrm{e}^{-%
\mathrm{i}\phi _{N}t}u_{N,-}\right) \mathrm{d}t,  \notag
\end{gather}%
where $\vec{\phi}=\left( \phi _{1},\ldots ,\phi _{N}\right) $. The operator $%
A_{T,\vec{\phi}}$ depends on the frequency vector $\vec{\phi}=\left( \phi
_{1},\ldots ,\phi _{N}\right) $. If $F\ $is a universal polynomial
nonlinearity, then $\ \left( A_{T,\vec{\phi}}F\right) _{j,\zeta }=F_{j,\zeta
}$ for any choice of frequencies $\phi _{1},\ldots ,\phi _{N}$. Note that
averaging 
\begin{equation}
G_{\text{av},j,\zeta }\left( \vec{u}\right) =\lim_{T\rightarrow \infty
}\left( A_{T,\vec{\phi}}G\right) _{j,\zeta }\left( \vec{u}\right)
\label{Gav}
\end{equation}%
is defined for any polynomial nonlinearity $G:\left( \mathbb{C}^{2}\right)
^{N}\rightarrow \left( \mathbb{C}^{2}\right) ^{N}$. If $\vec{\phi}$ is
generic,\emph{\ }then $G_{\text{av},j,\zeta }\left( \vec{u}\right) $ is
always a universal nonlinearity. In a general case $G_{\text{av},j,\zeta }$
for given frequencies $\vec{\phi}$ one obtains a weakly universal
nonlinearity which might be not universal.

Systems with universal nonlinearities have interesting properties which we
describe in the following Proposition and remark.

\begin{proposition}
\label{simple wave preservation}Let $\varrho =\beta $ and $\gamma _{2,j}=0$.
Then evolution governed by the first order system with a universal
nonlinearity (\ref{cansystext}) preserves simple wavepackets as defined by (%
\ref{wpint}).
\end{proposition}

\begin{proof}
Let $\vec{u}\left( \tau \right) $ be a solution of (\ref{cansystext}) \ for $%
0\leq \tau \leq \tau _{\ast }$. Using the property (\ref{invfi1}) we change
variables 
\begin{equation}
u_{j,\zeta }=\mathrm{e}^{\mathrm{i}\zeta \mathbf{k}_{\ast j}\cdot \mathbf{r}}%
\mathrm{e}^{-\mathrm{i}\frac{\zeta \gamma _{0,j}}{\varrho }\tau }\mathrm{e}%
^{-\mathrm{i}\frac{\gamma _{0j,\zeta }^{\prime }}{\beta }\tau }v_{j,\zeta
},\gamma _{0j,\zeta }^{\prime }=-\zeta \mathbf{\gamma }_{1j}\cdot \mathbf{k}%
_{\ast j}  \label{ujcan}
\end{equation}%
and obtain from (\ref{cansystext}) 
\begin{equation}
\partial _{\tau }v_{j,\zeta }=\frac{1}{\beta }\mathbf{\gamma }_{1j}\cdot
\nabla _{\mathbf{r}}v_{j,\zeta }+F_{j,\zeta }\left( \vec{v}\right) ,\ \left.
v_{j,\zeta }\right\vert _{\tau =0}=\Phi _{j,\zeta }\left( \beta \mathbf{r}%
\right) .  \label{ujcan1}
\end{equation}%
Changing variables 
\begin{equation}
v_{j,\zeta }\left( \mathbf{r},\tau \right) =w_{j,\zeta }\left( \beta \mathbf{%
r},\tau \right) ,\ \beta \mathbf{r}=\mathbf{z},  \label{movingcan}
\end{equation}%
we obtain from (\ref{ujcan1}) that $w_{j}$ is a solution of the following
system of differential equations 
\begin{equation}
\partial _{\tau }w_{j,\zeta }=\gamma _{1j}\cdot \nabla _{\mathbf{z}%
}w_{j,\zeta }+F_{j,\zeta }\left( \vec{w}\right) ,\ \left. w_{j,\zeta
}\right\vert _{\tau =0}=\Phi _{j,\zeta }\left( \mathbf{z}\right) ,
\label{redcan}
\end{equation}%
which does not depend on $\beta $. Then using (\ref{movingcan}) and (\ref%
{ujcan}) we observe that every component $u_{l}$ of the solution to (\ref%
{cansystext}) has the form of a simple wavepacket for every $\tau \in $ $%
\left[ 0,\tau _{\ast }\right] $, with an envelope $\hat{w}_{j}\left( \tau
\right) $.
\end{proof}

\begin{remark}
Equations (\ref{cansystext}) with universal nonlinearities allow special
solutions in the form of$\ \ u_{j,\zeta }=\mathrm{e}^{\mathrm{i}\mathbf{k}%
_{\ast j}\cdot \mathbf{r}}\mathrm{e}^{-\mathrm{i}\frac{\gamma _{0j}^{\prime }%
}{\beta }\tau }v_{j,\zeta }\left( \tau \right) $ where $v_{j,\zeta }\left(
\tau \right) $ do not depend on $\mathbf{r}$. If the initial data in (\ref%
{ujcan1}) are constants, $\Phi _{j,\zeta }\left( \beta \mathbf{r}\right)
=\Phi _{j,\zeta }\left( 0\right) $, then (\ref{ujcan1}) turns into a system
of ODE. This implies that every linear subspace of pure modal functions with
the basis $v_{j}\mathrm{e}^{\mathrm{i}\mathbf{k}_{\ast j}\cdot \mathbf{r}}$, 
$v_{j,-}\mathrm{e}^{-\mathrm{i}\mathbf{k}_{\ast j}\cdot \mathbf{r}},$ $%
j=1,..,N$ is invariant with respect to nonlinear equations (\ref{cansystext}%
). Another class of special solutions of (\ref{cansystext})\ are
time-harmonic solutions of the form $u_{j,\zeta }\left( \mathbf{r},\tau
\right) =\mathrm{e}^{-\mathrm{i}\zeta \omega _{j}\tau }v_{j,\zeta }\left( 
\mathbf{r}\right) $ where $v_{j,\zeta }$ solve a nonlinear eigenvalue
problem; for universal nonlinearities $\omega _{j}$ can be considered as
unknown nonlinear eigenvalue. Existence of such special solutions is a
special property of universal and weakly universal nonlinearities. It is
remarkable that original nonlinear equations might not have time harmonic
solutions whereas equations with canonical nonlinearities which approximate
evolution of wavepackets (see Theorem \ref{Theorem sumcanwave}) admit such
solutions.
\end{remark}

\subsection{ Invariance of excited modes for finite-dimensional ODE's}

Here we discuss the resonance invariance conditions imposed in Theorem \ref%
{Theorem sumwave} in a simpler case of finite-dimensional ODE's. In this
case one can also see the rise of universal nonlinearities in the process of
time averaging. As we already discussed in the introduction, a PDE system (%
\ref{difeqintr}) when restricted to constant functions turns into the
following system of ODE's 
\begin{equation}
\partial _{\tau }\mathbf{U}=-\frac{\mathrm{i}}{\varrho }\mathbf{L}^{0}%
\mathbf{U}+\mathbf{F}\left( \mathbf{U}\right) ,\ \left. \mathbf{U}\left(
\tau \right) \right\vert _{\tau =0}=\mathbf{h},\ \mathbf{h}\ \in \mathbb{C}%
^{2J},\ \mathbf{U}\in \mathbb{C}^{2J},  \label{UODE}
\end{equation}%
where $\mathbf{F}\left( \mathbf{U}\right) $ is a polynomial, $\mathbf{U}%
=\left( U_{1,+},U_{1,-},\ldots ,U_{J,+},U_{J,-}\right) \in \mathbb{C}^{2J}$.
We assume that the eigenvalues $\omega _{n,\zeta }\left( \mathbf{0}\right)
=\omega _{n,\zeta }^{0}$ of the Hermitian matrix $\mathbf{L}^{0}=\left. 
\mathbf{L}\left( \mathbf{k}\right) \right\vert _{\mathbf{k}=0}$ are distinct 
$\omega _{j,+}^{0}\neq \omega _{i,+}^{0}$ for $j\neq i$ and the symmetry
conditions (\ref{invsym}) take the form $\omega _{n,-\zeta }^{0}=-\omega
_{n,\zeta }^{0}$. We also assume that the eigenvectors of $\mathbf{L}^{0}$
coincide with the coordinate orts in $\mathbb{C}^{2J}$. The following limit
case of Theorem \ref{Theorem sumwave} with $\beta =0$ shows that solutions
to this system have the property to preserve the set of initially excited
modes.

\begin{theorem}
\label{Theorem sumwave1} Let the initial data $\mathbf{h}=\left(
h_{1,+},h_{1,-},\ldots ,h_{J,+},h_{J,-}\right) \in \mathbb{C}^{2J}\mathbf{\ }
$in (\ref{UODE}) have non-zero components $h_{j,\zeta }$ only for a\ subset $%
B$ of indices $j\in \left\{ 1,\ldots ,J\right\} $, and let $B^{\prime
}=\left\{ 1,\ldots ,J\right\} \setminus B$ be its complementary set. Assume
that $B$ is \emph{resonance invariant} in the sense that the resonance
equation 
\begin{equation}
\omega _{n^{\prime },\zeta }^{0}-\sum\nolimits_{j=1}^{m}\omega _{n_{j},\zeta
^{\left( j\right) }}^{0}=0,\text{where }n_{j}\in B,\;\zeta ^{\left( j\right)
}\in \left\{ +,-\right\}  \label{res0}
\end{equation}%
does not have solutions if $n^{\prime }\in B^{\prime }$(compare with
Definition \ref{Definition omclos} in the special case when all $\mathbf{k}%
_{\ast l}=0$). Then under the nonlinear evolution of (\ref{UODE}) modes with
indices $n^{\prime }\in B^{\prime }$ remain essentially unexcited in the
following sense%
\begin{equation}
\sup_{0\leq \tau \leq \tau _{\ast }}\left\vert U_{n^{\prime }}\left( \tau
\right) \right\vert \leq C\varrho \text{ for all }n^{\prime }\in B^{\prime }.
\label{Ur}
\end{equation}
\end{theorem}

Note that $\mathbf{F}\left( \mathbf{U}\right) $ provides a nonlinear
coupling between modes $U_{n_{j},\zeta ^{\left( j\right) }}$ with $n_{j}\in
B $ and $U_{n^{\prime },\zeta }$ with $n^{\prime }\in B^{\prime }$, but the
resulting interaction is not $O\left( 1\right) $ on a fixed time interval $%
\left[ 0,\tau _{\ast }\right] $ as one might expect, but rather of order $%
O\left( v\right) $ as (\ref{Ur}) shows. One way to prove Theorem \ref%
{Theorem sumwave1} is to follow the proofs of Theorems \ref{Theorem Dsmall}
and \ref{Theorem uminw} with obvious modifications and simplifications. In
particular, instead of (\ref{dfsF}) one has to consider\ the following
system with oscillatory coefficients 
\begin{equation}
\partial _{\tau }\mathbf{u}=\mathrm{e}^{\frac{\mathrm{i}\tau }{\varrho }%
\mathbf{L}^{0}}\mathbf{F}\left( \mathrm{e}^{\frac{-\mathrm{i}\tau }{\varrho }%
\mathbf{L}^{0}}\mathbf{u}\right) ,\ \left. \mathbf{u}\left( \tau \right)
\right\vert _{\tau =0}=\mathbf{h}.  \label{Uav}
\end{equation}%
Alternatively, Theorem \ref{Theorem sumwave1} can be derived directly from
the classical time averaging principle. Indeed, the time averaging of (\ref%
{Uav}) yields the following averaged system%
\begin{equation*}
\partial _{\tau }\mathbf{v}=\mathbf{F}_{\text{av}}\left( \mathbf{v}\right)
,\ \left. \mathbf{v}\left( \tau \right) \right\vert _{\tau =0}=\mathbf{h},
\end{equation*}%
where $\mathbf{F}_{\text{$\limfunc{av}$}}$ is defined as in \ (\ref{cana1}),
(\ref{Gav}) with the frequencies $\phi _{j}=\omega _{j,+}^{0}$. From the
Krylov-Bogolyubov averaging theorem (see \cite{BM}, \cite{MitrN}) one
obtains 
\begin{equation*}
\left\vert \mathbf{v}\left( \tau \right) -\mathbf{u}\left( \tau \right)
\right\vert \leq C\varrho ,\ 0\leq \tau \leq \tau _{\ast }.
\end{equation*}%
A straightforward examination shows that if $B$ is resonance invariant and $%
j\in B^{\prime }$ then the polynomial components $F_{\text{$\limfunc{av},$}%
j,\zeta }\left( \mathbf{v}\right) $ factorize into $F_{\text{$\limfunc{av},$}%
j,\zeta }\left( \mathbf{v}\right) =\sum_{j^{\prime }\in B^{\prime },\zeta
^{\prime }}F_{\text{$\limfunc{av},$}j^{\prime },\zeta ^{\prime }}^{1}\left( 
\mathbf{v}\right) v_{j^{\prime },\zeta ^{\prime }}$, implying (\ref{Ur})
since $v_{j,\zeta }\left( 0\right) =0$ for $j\in B^{\prime }$.

A stronger \emph{universal resonance invariance} condition in Definition \ref%
{Definition omclos} also takes a simpler form in the ODE case. Indeed, let
us collect the terms in (\ref{res0}) at different $\omega _{j,+}^{0}$ as in (%
\ref{rearr0}), namely 
\begin{equation}
\omega _{n^{\prime },\zeta }^{0}-\sum\nolimits_{j=1}^{m}\omega _{n_{j},\zeta
^{\left( j\right) }}^{0}=\sum\nolimits_{i=1}^{J}\delta _{i}\omega _{i,+}^{0},%
\text{ where }\delta _{i}\text{ are integers,}  \label{bi}
\end{equation}%
Similarly to Definition \ref{Definition omclos} we call $B$ universally
resonance invariant if every solution to the resonance equation (\ref{res0})
must have $n^{\prime }\in B$ and every\ coefficient $\delta _{i}$ in (\ref%
{bi}) for the solution is zero, i.e. $\delta _{i}=0,$ $i=1,\ldots ,J$.
Obviously, if all $\omega _{n,+}^{0}$ are rationally independent then it is
universally resonance invariant.

Now let us look how universal nonlinearities arise under time averaging.
Observe that if the entire set $\left\{ 1,\ldots ,J\right\} $ is universally
resonance invariant and $F_{j,\zeta }\left( \mathbf{v}\right) $ are
arbitrary polynomials, then the polynomials $F_{\text{$\limfunc{av},$}%
j,\zeta }\left( \mathbf{v}\right) $ are obtained by discarding the
"resonant" terms in $\mathrm{e}^{\frac{\mathrm{i}\tau }{\varrho }\mathbf{L}%
^{0}}\mathbf{F}\left( \mathrm{e}^{\frac{-\mathrm{i}\tau }{\varrho }\mathbf{L}%
^{0}}\mathbf{u}\right) $ yielding universal form (\ref{Ppm}), (\ref{Fjcan}).
For example, if $\mathbf{F}$ is an arbitrary cubic nonlinearity in $\mathbb{C%
}^{2N}$ then the time averaging yields NLS-like nonlinearity $\ \mathbf{F}_{%
\text{$\limfunc{av}$}}$ with components 
\begin{equation*}
F_{\text{$\limfunc{av},$}j,\zeta }\left( u_{1,+},u_{1,-},\ldots
,u_{N,+},u_{N,-}\right) =\sum\nolimits_{l=1}^{N}b_{l,j,\zeta
}u_{l,+}u_{l,-}u_{j,\zeta }.
\end{equation*}%
When $B$ is resonance invariant \ but not universally resonance invariant \
the averaging produces a weakly universal nonlinearity. A nonlinearity which
is weakly universal \ but not universal may include additional terms, for
example the cubic nonlinearity in classical four-wave interaction system \
where it is assumed that $\omega _{2,-}^{0}+\omega _{3,+}^{0}+\omega
_{4,+}^{0}=\omega _{1,+}^{0}$ \ (see \cite{Phillips} p. 201) in the equation
for $u_{1,+}$ in addition to NLS-like terms involves the product $%
u_{2,-}u_{3,+}u_{4,+}$ .

\section{Conditions and definitions}

In this section we formulate and discuss definitions and conditions under
which we study the nonlinear evolutionary system (\ref{difeqintr}) through
its modal, Fourier form (\ref{difeqfou}). Most of the conditions and
definitions are naturally formulated for the modal form (\ref{difeqfou}),
and this is one of the reasons we use it as the basic form.

\subsection{Linear part}

The basic properties of the linear part $\mathbf{L}\left( \mathbf{\mathbf{k}}%
\right) $ of the system (\ref{difeqfou}), which is a $2J\times 2J$ Hermitian
matrix with eigenvalues $\omega _{n,\zeta }\left( \mathbf{k}\right) $, has
been already discussed in the Introduction. To account for all needed
properties of $\mathbf{L}\left( \mathbf{\mathbf{k}}\right) $ we define the
singular set of points $\mathbf{\mathbf{k}}$.

\begin{definition}[band-crossing points]
\label{Definition band-crossing point} We call $\mathbf{k}_{0}$ a \emph{%
band-crossing point} for $\mathbf{L}\left( \mathbf{\mathbf{k}}\right) $ if $%
\omega _{n+1,\zeta }\left( \mathbf{k}_{0}\right) =\omega _{n,\zeta }\left( 
\mathbf{k}_{0}\right) $ \ for some $n,\zeta $ or $\mathbf{L}\left( \mathbf{%
\mathbf{k}}\right) $ is not continuous at $\mathbf{k}_{0}$ or if $\omega
_{1,\pm }\left( \mathbf{k}_{0}\right) =0$, we denote the set of such points
by $\sigma _{bc}$.
\end{definition}

In the next Condition we collect all constraints imposed on the linear
operator $\mathbf{L}\left( \mathbf{\mathbf{k}}\right) $.

\begin{condition}[linear part]
\label{clin}The linear part $\mathbf{L}\left( \mathbf{\mathbf{k}}\right) $
of the system (\ref{difeqfou}) is a $2J\times 2J$ Hermitian matrix with
eigenvalues $\omega _{n,\zeta }\left( \mathbf{k}\right) $ and corresponding
eigenvectors $\mathbf{g}_{n,\zeta }\left( \mathbf{k}\right) $ satisfying for 
$\mathbf{k}\notin \sigma _{bc}$ the basic relations (\ref{OmomL})-(\ref%
{invsym}). In addition to that we assume:

\begin{enumerate}
\item[(i)] the set of band-crossing points $\sigma _{bc}$ is a closed,
nowhere dense set in $\mathbb{R}^{d}$ and has zero Lebesgue measure;

\item[(ii)] the entries of the Hermitian matrix $\mathbf{L}\left( \mathbf{%
\mathbf{k}}\right) $ are infinitely differentiable in $\mathbf{k}$ for all $%
\mathbf{k}\notin \sigma _{bc}$ that readily implies via the spectral theory, 
\cite{Kato}, infinite differentiability of all eigenvalues $\omega
_{n}\left( \mathbf{k}\right) $ in $\mathbf{k}$ for all $\mathbf{k}\notin
\sigma $;

\item[(iii)] $\mathbf{L}\left( \mathbf{\mathbf{k}}\right) $ satisfies
polynomial bound 
\begin{equation}
\left\Vert \mathbf{L}\left( \mathbf{\mathbf{k}}\right) \right\Vert \leq
C\left( 1+\left\vert \mathbf{\mathbf{k}}\right\vert ^{p}\right) ,\ \mathbf{k}%
\in \mathbb{R}^{d},\ \text{for some }C>0\text{ and }p>0\text{.}  \label{Lpol}
\end{equation}
\end{enumerate}
\end{condition}

\begin{remark}[dispersion relations symmetry]
\label{Remark symmetry}The symmetry condition (\ref{invsym}) on the
dispersion relations naturally arise in many physical problems, for example
Maxwell equations in periodic media, see \cite{BF1}-\cite{BF3}, \cite{BF5},
or when $\mathbf{L}\left( \mathbf{k}\right) $ originates from a Hamiltonian. 
\emph{We would like to stress that this symmetry conditions are not imposed
to simplify studies }but rather to take into account fundamental symmetries
of physical media. In fact, the opposite case when ( (\ref{invsym}) is
assumed not to hold is much simpler. The symmetry\ creates resonant
nonlinear interactions, which makes studies more intricate. Interestingly,
many problems without symmetries can be put into the framework with symmetry
by an extension of the relevant system (see Section 4).
\end{remark}

\begin{remark}[band-crossing points]
Band-crossing points are discussed in more details in \cite[Section 5.4]{BF1}%
, \cite[Sections 4.1, 4.2]{BF2}. In particular, generically the set $\sigma
_{bc}$ of band-crossing point is a manifold of the dimension $d-2$. Notice,
that there is an natural ambiguity in the definition of the normalized
eigenvectors $\mathbf{g}_{n,\zeta }\left( \mathbf{k}\right) $ of $\mathbf{L}%
\left( \mathbf{\mathbf{k}}\right) $ which is defined up to a complex number $%
\xi $ with $\left\vert \xi \right\vert =1$. This ambiguity may not allow an
eigenvector $\mathbf{g}_{n,\zeta }\left( \mathbf{k}\right) $ which can be a
locally smooth function in $\mathbf{k}$\ to be a uniquely defined continuous
function in $\mathbf{k}$ globally for all $\mathbf{k}\notin \sigma _{bc}$
because of a possibility of branching. But, importantly, the orthogonal
projector $\Pi _{n,\zeta }\left( \mathbf{\mathbf{k}}\right) $ on $\mathbf{g}%
_{n,\zeta }\left( \mathbf{k}\right) $ as defined by (\ref{Pin}) is uniquely
defined and, consequently, infinitely differentiable in $\mathbf{k}$ via the
spectral theory, \cite{Kato}, for all $\mathbf{k}\notin \sigma _{bc}$. Since
we consider $\mathbf{\hat{U}}\left( \mathbf{k}\right) $ as an element of the
space $L^{1}$ and $\sigma _{bc}$ is of zero Lebesgue measure considering $%
\mathbf{k}\notin \sigma _{bc}$ is sufficient for us.
\end{remark}

We introduce for vectors $\mathbf{\hat{u}}\in \mathbb{C}^{2J}$ their
expansion with respect to the orthonormal basis$\left\{ \mathbf{g}_{n,\zeta
}\left( \mathbf{k}\right) \right\} $: 
\begin{equation}
\mathbf{\hat{u}}\left( \mathbf{k}\right) =\sum_{n=1}^{J}\sum_{\zeta =\pm }%
\hat{u}_{n,\zeta }\left( \mathbf{k}\right) \mathbf{g}_{n,\zeta }\left( 
\mathbf{k}\right) =\sum_{n=1}^{J}\sum_{\zeta =\pm }\mathbf{\hat{u}}_{n,\zeta
}\left( \mathbf{k}\right) ,\ \mathbf{\hat{u}}_{n,\zeta }\left( \mathbf{k}%
\right) =\Pi _{n,\zeta }\left( \mathbf{\mathbf{k}}\right) \mathbf{\hat{u}}%
\left( \mathbf{k}\right)  \label{Uboldj}
\end{equation}%
and we refer to it as the \emph{modal decomposition} of $\mathbf{\hat{u}}%
\left( \mathbf{k}\right) $ and to $\hat{u}_{n,\zeta }\left( \mathbf{k}%
\right) $ as the \emph{modal coefficients} of $\mathbf{\hat{u}}\left( 
\mathbf{k}\right) $. Evidently 
\begin{equation}
\sum\nolimits_{n=1}^{j}\sum\nolimits_{\zeta =\pm }\Pi _{n,\zeta }\left( 
\mathbf{\mathbf{k}}\right) =I_{2J},\text{ where }I_{2J}\text{ is the }%
2J\times 2J\text{ identity matrix.}  \label{sumPi}
\end{equation}%
Notice that in view of the polynomial bound \ref{Lpol}) we can define the
action of the operator $\mathbf{L}\left( -\mathrm{i}\nabla _{\mathbf{r}%
}\right) $ on any Schwartz function $\mathbf{Y}\left( \mathbf{r}\right) $ by
the formula%
\begin{equation}
\widehat{\mathbf{L}\left( -\mathrm{i}\nabla _{\mathbf{r}}\right) \mathbf{Y}}%
\left( \mathbf{\mathbf{k}}\right) =\mathbf{L}\left( \mathbf{\mathbf{k}}%
\right) \mathbf{\hat{Y}}\left( \mathbf{\mathbf{k}}\right) ,\text{ where the
order of }\mathbf{L}\text{ does not exceed }p.  \label{Ldiff}
\end{equation}%
In a special case when all the entries of $\mathbf{L}\left( \mathbf{\mathbf{k%
}}\right) $ are polynomials (\ref{Ldiff}) turns into the action of the
differential operator with constant coefficients of order not exceeding $p$.

\subsection{Nonlinear part}

The nonlinear term $\hat{F}$ in (\ref{difeqfou}) is assumed to be a general
functional polynomial of the form 
\begin{gather}
\hat{F}\left( \mathbf{\hat{U}}\right) =\sum\nolimits_{m\in \mathfrak{M}_{F}}%
\hat{F}^{\left( m\right) }\left( \mathbf{\hat{U}}^{m}\right) ,\text{ where }%
\hat{F}^{\left( m\right) }\text{ is }m\text{-homogeneous polylinear operator,%
}  \label{Fseries} \\
\mathfrak{M}_{F}=\left\{ m_{1},\ldots ,m_{p}\right\} \subset \left\{
2,3,\ldots \right\} \text{ is a finite set, and }m_{F}=\max \left\{ m:m\in 
\mathfrak{M}_{F}\right\} .  \label{Fseries1}
\end{gather}%
The integer $m_{F}$ in (\ref{Fseries1}) is called the \emph{degree of the
functional polynomial} $\hat{F}$. For instance, if $\mathfrak{M}_{F}=\left\{
2\right\} $ or $\mathfrak{M}_{F}=\left\{ 3\right\} $ the polynomial $\hat{F}%
\,$\ is respectively homogeneous quadratic or cubic. Every $m$-linear
operator $\hat{F}^{\left( m\right) }$ in (\ref{Fseries}) is assumed to be of
the form of a convolution%
\begin{gather}
\hat{F}^{\left( m\right) }\left( \mathbf{\hat{U}}_{1},\ldots ,\mathbf{\hat{U}%
}_{m}\right) \left( \mathbf{k},\tau \right) =\int_{\mathbb{D}_{m}}\chi
^{\left( m\right) }\left( \mathbf{\mathbf{k}},\vec{k}\right) \mathbf{\hat{U}}%
_{1}\left( \mathbf{k}^{\prime }\right) \ldots \mathbf{\hat{U}}_{m}\left( 
\mathbf{k}^{\left( m\right) }\left( \mathbf{k},\vec{k}\right) \right) \,%
\mathrm{\tilde{d}}^{\left( m-1\right) d}\vec{k},  \label{Fmintr} \\
\text{where }\mathbb{D}_{m}=\mathbb{R}^{\left( m-1\right) d},\ \mathrm{%
\tilde{d}}^{\left( m-1\right) d}\vec{k}=\frac{\mathrm{d}\mathbf{k}^{\prime
}\ldots \,\mathrm{d}\mathbf{k}^{\left( m-1\right) }}{\left( 2\pi \right)
^{\left( m-1\right) d}},  \notag \\
\mathbf{k}^{\left( m\right) }\left( \mathbf{k},\vec{k}\right) =\mathbf{k}-%
\mathbf{k}^{\prime }-\ldots -\mathbf{k}^{\left( m-1\right) },\ \vec{k}%
=\left( \mathbf{k}^{\prime },\ldots ,\mathbf{k}^{\left( m\right) }\right) .
\label{conv}
\end{gather}%
indicating that the nonlinear operator $F^{\left( m\right) }\left( \mathbf{U}%
_{1},\ldots ,\mathbf{U}_{m}\right) $ is translation invariant (it may be
local or non-local). The quantities $\chi ^{\left( m\right) }$ in (\ref%
{Fmintr}) are called \emph{susceptibilities}. For numerous examples of
nonlinearities of the form similar to (\ref{Fseries}), (\ref{Fmintr}) see 
\cite{BF1}-\cite{BF7} and references therein. In what follows the nonlinear
term $\hat{F}$ in (\ref{difeqfou}) will satisfy the following conditions.

\begin{condition}[nonlinearity]
\label{cnonbound}The nonlinearity $\hat{F}\left( \mathbf{\hat{U}}\right) $
is assumed to be of the form (\ref{Fseries})-(\ref{Fmintr}). The
susceptibility $\chi ^{\left( m\right) }\left( \mathbf{\mathbf{k}},\mathbf{k}%
^{\prime },\ldots ,\mathbf{k}^{\left( m\right) }\right) $ is infinitely
differentiable for all $\mathbf{\mathbf{k}}$ and $\mathbf{k}^{\left(
j\right) }$ which are not band-crossing points, and is bounded, namely 
\begin{equation}
\left\Vert \chi ^{\left( m\right) }\right\Vert =\left( 2\pi \right)
^{-\left( m-1\right) d}\sup_{\mathbf{\mathbf{k}},\mathbf{k}^{\prime },\ldots
,\mathbf{k}^{\left( m\right) }\in \mathbb{R}^{d}\setminus \sigma
_{bc}}\left\vert \chi ^{\left( m\right) }\left( \mathbf{\mathbf{k}},\mathbf{k%
}^{\prime },\ldots ,\mathbf{k}^{\left( m\right) }\right) \right\vert \leq
C_{\chi },\ m\in \mathfrak{M}_{F},  \label{chiCR}
\end{equation}%
where the norm $\left\vert \chi ^{\left( m\right) }\left( \mathbf{k},\vec{k}%
\right) \right\vert $ of $m$-linear tensor $\chi ^{\left( m\right) }:\left( 
\mathbb{C}^{2J}\right) ^{m}\rightarrow \left( \mathbb{C}^{2J}\right) ^{m}$
for fixed $\mathbf{k},\vec{k}$ is defined by 
\begin{equation}
\left\vert \chi ^{\left( m\right) }\left( \mathbf{k},\vec{k}\right)
\right\vert =\sup_{\left\vert \mathbf{x}_{j}\right\vert \leq 1}\left\vert
\chi _{\ }^{\left( m\right) }\left( \mathbf{k},\vec{k}\right) \left( \mathbf{%
x}_{1},\ldots ,\mathbf{x}_{m}\right) \right\vert ,\text{ where }\left\vert 
\mathbf{x}\right\vert \text{ is the Euclidean norm}.  \label{normchi0}
\end{equation}
\end{condition}

When $\chi ^{\left( m\right) }\left( \mathbf{k},\vec{k}\right) $ depend on
small $\varrho $ or, more generally, on $\varrho ^{q}$, $q>0$, \ we
similarly have $\chi ^{\left( m\right) }\left( \mathbf{k},\vec{k},\varrho
^{q}\right) $. Many results of this paper extend to this case, in particular
if $\left\Vert \chi ^{\left( m\right) }\left( \mathbf{k},\vec{k},\varrho
^{q}\right) -\chi ^{\left( m\right) }\left( \mathbf{k},\vec{k},0\right)
\right\Vert \leq C_{\chi }^{\prime }\varrho ^{q}$ for $\varrho \leq 1$ \
then conditions of Corollary \ref{Corollary Fofrho} \ are fulfilled.

Note that since the tensors $\chi ^{\left( m\right) }\left( \mathbf{k},\vec{k%
}\right) $ are bounded, the dependence on $\left( \mathbf{k},\vec{k}\right) $
cannot be polynomial, therefore the original equation (\ref{difeqintr}) does
not include spatial derivatives but rather includes bounded
"pseudodifferential" operators. Note that this type of susceptibilities with
spatial dispersion \ is common in nonlinear optics, see \cite{ButCot}, \cite%
{MolNew}, \cite{SipeBhat}.

\subsection{Resonance invariant $nk$-spectrum}

In this section, relying on given dispersion relations $\omega _{n}\left( 
\mathbf{k}\right) \geq 0$, $n\in \left\{ 1,...,J\right\} $, \ we consider
resonance properties of $nk$-spectra $S$ and the corresponding $k$-spectra $%
K_{S}$ as defined in Definition \ref{dmwavepack}, i.e. 
\begin{equation}
S=\left\{ \left( n_{l},\mathbf{k}_{\ast l}\right) ,\ l=1,...,N\right\}
\subset \Sigma =\left\{ 1,...,J\right\} \times \mathbb{R}^{d},\ K_{S}\text{ }%
=\left\{ \mathbf{k}_{\ast _{l_{i}}},\ i=1,...,\left\vert K_{S}\right\vert
\right\} .  \label{ssp1}
\end{equation}%
We precede the formal description of the \emph{resonance invariance} (see
Definition \ref{Definition omclos}) with the following guiding physical
picture. Initially at $\tau =0$ the wave is a multi-wavepacket composed of
modes from a small vicinity of the $nk$-spectrum $S$. As the wave evolves
according to (\ref{difeqfou}) the polynomial nonlinearity inevitably
involves \ a larger set of modes $\left[ S\right] _{\text{out}}\supseteq S$,
but not all modes in $\left[ S\right] _{\text{out}}$\ are "equal" in
developing significant amplitudes. The qualitative picture is that whenever
certain interaction phase function (see (\ref{phim}) below) is not zero, the
fast time oscillations weaken effective nonlinear mode interaction and the
energy transfer from the original modes in $S$ to relevant modes from $\left[
S\right] _{\text{out}}$, keeping their magnitudes vanishingly small as $%
\beta ,\varrho \rightarrow 0$. There is a smaller set of modes $\left[ S%
\right] _{\text{out}}^{\text{res}}$ which can interact with modes from $S$
rather effectively and develop significant amplitudes. Now, 
\begin{equation}
\text{if }\left[ S\right] _{\text{out}}^{\text{res}}\subseteq S\text{ then }S%
\text{ is called resonance invariant.}  \label{ssp2}
\end{equation}%
In simpler situations the resonance invariance conditions turns into the
well-known in nonlinear optics phase and frequency matching conditions. For
instance, if $S$ contains $\left( n_{0},\mathbf{k}_{\ast l_{0}}\right) $ and
the dispersion relations allow for the second harmonic generation in another
band $n_{1}$ so that $2\omega _{n_{0}}\left( \mathbf{k}_{\ast l_{0}}\right)
=\omega _{n_{1}}\left( 2\mathbf{k}_{\ast l_{0}}\right) $, then for $S$ to be
resonance invariant it must contain $\left( n_{1},2\mathbf{k}_{\ast
l_{0}}\right) $ too.

Let us turn now to the rigorous constructions. First we introduce necessary
notations. Let $m\geq 2$ be an integer, $\vec{l}=\left(
l_{1},..,l_{m}\right) $,\ $l_{j}\in \left\{ 1,...,N\right\} $ be an integer
vector from $\left\{ 1,...,N\right\} ^{m}$ and $\vec{\zeta}=\left( \zeta
^{\left( 1\right) },,..,\zeta ^{\left( m\right) }\right) $, $\zeta ^{\left(
j\right) }\in \left\{ +1,-1\right\} $ be a binary vector from $\left\{
+1,-1\right\} ^{m}$. \emph{\ }Note that a pair $\left( \vec{\zeta},\vec{l}%
\right) $ naturally labels a sample string of the length $m$ composed of
elements $\left( \zeta ^{\left( j\right) },n_{l_{j}},\mathbf{k}_{\ast
l_{j}}\right) $ from the set $\left\{ +1,-1\right\} \times S$. Let us
introduce the sets 
\begin{gather}
\Lambda =\left\{ \left( \zeta ,l\right) :l\in \left\{ 1,...,N\right\} ,\
\zeta \in \left\{ +1,-1\right\} \right\} ,  \label{setlam} \\
\Lambda ^{m}=\left\{ \vec{\lambda}=\left( \lambda _{1},...,\lambda
_{m}\right) ,\ \lambda _{j}\in \Lambda ,\ j=1,...,m\right\} .  \notag
\end{gather}%
There is a natural one-to-one correspondence between $\Lambda ^{m}$ and $%
\left\{ -1,1\right\} ^{m}\times \left\{ 1,...,N\right\} ^{m}$ \ and we will
write, exploiting this correspondence 
\begin{equation}
\vec{\lambda}=\left( \left( \zeta ^{\prime },l_{1}\right) ,...,\left( \zeta
^{\left( m\right) },l_{m}\right) \right) =\left( \vec{\zeta},\vec{l}\right)
,\ \vec{\vartheta}\in \left\{ -1,1\right\} ^{m},\ \vec{l}\in \left\{
1,...,N\right\} ^{m}\text{ for }\vec{\lambda}\in \Lambda ^{m}.
\label{lamprop}
\end{equation}%
Let us introduce the following linear combination \ 
\begin{equation}
\varkappa _{m}\left( \vec{\lambda}\right) =\varkappa _{m}\left( \vec{\zeta},%
\vec{l}\right) =\sum\nolimits_{j=1}^{m}\zeta ^{\left( j\right) }\mathbf{k}%
_{\ast l_{j}}\text{ with }\zeta ^{\left( j\right) }\in \left\{ +1,-1\right\}
,  \label{kapzel}
\end{equation}%
and let $\left[ S\right] _{K,\text{out}}$ be the set of all its values as $%
\mathbf{k}_{\ast l_{j}}\in K_{S}$, $\vec{\lambda}\in \Lambda ^{m}$, namely 
\begin{equation}
\left[ S\right] _{K,\text{out}}=\dbigcup\nolimits_{m\in \mathfrak{M}%
_{F}}\dbigcup\nolimits_{\vec{\lambda}\in \Lambda ^{m}}\left\{ \varkappa
_{m}\left( \vec{\lambda}\right) \right\} .  \label{K1}
\end{equation}%
\ We call $\left[ S\right] _{K,\text{out}}$ \emph{output }$k$\emph{-spectrum}
of $K_{S}$. Everywhere in this paper \ we consider $nk$-spectra $S$ which
satisfy the following condition 
\begin{equation}
\left[ S\right] _{K,\text{out}}\dbigcap \sigma _{bc}=\varnothing .
\label{Skempty}
\end{equation}%
We also define the \emph{output }$nk$-\emph{spectrum} of $S$ by 
\begin{equation}
\left[ S\right] _{\text{out}}=\left\{ \left( n,\mathbf{k}\right) \in \left\{
1,...,J\right\} \times \mathbb{R}^{d}:n\in \left\{ 1,...,J\right\} ,\ 
\mathbf{k}\in \left[ S\right] _{K,\text{out}}\right\} .  \label{Sout}
\end{equation}%
We introduce the following functions%
\begin{equation}
\Omega _{1,m}\left( \vec{\lambda}\right) \left( \vec{k}_{\ast }\right)
=\sum\nolimits_{j=1}^{m}\zeta ^{\left( j\right) }\omega _{l_{j}}\left( 
\mathbf{k}_{\ast l_{j}}\right) ,\ \vec{k}_{\ast }=\left( \mathbf{k}_{\ast
1},...,\mathbf{k}_{\ast \left\vert K_{S}\right\vert }\right) ,\ \text{where }%
\mathbf{k}_{\ast l_{j}}\in K_{S},  \label{Om1}
\end{equation}%
\begin{equation}
\Omega \left( \zeta ,n,\vec{\lambda}\right) \left( \mathbf{k}_{\ast \ast },%
\vec{k}_{\ast }\right) =-\zeta \omega _{n}\left( \mathbf{k}_{\ast \ast
}\right) +\Omega _{1,m}\left( \vec{\lambda}\right) \left( \vec{k}_{\ast
}\right) ,  \label{Omzet}
\end{equation}%
where $\zeta =\pm 1$,$\ m\in \mathfrak{M}_{F}\ $as in (\ref{Fseries}). We
introduce these functions to apply later to phase functions (\ref{phim}).

Now we introduce the \emph{resonance equation} 
\begin{equation}
\Omega \left( \zeta ,n,\vec{\lambda}\right) \left( \zeta \varkappa
_{m}\left( \vec{\lambda}\right) ,\vec{k}_{\ast }\right) =0,\ \vec{l}\in
\left\{ 1,...,N\right\} ^{m},\ \vec{\zeta}\in \left\{ -1,1\right\} ^{m},
\label{Omeq0}
\end{equation}%
denoting by $P\left( S\right) $ the set of its solutions $\left( m,\zeta ,n,%
\vec{\lambda}\right) $. Such a solution is called $S$\emph{-internal }if%
\begin{equation*}
\left( n,\zeta \varkappa _{m}\left( \vec{\lambda}\right) \right) \in S,\text{
that is }n=n_{l_{0}},\ \zeta \varkappa _{m}\left( \vec{\lambda}\right) =%
\mathbf{k}_{\ast l_{0}},\ l_{0}\in \left\{ 1,...,N\right\} ,
\end{equation*}%
and we denote the corresponding $l_{0}=I\left( \vec{\lambda}\right) $. We
also denote by $P_{\text{int}}\left( S\right) \subset P\left( S\right) $ the
set of all $S$-internal solutions to (\ref{Omeq0}).

Now we consider the simplest solutions to (\ref{Omeq0}) which play an
important role. Keeping in mind that the string $\vec{l}$ can contain
several copies of a single value $l$, we can recast the sum in (\ref{Om1})
as follows: 
\begin{gather}
\Omega _{1,m}\left( \vec{\lambda}\right) =\Omega _{1,m}\left( \vec{\zeta},%
\vec{l}\right) =\sum\nolimits_{l=1}^{N}\delta _{l}\omega _{l}\left( \mathbf{k%
}_{\ast l}\right) ,\ \text{where }\delta _{l}=\left\{ 
\begin{array}{ccc}
\sum_{j\in \vec{l}^{-1}\left( l\right) }\zeta ^{\left( j\right) } & \text{if}
& \vec{l}^{-1}\left( l\right) \neq \varnothing \\ 
0 & \text{if} & \vec{l}^{-1}\left( l\right) =\varnothing%
\end{array}%
\right. ,  \label{rearr0} \\
\vec{l}^{-1}\left( l\right) =\left\{ j:l_{j}=l,\ 1\leq j\leq m\right\} ,\ 
\vec{l}=\left( l_{1},\ldots ,l_{m}\right) ,\ 1\leq l\leq N.  \notag
\end{gather}%
Let us call a solution $\left( m,\zeta ,n,\vec{\lambda}\right) \in P\left(
S\right) $ of (\ref{Omeq0}) \emph{universal} if it has the following
properties: (i) only a single coefficient out of all $\delta _{l}$ in (\ref%
{rearr0}) is nonzero, namely for some $I_{0}$ we have $\delta _{I_{0}}=\pm 1$
and $\delta _{l}=0$ for $l\neq I_{0}$; (ii) $n=n_{I_{0}}$ \ and $\zeta
=\delta _{I_{0}}$. \emph{A justification for calling such a solution
universal comes from the fact that if it is a solution for one }$\vec{k}%
_{\ast }$ \emph{it is a solution for any other }$\vec{k}_{\ast }\in \mathbb{R%
}^{d}$\emph{. }We denote the set of universal solutions to (\ref{Omeq0}) by $%
P_{\text{univ}}\left( S\right) $, and note that a universal solution is a $S$%
-internal solution with $I\left( \vec{\lambda}\right) =I_{0}$ implying 
\begin{equation}
P_{\text{univ}}\left( S\right) \subseteq P_{\text{int}}\left( S\right) .
\label{gint}
\end{equation}%
Indeed, observe that for $\delta _{l}$ as in (\ref{rearr0} 
\begin{equation}
\varkappa _{m}\left( \vec{\lambda}\right) =\varkappa _{m}\left( \vec{\zeta},%
\vec{l}\right) =\sum\nolimits_{j=1}^{m}\zeta ^{\left( j\right) }\mathbf{k}%
_{\ast l_{j}}=\sum\nolimits_{l=1}^{N}\delta _{l}\mathbf{k}_{\ast l}
\label{rearr1}
\end{equation}%
implying $\varkappa _{m}\left( \vec{\lambda}\right) =\delta _{I_{0}}\mathbf{k%
}_{\ast I_{0}}$ and$\ \ \zeta \varkappa _{m}\left( \vec{\lambda}\right)
=\delta _{I_{0}}^{2}\mathbf{k}_{\ast I_{0}}=\mathbf{k}_{\ast I_{0}}$. Then
equation (\ref{Omeq0}) \ is obviously satisfied and $\left( n,\zeta
\varkappa _{m}\left( \vec{\lambda}\right) \right) =\left( n_{I_{0}},\mathbf{k%
}_{\ast I_{0}}\right) \in S$.

\begin{example}[Universal solutions]
\label{Example counterprop}Suppose there is just a single band, i.e. $J=1$,
a symmetric dispersion relation $\omega _{1}\left( -\mathbf{k}\right)
=\omega _{1}\left( \mathbf{k}\right) $, a cubic nonlinearity $F$ with $%
\mathfrak{M}_{F}=\left\{ 3\right\} $. First let us take the simplest $nk$%
-spectrum\emph{\ } $S_{1}=\left\{ \left( 1,\mathbf{k}_{\ast }\right)
\right\} $, that is $N=1$. Then $\Omega _{1,3}\left( \vec{\lambda}\right)
\left( \vec{k}_{\ast }\right) =\delta _{1}\omega _{1}\left( \mathbf{k}_{\ast
}\right) $\ and $\varkappa _{m}\left( \vec{\lambda}\right) =\delta _{1}%
\mathbf{k}_{\ast }$ where we use notation (\ref{rearr0}). The universal
solution set has the form $P_{\text{univ}}\left( S_{1}\right) =\left\{
\left( 3,\zeta ,1,\vec{\lambda}\right) :\vec{\lambda}\in \Lambda _{\zeta },\
\zeta =\pm \right\} $ where $\Lambda _{+}$ consists of vectors $\left(
\lambda _{1},\lambda _{2},\lambda _{3}\right) $ of the form $\left( \left(
-,1\right) ,\left( +,1\right) ,\left( +,1\right) \right) $, $\ \left( \left(
+,1\right) ,\left( -,1\right) ,\left( +,1\right) \right) $ and $\left(
\left( +,1\right) ,\left( +,1\right) ,\left( -,1\right) \right) $.
Obviously, $P_{\text{univ}}\left( S_{1}\right) =P_{\text{int}}\left(
S_{1}\right) $. In the next example we take the $nk$-spectrum\emph{\ }$%
S=\left\{ \left( 1,\mathbf{k}_{\ast }\right) ,\left( 1,-\mathbf{k}_{\ast
}\right) \right\} $, that is $N=2$ and $\mathbf{k}_{\ast 1}=\mathbf{k}_{\ast
},\mathbf{k}_{\ast 2}=-\mathbf{k}_{\ast }$. This example is typical for two
counterpropagating waves. Then $\Omega _{1,3}\left( \vec{\lambda}\right)
\left( \vec{k}_{\ast }\right) =\sum\nolimits_{j=1}^{3}\zeta ^{\left(
j\right) }\omega _{l_{j}}\left( \mathbf{k}_{\ast l_{j}}\right) =\left(
\delta _{1}+\delta _{2}\right) \omega _{1}\left( \mathbf{k}_{\ast }\right) $
\ and $\varkappa _{m}\left( \vec{\lambda}\right)
=\sum\nolimits_{j=1}^{m}\zeta ^{\left( j\right) }\mathbf{k}_{\ast
l_{j}}=\delta _{1}\mathbf{k}_{\ast 1}+\delta _{2}\mathbf{k}_{\ast 2}=\left(
\delta _{1}-\delta _{2}\right) \mathbf{k}_{\ast }$ where we use notation (%
\ref{rearr0}). The universal solution set has the form $P_{\text{univ}%
}\left( S\right) =\left\{ \left( 3,\zeta ,1,\vec{\lambda}\right) :\vec{%
\lambda}\in \Lambda _{\zeta },\zeta =\pm \right\} $ where $\Lambda _{+}$
consists of vectors $\left( \lambda _{1},\lambda _{2},\lambda _{3}\right) $ $%
\ $of the form $\left( \left( +,1\right) ,\left( -,1\right) ,\left(
+,1\right) \right) $, $\left( \left( +,1\right) ,\left( -,1\right) ,\left(
+,2\right) \right) $, $\left( \left( +,2\right) ,\left( -,2\right) ,\left(
+,1\right) \right) $, $\left( \left( +,2\right) ,\left( -,2\right) ,\left(
+,2\right) \right) $, and vectors obtained from the listed ones by
permutations of coordinates $\lambda _{1},\lambda _{2},\lambda _{3}$. The
solutions from $P_{\text{int}}\left( S\right) $ have to satisfy $\left\vert
\delta _{1}-\delta _{2}\right\vert =1$ and $\left\vert \delta _{1}+\delta
_{2}\right\vert =1$ which is possible only if $\delta _{1}\delta _{2}=0$.
Since $\zeta =\delta _{1}+\delta _{2}$ we have $\zeta \varkappa _{m}\left( 
\vec{\lambda}\right) =\left( \delta _{1}^{2}-\delta _{2}^{2}\right) \mathbf{k%
}_{\ast }$ and $\zeta \varkappa _{m}\left( \vec{\lambda}\right) =\mathbf{k}%
_{\ast 1}$ if $\left\vert \delta _{1}\right\vert =1$ or $\zeta \varkappa
_{m}\left( \vec{\lambda}\right) =\mathbf{k}_{\ast 2}$ if $\left\vert \delta
_{2}\right\vert =1$. Hence $P_{\text{int}}\left( S\right) =P_{\text{univ}%
}\left( S\right) $ in this case. Note that if we set $S_{2}=\left\{ \left(
1,-\mathbf{k}_{\ast }\right) \right\} $ then $S=S_{1}\cup S_{2}$ but $P_{%
\text{int}}\left( S\right) $ is larger than $P_{\text{int}}\left(
S_{1}\right) \cup P_{\text{int}}\left( S_{2}\right) $. This can be
interpreted as follows. When only \ modes from $S_{1}$ are excited, the
modes from $S_{2}$ remain non-excited. But when the both $S_{1}$ and $S_{2}$
are excited, there is a resonance effect of $S_{1}$ onto $S_{2}$,
represented, for example, by $\vec{\lambda}=\left( \left( +,1\right) ,\left(
-,1\right) ,\left( +,2\right) \right) $, which involves the mode $\zeta
\varkappa _{m}\left( \vec{\lambda}\right) =\mathbf{k}_{\ast 2}$.
\end{example}

Now we are ready to define resonance invariant spectra. First, we introduce
a subset $\left[ S\right] _{\text{out}}^{\text{res}}$ of $\left[ S\right] _{%
\text{out}}$ by the formula%
\begin{eqnarray}
\left[ S\right] _{\text{out}}^{\text{res}} &=&\left\{ \left( n,\mathbf{k}%
_{\ast \ast }\right) \in \left[ S\right] _{\text{out}}:\mathbf{k}_{\ast \ast
}=\zeta \varkappa _{m}\left( \vec{\lambda}\right) ,\ m\in \mathfrak{M}_{F},%
\text{ where}\right.  \label{Sresout} \\
&&\left. \left( m,\zeta ,n,\vec{\lambda}\right) \text{ is a solution of (\ref%
{Omeq0})}\right\} ,  \notag
\end{eqnarray}%
\ calling it \emph{resonant output spectrum} of $S$, and then we define%
\begin{equation}
\text{resonance selection operation }\mathcal{R}\left( S\right) =S\cup \left[
S\right] _{\text{out}}^{\text{res}}.  \label{hull}
\end{equation}

\begin{definition}[resonance invariant $nk$-spectrum]
\label{Definition omclos} The $nk$-spectrum $S$ is called \emph{resonance
invariant} if $\mathcal{R}\left( S\right) =S$ or, equivalently, $\left[ S%
\right] _{\text{out}}^{\text{res}}\subseteq S$. The $nk$-spectrum\ $S$ is
called \emph{universally} \emph{resonance invariant }if $\mathcal{R}\left(
S\right) =S$ and $P_{\text{univ}}\left( S\right) =P_{\text{int}}\left(
S\right) $.
\end{definition}

It is worth noticing that even when a $nk$-spectrum is not resonance
invariant often it can be easily extended to a resonance invariant one.
Namely, if $\mathcal{R}^{j}\left( S\right) \cap \sigma _{bc}=\varnothing $
for all $j$ then the set 
\begin{equation*}
\mathcal{R}^{\infty }\left( S\right) =\dbigcup\nolimits_{j=1}^{\infty }%
\mathcal{R}^{j}\left( S\right) \subset \Sigma =\left\{ 1,...,J\right\}
\times \mathbb{R}^{d}\text{ }
\end{equation*}%
is resonance invariant. In addition to that, $\mathcal{R}^{\infty }\left(
S\right) $ is always at most countable. Usually it is finite i.e. $\mathcal{R%
}^{\infty }\left( S\right) =\mathcal{R}^{p}\left( S\right) $ for a finite $p$%
, see examples below and we also show below that $\mathcal{R}^{\infty
}\left( S\right) =S$ \ for generic $K_{S}$.

\begin{example}[resonance invariant $nk$-spectra \ for quadratic nonlinearity%
]
\label{Example quad} Suppose there is a single band, i.e. $J=1$, with a
symmetric dispersion relation, and a quadratic nonlinearity $F$, that is $%
\mathfrak{M}_{F}=\left\{ 2\right\} $. Let us assume that $\mathbf{k}_{\ast
}\neq 0$, $\mathbf{k}_{\ast },2\mathbf{k}_{\ast },\mathbf{0}$ are not
band-crossing points and look at two examples. First, suppose that $2\omega
_{1}\left( \mathbf{k}_{\ast }\right) \neq \omega _{1}\left( 2\mathbf{k}%
_{\ast }\right) $ (no second harmonic generation) and $\omega _{1}\left( 
\mathbf{0}\right) \neq 0$. Let us set the $nk$-spectrum to be the set $%
S_{1}=\left\{ \left( 1,\mathbf{k}_{\ast }\right) \right\} $, then $S_{1}$ is
resonance invariant. Indeed, $K_{S_{1}}=\left\{ \mathbf{k}_{\ast }\right\} $%
, $\left[ S_{1}\right] _{K,\text{out}}=\left\{ \mathbf{0},2\mathbf{k}_{\ast
},-2\mathbf{k}_{\ast }\right\} $, $\left[ S_{1}\right] _{\text{out}}=\left\{
\left( 1,\mathbf{0}\right) ,\left( 1,2\mathbf{k}_{\ast }\right) ,\left( 1,-2%
\mathbf{k}_{\ast }\right) \right\} $ and an elementary examination shows
that $\left[ S_{1}\right] _{\text{out}}^{\text{res}}=\varnothing \subset
S_{1}$ implying $\mathcal{R}\left( S_{1}\right) =S_{1}$. For the second
example let us assume $\omega _{1}\left( \mathbf{0}\right) \neq 0$ and $%
2\omega _{1}\left( \mathbf{k}_{\ast }\right) =\omega _{1}\left( 2\mathbf{k}%
_{\ast }\right) $, that is the second harmonic generation is allowed. Here $%
\left[ S_{1}\right] _{\text{out}}^{\text{res}}=\left\{ \left( 1,2\mathbf{k}%
_{\ast }\right) \right\} $ and $\mathcal{R}\left( S_{1}\right) =\left\{
\left( 1,\mathbf{k}_{\ast }\right) ,\left( 1,2\mathbf{k}_{\ast }\right)
\right\} $ implying $\mathcal{R}\left( S_{1}\right) \neq S_{1}$ and, hence, $%
S_{1}$ is not resonance invariant. Suppose now that $4\mathbf{k}_{\ast },3%
\mathbf{k}_{\ast }\notin \sigma _{bc}$ and \ $\omega _{1}\left( \mathbf{0}%
\right) \neq 0$, $\omega _{1}\left( 4\mathbf{k}_{\ast }\right) \neq 2\omega
_{1}\left( 2\mathbf{k}_{\ast }\right) $, $\omega _{1}\left( 3\mathbf{k}%
_{\ast }\right) \neq \omega _{1}\left( \mathbf{k}_{\ast }\right) +\omega
_{1}\left( 2\mathbf{k}_{\ast }\right) $ and let us set $S_{2}=\left\{ \left(
1,\mathbf{k}_{\ast }\right) ,\left( 1,2\mathbf{k}_{\ast }\right) \right\} $.
An elementary examination shows that $S_{2}$ is resonance invariant. Note
that $S_{2}$ can be obtained by iterating the resonance selection operator,
namely $S_{2}=\mathcal{R}\left( \mathcal{R}\left( S_{1}\right) \right) $.
Note also that $P_{\text{univ}}\left( S_{2}\right) \neq P_{\text{int}}\left(
S_{2}\right) $. Notice that $\omega _{1}\left( \mathbf{0}\right) =0$ is a
special case since $\mathbf{k}=\mathbf{0}$ is a band-crossing point, and it
requires a special treatment.
\end{example}

\begin{example}[resonance invariant $nk$-spectra for cubic nonlinearity]
\label{Example cube}Let us consider one-band case with symmetric dispersion
relation and a cubic \ nonlinearity that is $\mathfrak{M}_{F}=\left\{
3\right\} $. First we take $S_{1}=\left\{ \left( 1,\mathbf{k}_{\ast }\right)
\right\} $, we assume that $\mathbf{k}_{\ast },3\mathbf{k}_{\ast }$\ are not
band-crossing points, implying $\left[ S_{1}\right] _{K,\text{out}}=\left\{ 
\mathbf{k}_{\ast },-\mathbf{k}_{\ast },3\mathbf{k}_{\ast },-3\mathbf{k}%
_{\ast }\right\} $. We have $\Omega _{1,3}\left( \vec{\lambda}\right) \left( 
\vec{k}_{\ast }\right) =\sum\nolimits_{j=1}^{3}\zeta ^{\left( j\right)
}\omega _{1}\left( \mathbf{k}_{\ast }\right) =\delta _{1}\omega _{1}\left( 
\mathbf{k}_{\ast }\right) $ \ and $\varkappa _{m}\left( \vec{\lambda}\right)
=\delta _{1}\mathbf{k}_{\ast }$ where we use notation (\ref{rearr0}), $%
\delta _{1}$ takes values $1,-1,3,-3$ . If $3\omega _{1}\left( \mathbf{k}%
_{\ast }\right) \neq \omega _{1}\left( 3\mathbf{k}_{\ast }\right) $ then (%
\ref{Omeq0}) has a solution only if $\ \left\vert \delta _{1}\right\vert =1$
and $\delta _{1}=\zeta $, hence $\zeta \varkappa _{m}\left( \vec{\lambda}%
\right) =\mathbf{k}_{\ast }$ and every solution is internal. \ \ Therefore, $%
\left[ S_{1}\right] _{\text{out}}^{\text{res}}=\varnothing $ and $\mathcal{R}%
\left( S_{1}\right) =S_{1}$.\ Now consider the case associated with the
third harmonic generation, namely $3\omega _{1}\left( \mathbf{k}_{\ast
}\right) =\omega _{1}\left( 3\mathbf{k}_{\ast }\right) $ and assume that $%
\omega _{1}\left( 3\mathbf{k}_{\ast }\right) +2\omega _{1}\left( \mathbf{k}%
_{\ast }\right) \neq \omega _{1}\left( 5\mathbf{k}_{\ast }\right) $, $%
3\omega _{1}\left( 3\mathbf{k}_{\ast }\right) \neq \omega _{1}\left( 9%
\mathbf{k}_{\ast }\right) $, $2\omega _{1}\left( 3\mathbf{k}_{\ast }\right)
+\omega _{1}\left( \mathbf{k}_{\ast }\right) \neq \omega _{1}\left( 7\mathbf{%
k}_{\ast }\right) $, $2\omega _{1}\left( 3\mathbf{k}_{\ast }\right) -\omega
_{1}\left( \mathbf{k}_{\ast }\right) \neq \omega _{1}\left( 5\mathbf{k}%
_{\ast }\right) $. An\ elementary examination shows that the set $%
S_{4}=\left\{ \left( 1,3\mathbf{k}_{\ast }\right) ,\left( 1,\mathbf{k}_{\ast
}\right) ,\left( 1,-\mathbf{k}_{\ast }\right) \left( 1,-3\mathbf{k}_{\ast
}\right) \right\} $ satisfies $\mathcal{R}\left( S_{4}\right) =S_{4}$.
Consequently, a multiwavepacket having $S_{4}$ as its resonance invariant $%
nk $-spectrum involves the third harmonic generation and, according to
Theorem \ref{Theorem invarwave}, it is preserved under nonlinear evolution.
\end{example}

The above examples indicate that in simple cases the conditions on $\ 
\mathbf{k}_{\ast }$ which can make $S$ non-invariant with respect to $%
\mathcal{R}$ have a form of several algebraic equations, therefore, for
almost all $\mathbf{k}_{\ast }$ such spectra $S$ are resonance invariant.
The examples also show that if we fix $S$\ and dispersion relations then we
can include $S$ in larger spectrum $S^{\prime }=\mathcal{R}^{p}\left(
S\right) $ \ using repeated application of the operation $\mathcal{R}$ to $S$%
, and often the resulting extended $nk$-spectrum $S^{\prime }$ is resonance
invariant. We show in the following section that $nk$-spectrum $S$ \ with
generic $K_{S}$ is universally resonance invariant.

Note that the concept of resonance invariant $nk$-spectrum gives a
mathematical description of such fundamental concepts of nonlinear optics as
phase matching, frequency matching, four wave interaction in cubic media and
three wave interaction in quadratic media. If a multi-wavepacket has a
resonance invariant spectrum, all these phenomena may take place in the
internal dynamics of \ the multi-wavepacket, but do not lead to resonant \
interactions with continuum of all remaining modes.

\subsection{Genericity of the $nk$-spectrum invariance condition}

In simpler situations, when the number of bands $J$ and wavepackets $N$ are
not too large, the resonance invariance of $nk$- spectrum can be easily
verified as above in Examples \ref{Example quad}, \ref{Example cube}, but
what one can say if $J$ or $N$ are large, or if the dispersion relations are
not explicitly given? We show below that in properly defined non-degenerate
cases a small variation of $K_{S}$ makes $S$ universally resonance
invariant, i.e. the resonance invariance is a generic phenomenon..

Assume that the dispersion relations $\omega _{n}\left( \mathbf{k}\right)
\geq 0$, $n\in \left\{ 1,...,J\right\} $ are given. Observe then that $%
\Omega _{m}\left( \zeta ,n,\vec{\lambda}\right) =\Omega _{m}\left( \zeta ,n,%
\vec{\lambda}\right) \left( \mathbf{k}_{\ast 1},...,\mathbf{k}_{\ast
\left\vert K_{S}\right\vert }\right) $ defined by (\ref{Omzet}) is a
continuous function of $\mathbf{k}_{\ast l}\notin \sigma _{bc}$ for every $%
m,\zeta ,n,\vec{\lambda}$.

\begin{definition}[$\protect\omega $-degenerate dispersion relations]
\label{Definition degen}We call dispersion relations $\omega _{n}\left( 
\mathbf{k}\right) $, $n=1,...,J$, $\omega $-degenerate if there exists such
a point $\mathbf{k}_{\ast }\in \mathbb{R}^{d}\setminus \sigma _{bc}$ that
for all $\mathbf{k}$ \ in a neighborhood of $\mathbf{k}_{\ast }$ at least
one of the following four conditions holds: (i) the relations are linearly
dependent, namely $\sum_{n=0}^{J}C_{n}\omega _{n}\left( \mathbf{k}\right)
=c_{0},$ where \ all $C_{n}$ are integers, one of which is nonzero, and the $%
c_{0}$ is a constant; (ii) at least one of $\omega _{n}\left( \mathbf{k}%
\right) $ is a linear function; (iii) at least one of $\omega _{n}\left( 
\mathbf{k}\right) $ satisfies equation $C\omega _{n}\left( \mathbf{k}\right)
=\omega _{n}\left( C\mathbf{k}\right) $ with some $n$ and integer $C\neq \pm
1$; (iv) at least one of $\omega _{n}\left( \mathbf{k}\right) $ satisfies
equation $\omega _{n}\left( \mathbf{k}\right) =\omega _{n^{\prime }}\left( -%
\mathbf{k}\right) $ where $n^{\prime }\neq n$.
\end{definition}

Note that fulfillment of any of the four conditions in Definition \ref%
{Definition degen} makes impossible turning some non resonance invariant
sets into resonance invariant ones by a variation of $\mathbf{k}_{\ast l}$.
For instance, if $\mathfrak{M}_{F}=\left\{ 2\right\} $ as in Example \ref%
{Example quad} and $2\omega _{1}\left( \mathbf{k}\right) =\omega _{1}\left( 2%
\mathbf{k}\right) $ for all $\mathbf{k}$ in an open set $G$ then the set $%
\left\{ \left( 1,\mathbf{k}_{\ast }\right) \right\} $ with $\mathbf{k}_{\ast
}\in G$ cannot be made resonance invariant by a small variation of $\mathbf{k%
}_{\ast }$. Below we show that if dispersion relations are not $\omega $%
-degenerate, then a small variation of $\mathbf{k}_{\ast l}$ turns non
resonance invariant sets into resonance invariant.

\begin{theorem}
\label{Theorem alt}If $\Omega _{m}\left( \zeta ,n_{0},\vec{\lambda}\right)
\left( \mathbf{k}_{\ast 1}^{\prime },...,\mathbf{k}_{\ast \left\vert
K_{S}\right\vert }^{\prime }\right) =0$ on a cylinder $G$\ in $\left( 
\mathbb{R}^{d}\setminus \sigma _{bc}\right) ^{\left\vert K_{S}\right\vert }$
which is a product of small balls $G_{i}\subset \left( \mathbb{R}%
^{d}\setminus \sigma _{bc}\right) $ then either $\left( m,\zeta ,n_{0},\vec{%
\lambda}\right) \in P_{\text{univ}}\left( S\right) $ or dispersive relations 
$\omega _{n}\left( \mathbf{k}\right) $ are $\omega $-degenerate as in
Definition \ref{Definition degen}.
\end{theorem}

\begin{proof}
Collecting similar terms in (\ref{Omeq0})\ we obtain the following equation
for $\mathbf{k}_{i}$ from $G_{i}$: 
\begin{equation}
\sum\nolimits_{n=1}^{J}\sum\nolimits_{i=1}^{\left\vert K_{S}\right\vert
}\delta _{in}^{\prime }\omega _{n}\left( \mathbf{k}_{i}\right) =\zeta \omega
_{n_{0}}\left( \sum\nolimits_{i=1}^{\left\vert K_{S}\right\vert }\delta
_{i}^{\prime }\mathbf{k}_{i}\right) \text{ where }\delta _{in}^{\prime
},\delta _{i}^{\prime }\text{ are integers.}  \label{eqom}
\end{equation}%
Comparing (\ref{eqom}) with (\ref{rearr0}) we see that $\delta _{in}^{\prime
}$\ may be non-zero only if $\left( n,\mathbf{k}_{i}\right) \in S$, that is $%
\left( n,\mathbf{k}_{i}\right) =\left( n_{l},\mathbf{k}_{l}\right) $ with $%
l\in \left\{ 1,...,N\right\} $, where $l=l\left( i,n\right) $ is uniquely
determined and $\delta _{in}^{\prime }=\delta _{l}$ with $\delta _{l}$ as in
(\ref{rearr0}). If there are two nonzero coefficients $\delta _{i}$ in (\ref%
{eqom}) we use an elementary Proposition \ref{Proposition sum} below
noticing that we are in the case (ii) of Definition \ref{Definition degen}.
If we do not have two nonzero $\delta _{i}^{\prime }$ then either all $%
\delta _{i}^{\prime }=0$ or only one $\delta _{i}^{\prime }=\delta
_{i_{0}}^{\prime }\neq 0$. If all $\delta _{i}^{\prime }=0$ then the
right-hand side of (\ref{eqom}) turns into $\omega _{n_{0}}\left( 0\right) $
and, $G_{i}\subset \left( \mathbb{R}^{d}\setminus \sigma _{bc}\right) $, $%
\omega _{n_{0}}\left( 0\right) \neq 0$. Hence, for every $i$ the sum $%
\sum_{n=1}^{J}\delta _{in}^{\prime }\omega _{n}\left( \mathbf{k}_{i}\right)
\ $is constant, one of $\delta _{in}^{\prime }$ is non-zero and we are in
the case (i) of Definition \ref{Definition degen}. If only one $\delta
_{i}^{\prime }\neq 0$ with $i=i_{0}$ we have 
\begin{equation}
\sum\nolimits_{n=1}^{J}\sum\nolimits_{i=1}^{\left\vert K_{S}\right\vert
}\delta _{in}^{\prime }\omega _{n}\left( \mathbf{k}_{i}\right) =\zeta \omega
_{n_{0}}\left( \delta _{i_{0}}^{\prime }\mathbf{k}_{i_{0}}\right) \text{ for
all }\mathbf{k}_{i}\in G_{i},\ \mathbf{k}_{i_{0}}\in G_{i_{0}},  \label{oneb}
\end{equation}%
implying linear dependence of the dispersion relations, namely 
\begin{equation*}
\sum\nolimits_{n=1}^{J}\delta _{in}^{\prime }\omega _{n}\left( \mathbf{k}%
_{i}\right) =C_{i},\ i\neq i_{0},\text{ where }C_{i}\text{ are constant.}
\end{equation*}%
The above equations would not imply linear dependence as in case (i) of
Definition \ref{Definition degen} only if 
\begin{equation}
\delta _{in}^{\prime }=0,\ i\neq i_{0},\ n=1,...,J,  \label{bprime0}
\end{equation}%
and in this case the equality (\ref{oneb}) takes the form 
\begin{equation}
\sum\nolimits_{n=1}^{J}\delta _{i_{0}n}^{\prime }\omega _{n}\left( \mathbf{k}%
_{i_{0}}\right) =\zeta \omega _{n_{0}}\left( \delta _{i_{0}}^{\prime }%
\mathbf{k}_{i_{0}}\right) \text{ for all }\mathbf{k}_{i_{0}}\in G_{i_{0}}.
\label{onec}
\end{equation}%
Note that in this case we deduce from (\ref{kapzel}) and (\ref{Om1}) that $%
\sum_{n=1}^{J}\delta _{i_{0}n}^{\prime }=\delta _{i_{0}}^{\prime }$. If $%
\left\vert \delta _{i_{0}}^{\prime }\right\vert \neq 1$ we are in the case
(iii) of Definition \ref{Definition degen}, whereas if $\left\vert \delta
_{i_{0}}^{\prime }\right\vert =1$ and $n\neq n_{0}$ we are in the case (iv)
of Definition \ref{Definition degen}. If $\left\vert \delta _{i_{0}}^{\prime
}\right\vert =1$ and $n=n_{0}$ (\ref{onec}) turns into $\delta
_{i_{0}}^{\prime }\omega _{n_{0}}\left( \mathbf{k}_{i_{0}}\right) =\zeta
\omega _{n_{0}}\left( \delta _{i_{0}}^{\prime }\mathbf{k}_{i_{0}}\right) $.
Since $\omega _{n_{0}}>0$ it implies $\delta _{i_{0}}^{\prime }=\zeta $ and $%
\omega _{n_{0}}\left( \mathbf{k}_{i_{0}}\right) =\omega _{n_{0}}\left( \zeta
\delta _{i_{0}}^{\prime }\mathbf{k}_{i_{0}}\right) $. Hence, in this case $%
\left( m,\zeta ,n_{0},\vec{\lambda}\right) \in P_{\text{univ}}\left(
S\right) $, and since all possibilities are exhausted the proof is complete.
\end{proof}

\begin{theorem}[genericity of resonance invariance]
\label{Proposition generic} Assume that dispersive relations $\omega
_{n}\left( \mathbf{k}\right) $ are continuous and not $\omega $-degenerate
as in Definition \ref{Definition degen}. Let $\mathcal{K}_{\text{rinv}}$ be
a set of points $\left( \mathbf{k}_{\ast 1},...,\mathbf{k}_{\ast \left\vert
K_{S}\right\vert }\right) $ such that there exists a universally resonance
invariant $nk$-spectrum $S$ for which its $k$-spectrum $K_{S}=\left\{ 
\mathbf{k}_{\ast 1},...,\mathbf{k}_{\ast \left\vert K_{S}\right\vert
}\right\} $. Then $\mathcal{K}_{\text{rinv}}$ is open and everywhere dense
set in $\left( \mathbb{R}^{d}\setminus \sigma _{bc}\right) ^{\left\vert
K_{S}\right\vert }$.
\end{theorem}

\begin{proof}
The fact that $\mathcal{K}_{\text{rinv}}$ is open follows from the
Definition \ref{Definition omclos} and the continuity in $\mathbf{k}$ of the
dispersion relations $\omega _{n}\left( \mathbf{k}\right) $. Let $G$ be a
small open ball such that its closure$\ \bar{G}\subset \left( \mathbb{R}%
^{d}\setminus \sigma _{bc}\right) ^{\left\vert K_{S}\right\vert }$. It
suffices to prove that $\bar{G}\cap \mathcal{K}_{\text{rinv}}$ contains at
least one point $\left( \mathbf{k}_{\ast 1},...,\mathbf{k}_{\ast \left\vert
K_{S}\right\vert }\right) $. For a given finite set $\mathfrak{M}_{F}$ let
us consider all possible%
\begin{equation*}
\left( m,\zeta ,n_{0},\vec{\lambda}\right) \in \dbigcup\nolimits_{m\in 
\mathfrak{M}_{F}}\times \left\{ -1,1\right\} \times \left\{ 1,...,J\right\}
\times \Lambda ^{m}
\end{equation*}%
which are not universal solutions to (\ref{Omeq0}), and for a given $\left(
m,\zeta ,n_{0},\vec{\lambda}\right) $ let $G_{0}\left( m,\zeta ,n_{0},\vec{%
\lambda}\right) $ be a set of solutions $\left( \mathbf{k}_{1},...,\mathbf{k}%
_{\left\vert K_{S}\right\vert }\right) $ to (\ref{Omeq0}) in $\bar{G}$, and
notice that it is a closed set.\ Let now $G_{0}\left( S\right) \subset \bar{G%
}$ be the union of the sets $G_{0}\left( m,\zeta ,n_{0},\vec{\lambda}\right) 
$ over all $\left( m,\zeta ,n_{0},\vec{\lambda}\right) \in P\left( S\right)
\setminus P_{\text{univ}}\left( S\right) $ and let us show that $G_{0}\left(
S\right) \neq G$. Indeed, suppose that $G_{0}\left( S\right) =G$ and hence $%
G\ $is a finite union of closed sets. According to the Baire's theorem one
of the sets $G_{0}\left( m,\zeta ,n_{0},\vec{\lambda}\right) $ with $\left(
m,\zeta ,n_{0},\vec{\lambda}\right) \in P\left( S\right) \setminus P_{\text{%
univ}}\left( S\right) $ must have a nonempty interior. Then, according to
Theorem \ref{Theorem alt}, the dispersion relations $\omega _{n}\left( 
\mathbf{k}\right) $ are $\omega $-degenerate as in Definition \ref%
{Definition degen} contradicting the conditions of the theorem. Hence, there
is always a point $\left( \mathbf{k}_{\ast 1},...,\mathbf{k}_{\ast
\left\vert K_{S}\right\vert }\right) \in P\left( S\right) \setminus P_{\text{%
univ}}\left( S\right) $ that completes the proof.
\end{proof}

The proof of the next statement is elementary and we skip it.

\begin{proposition}
\label{Proposition sum}Let $f_{1}\left( \mathbf{k}\right) $, $f_{2}\left( 
\mathbf{k}\right) $, $f_{3}\left( \mathbf{k}\right) $ be real-valued and
continuous functions respectively in neighborhoods of $\ \mathbf{k}_{\ast 1}$%
, $\mathbf{k}_{\ast 2}$, $\mathbf{k}_{\ast 1}+\mathbf{k}_{\ast 2}\mathbf{\ }$%
in $\mathbb{R}^{d}$. Assume that the following equation 
\begin{equation*}
f_{1}\left( \mathbf{k}_{1}\right) +f_{2}\left( \mathbf{k}_{2}\right)
=f_{3}\left( \delta _{1}\mathbf{k}_{1}+\delta _{2}\mathbf{k}_{2}\right)
+C_{0}
\end{equation*}%
\ holds in these neighborhoods where $C_{0},\delta _{1},\delta _{2}$ are
constants and $\delta _{1}\delta _{2}\neq 0$. Then all three functions $%
f_{1}\left( \mathbf{k}\right) $, $f_{2}\left( \mathbf{k}\right) $, $%
f_{3}\left( \mathbf{k}\right) $ are linear in neighborhoods of $\mathbf{k}%
_{\ast 1},$ $\mathbf{k}_{\ast 2}$, $\mathbf{k}_{\ast 1}+\mathbf{k}_{\ast 2}$
respectively.
\end{proposition}

\section{Reduction to a standard framework}

Many well known nonlinear evolutionary equations and systems can be easily
reduced to the framework of (\ref{difeqintr}), (\ref{difeqfou}) involving
two small parameters $\varrho $ and $\beta $ and characterized by the
following properties: (i) the linear part is $\mathbf{L}$ has large factor $%
\frac{1}{\varrho }$ before it; (ii) the nonlinearity $\mathbf{F}\left( 
\mathbf{U}\right) $ is independent of $\varrho ,\beta $ or depends on $%
\varrho $ regularly; (iii) the initial data depend on $\beta $ so that they
do not vanish as $\beta \rightarrow 0$; (iv) the solutions are considered on
the time interval $0\leq \tau \leq \tau _{\ast }$ where $\tau _{\ast }>0$
does not depend on $\varrho ,\beta $. Notice that solutions to (\ref%
{difeqintr}), (\ref{difeqfou}) under the above conditions exhibit nonlinear
effects uniformly with respect to small $\varrho ,\beta $\ on the time
interval $0\leq \tau \leq \tau _{\ast }$.

There are important classes of problems which can be readily reduced to the
framework of (\ref{difeqintr}), (\ref{difeqfou}) by a simple rescaling.

\subparagraph{Systems with a small factor before the nonlinearity.}

Consider a problem of the form 
\begin{equation}
\partial _{t}\mathbf{v}=-\mathrm{i}\mathbf{Lv}+\alpha \mathbf{f}\left( 
\mathbf{v}\right) ,\ \left. \mathbf{v}\right\vert _{t=0}=\mathbf{h},\
0<\alpha \ll 1,  \label{smallnon}
\end{equation}%
where initial data are bounded uniformly in $\alpha $. \ Such problems are
reduced to (\ref{difeqintr}) by the time rescaling $\tau =t\alpha .\ $Note
that now $\varrho =\alpha $ \ and the finite time interval $0\leq \tau \leq
\tau _{\ast }$ corresponds to the long time interval $0\leq t\leq \tau
_{\ast }/\alpha $.

\subparagraph{Systems with small initial data on long time intervals.}

The equation here is 
\begin{gather}
\partial _{t}\mathbf{v}=-\mathrm{i}\mathbf{Lv}+\mathbf{f}_{0}\left( \mathbf{v%
}\right) ,\ \left. \mathbf{v}\right\vert _{t=0}=\alpha _{0}\mathbf{h},\
0<\alpha _{0}\ll 1,\text{ where}  \label{smallinit} \\
\mathbf{f}_{0}\left( \mathbf{v}\right) =\mathbf{f}_{0}^{\left( m\right)
}\left( \mathbf{v}\right) +\mathbf{f}_{0}^{\left( m+1\right) }\left( \mathbf{%
v}\right) +\ldots ,  \notag
\end{gather}%
and $\mathbf{f}^{\left( m\right) }\left( \mathbf{v}\right) $ is a \
homogeneous polynomial of degree $m\geq 2$. \ After the rescaling $\mathbf{v}%
=\alpha _{0}\mathbf{V}$ we obtain the following equation \ with a small
nonlinearity 
\begin{equation}
\partial _{t}\mathbf{V}=-\mathrm{i}\mathbf{LV}+\alpha _{0}^{m-1}\left[ 
\mathbf{f}_{0}^{\left( m\right) }\left( \mathbf{V}\right) +\alpha _{0}%
\mathbf{f}^{0\left( m+1\right) }\left( \mathbf{V}\right) +\ldots \right] ,\
\left. \mathbf{V}\right\vert _{t=0}=\mathbf{h},  \label{smin1}
\end{equation}%
which is of the form of (\ref{smallnon}) with $\alpha =\alpha _{0}^{m-1}$.
Note that nonlinearities $\mathbf{f}$ in (\ref{smallnon}) which are obtained
from problems with small initial data and regular nonlinearities $\mathbf{f}%
_{0}\left( \mathbf{v}\right) $ have a special form. Namely, they are almost
homogeneous, $\mathbf{f}\left( \mathbf{V}\right) =\mathbf{f}_{0}^{\left(
m\right) }\left( \mathbf{V}\right) +\alpha \left[ \ldots \right] \ \ $with
leading term $\mathbf{f}_{0}^{\left( m\right) }\left( \mathbf{V}\right) $.
Introducing the slow time variable $\tau =t\alpha _{0}^{m-1}$ we get from
the above an equation of the form (\ref{difeqintr}), namely%
\begin{equation}
\partial _{\tau }\mathbf{V}=-\frac{\mathrm{i}}{\alpha _{0}^{m-1}}\mathbf{LV}+%
\left[ \mathbf{f}^{\left( m\right) }\left( \mathbf{V}\right) +\alpha _{0}%
\mathbf{f}^{\left( m+1\right) }\left( \mathbf{V}\right) +\ldots \right] ,\
\left. \mathbf{V}\right\vert _{t=0}=\mathbf{h},  \label{smin2}
\end{equation}%
where the nonlinearity does not vanish as $\alpha _{0}\rightarrow 0$. In
this case $\varrho =\alpha _{0}^{m-1}$ and the finite time interval $0\leq
\tau \leq \tau _{\ast }$ corresponds to the long time interval $0\leq t\leq 
\frac{\tau _{\ast }}{\alpha _{0}^{m-1}}$ \ with small $\alpha _{0}\ll 1$.
Note that Corollary \ref{Corollary Fofrho} for $\varrho $-dependent
nonlinearities can be applied to this case. This allows, in particular, to
apply results of this paper to Sine-Gordon equation where $\mathbf{f}%
_{0}\left( v\right) =\sin v$.

\subparagraph{High-frequency carrier waves.}

Sometimes high spatial frequency of carrier waves in the initial wavepackets
after a rescaling creates a large parameter $\frac{1}{\varrho }$ at the
linear part. For example, Nonlinear Schrodinger equation 
\begin{equation}
\partial _{\tau }U=-\mathrm{i}\partial _{x}^{2}U+\mathrm{i}\alpha \left\vert
U\right\vert ^{2}U,\ \left. U\right\vert _{\tau =0}=h_{1}\left( \beta
x\right) e^{\mathrm{i}Mk_{\ast 1}x}+h_{2}\left( \beta x\right) e^{\mathrm{i}%
Mk_{\ast 2}x}+c.c.,  \label{NLSh}
\end{equation}%
where $c.c.$ stands for complex conjugate of the prior term, and $M\gg 1$ is
a large parameter, can be recast in the form (\ref{difeqintr}). Indeed,
changing variables $y=Mx$ in the above equation we obtain 
\begin{equation*}
\partial _{\tau }U=-\mathrm{i}\frac{1}{\varrho }\partial _{y}^{2}U+\mathrm{i}%
\alpha \left\vert U\right\vert ^{2}U,\ \left. U\right\vert _{\tau
=0}=h_{1}\left( \beta _{1}y\right) e^{\mathrm{i}k_{\ast 1}y}+h_{2}\left(
\beta _{1}y\right) e^{\mathrm{i}k_{\ast 2}y}+c.c.,
\end{equation*}%
where $\beta _{1}=\frac{\beta }{M}\ll 1$, $\varrho =\frac{1}{M^{2}}\ll 1$.
Note that though the nonlinearity$\left\vert U\right\vert ^{2}U$ in (\ref%
{NLSh}) is not complex homogeneous, it can be considered as a restriction of
a system with a complex homogeneous nonlinearity as (\ref{cansyst1}) is a
restriction of (\ref{cansystext}).

\subparagraph{ First order hyperbolic equations and systems.}

Consider now the system (\ref{toy1}), (\ref{toy2}) for which the symmetry (%
\ref{invsym}) does not hold. \ The system can be put into the standard
framework by formally adding two more equations 
\begin{gather}
\partial _{\tau }w_{1}=\frac{c_{1}}{\varrho }\partial _{x}w_{1}+F_{1}\left(
w_{1},w_{2}\right) ,\ \partial _{\tau }w_{2}=\frac{c_{2}}{\varrho }\partial
_{x}w_{2}+F_{2}\left( w_{1},w_{2}\right) ,  \label{toyext} \\
\left. w_{1}\right\vert _{\tau =0}=0,\ \left. w_{2}\right\vert _{\tau =0}=0,
\notag
\end{gather}%
which have only trivial solution $w_{1}=w_{2}=0$ not affecting the solutions
to the original system (\ref{toy1}), (\ref{toy2}). The extended system has
the linear part with two-band dispersion relations $\omega _{1,\zeta }\left(
k\right) =c_{1}\zeta \left\vert k\right\vert $, $\omega _{2,\zeta }\left(
k\right) =c_{2}\zeta \left\vert k\right\vert $ ,$\zeta =\pm ,$ satisfying
evidently (\ref{invsym}).

\section{Integrated evolution equation}

Using the variation of constants formula we recast the modal evolution
equation (\ref{difeqfou}) into the following equivalent integral form 
\begin{equation}
\mathbf{\hat{U}}\left( \mathbf{k},\tau \right) =\int_{0}^{\tau }\mathrm{e}^{%
\frac{-\mathrm{i}\left( \tau -\tau ^{\prime }\right) }{\varrho }\mathbf{L}%
\left( \mathbf{k}\right) }\hat{F}\left( \mathbf{\hat{U}}\right) \left( 
\mathbf{k},\tau \right) \,\mathrm{d}\tau ^{\prime }+\mathrm{e}^{\frac{-%
\mathrm{i}\zeta \tau }{\varrho }\mathbf{L}\left( \mathbf{k}\right) }\mathbf{%
\hat{h}}\left( \mathbf{k}\right) ,\ \tau \geq 0.  \label{varc}
\end{equation}%
Then we factor $\mathbf{\hat{U}}\left( \mathbf{k},\tau \right) $ into the
slow variable $\mathbf{\hat{u}}\left( \mathbf{k},\tau \right) $ and the fast
oscillatory term as in (\ref{Uu0}), namely 
\begin{equation}
\mathbf{\hat{U}}\left( \mathbf{k},\tau \right) =\mathrm{e}^{-\frac{\mathrm{i}%
\tau }{\varrho }\mathbf{L}\left( \mathbf{k}\right) }\mathbf{\hat{u}}\left( 
\mathbf{k},\tau \right) ,\ \mathbf{\hat{U}}_{n,\zeta }\left( \mathbf{k},\tau
\right) =\mathbf{\hat{u}}_{n,\zeta }\left( \mathbf{k},\tau \right) \mathrm{e}%
^{-\frac{\mathrm{i}\tau }{\varrho }\zeta \omega _{n}\left( \mathbf{k}\right)
},  \label{Uu}
\end{equation}%
where $\mathbf{\hat{u}}_{n,\zeta }\left( \mathbf{k},\tau \right) $ are the
modal coefficients of $\mathbf{\hat{u}}\left( \mathbf{k},\tau \right) $ as
in (\ref{Uboldj}). Notice that $\mathbf{\hat{u}}_{n,\zeta }\left( \mathbf{k}%
,\tau \right) $ in (\ref{Uu}) may depend on $\varrho $ and (\ref{Uu}) is
just a change of variables and not an assumption. Consequently we obtain the
following \emph{integrated evolution equation} for $\mathbf{\hat{u}}=\mathbf{%
\hat{u}}\left( \mathbf{k},\tau \right) $, $\tau \geq 0$, 
\begin{gather}
\mathbf{\hat{u}}\left( \mathbf{k},\tau \right) =\mathcal{F}\left( \mathbf{%
\hat{u}}\right) \left( \mathbf{k},\tau \right) +\mathbf{\hat{h}}\left( 
\mathbf{k}\right) ,\ \mathcal{F}\left( \mathbf{\hat{u}}\right)
=\sum\nolimits_{m\in \mathfrak{M}_{F}}\mathcal{F}^{\left( m\right) }\left( 
\mathbf{\hat{u}}^{m}\left( \mathbf{k},\tau \right) \right) ,  \label{varcu}
\\
\mathcal{F}^{\left( m\right) }\left( \mathbf{\hat{u}}^{m}\right) \left( 
\mathbf{k},\tau \right) =\int_{0}^{\tau }\mathrm{e}^{\frac{\mathrm{i}\tau
^{\prime }}{\varrho }\mathbf{L}\left( \mathbf{k}\right) }\hat{F}_{m}\left(
\left( \mathrm{e}^{\frac{-\mathrm{i}\tau ^{\prime }}{\varrho }\mathbf{L}%
\left( \cdot \right) }\mathbf{\hat{u}}\right) ^{m}\right) \left( \mathbf{k}%
,\tau ^{\prime }\right) \,\mathrm{d}\tau ^{\prime },  \label{Fu}
\end{gather}%
where $\hat{F}_{m}$ are defined by (\ref{Fseries}) and (\ref{Fmintr}) in
terms of the susceptibilities $\chi ^{\left( m\right) }$, and $\mathcal{F}%
^{\left( m\right) }$ are bounded as in the following lemma.

\begin{lemma}[boundness of multilinear operators]
\label{Lemma bound}$\mathcal{F}^{\left( m\right) }$ defined by (\ref{Fmintr}%
), (\ref{Fu}) is bounded operator from $E=C\left( \left[ 0,\tau _{\ast }%
\right] ,L^{1}\right) $ into $C^{1}\left( \left[ 0,\tau _{\ast }\right]
,L^{1}\right) $ satisfying 
\begin{equation}
\left\Vert \mathcal{F}^{\left( m\right) }\left( \mathbf{\hat{u}}_{1}\ldots 
\mathbf{\hat{u}}_{m}\right) \right\Vert _{E}\leq \tau _{\ast }\left\Vert
\chi ^{\left( m\right) }\right\Vert \dprod\nolimits_{j=1}^{m}\left\Vert 
\mathbf{\hat{u}}_{j}\right\Vert _{E},  \label{dtf}
\end{equation}%
\begin{equation}
\left\Vert \partial _{\tau }\mathcal{F}^{\left( m\right) }\left( \mathbf{%
\hat{u}}_{1}\ldots \mathbf{\hat{u}}_{m}\right) \right\Vert _{E}\leq
\left\Vert \chi ^{\left( m\right) }\right\Vert \dprod\nolimits_{j}\left\Vert 
\mathbf{\hat{u}}_{j}\right\Vert _{E}.  \label{dtf1}
\end{equation}
\end{lemma}

\begin{proof}
Notice that since $\mathbf{L}\left( \mathbf{\mathbf{k}}\right) $ is
Hermitian, $\left\Vert \exp \left\{ -\mathrm{i}\mathbf{L}\left( \mathbf{%
\mathbf{k}}\right) \frac{\tau _{1}}{\varrho }\right\} \right\Vert =1$. Using
the Young inequality 
\begin{equation}
\left\Vert \mathbf{\hat{u}}\ast \mathbf{\hat{v}}\right\Vert _{L^{1}}\leq
\left\Vert \mathbf{\hat{u}}\right\Vert _{L^{1}}\left\Vert \mathbf{\hat{v}}%
\right\Vert _{L^{1}}.  \label{Yconv}
\end{equation}%
together with (\ref{Fmintr}), (\ref{Fu}) we obtain 
\begin{gather*}
\left\Vert \mathcal{F}^{\left( m\right) }\left( \mathbf{\hat{u}}_{1}\ldots 
\mathbf{\hat{u}}_{m}\right) \left( \mathbf{\cdot },\tau \right) \right\Vert
_{L^{1}}\leq \sup_{\ \mathbf{\mathbf{k}},\vec{k}}\left\vert \chi _{\
}^{\left( m\right) }\left( \mathbf{\mathbf{k}},\vec{k}\right) \right\vert \\
\int_{\mathbb{R}^{d}}\int_{0}^{\tau }\int_{\mathbb{D}_{m}}\left\vert \mathbf{%
\hat{u}}_{1}\left( \mathbf{k}^{\prime }\right) \right\vert \ldots \left\vert 
\mathbf{\hat{u}}_{m}\left( \mathbf{k}^{\left( m\right) }\left( \mathbf{k},%
\vec{k}\right) \right) \right\vert \,\mathrm{d}\mathbf{k}^{\prime }\ldots 
\mathrm{d}\mathbf{k}^{\left( m-1\right) }\mathrm{d}\tau _{1}\mathrm{d}%
\mathbf{k}\leq \\
\left\Vert \chi ^{\left( m\right) }\right\Vert \int_{0}^{\tau }\left\Vert 
\mathbf{\hat{u}}_{1}\left( \tau _{1}\right) \right\Vert _{L^{1}}\ldots
\left\Vert \mathbf{\hat{u}}_{m}\left( \tau _{1}\right) \right\Vert _{L^{1}}\,%
\mathrm{d}\tau _{1}\leq \tau _{\ast }\left\Vert \chi ^{\left( m\right)
}\right\Vert \left\Vert \mathbf{\hat{u}}_{1}\right\Vert _{E}\ldots
\left\Vert \mathbf{\hat{u}}_{m}\right\Vert _{E}.
\end{gather*}%
proving (\ref{dtf}). Similarly we prove (\ref{dtf1}) by%
\begin{gather*}
\left\Vert \partial _{\tau }\mathcal{F}^{\left( m\right) }\left( \mathbf{%
\hat{u}}_{1}\ldots \mathbf{\hat{u}}_{m}\right) \left( \mathbf{\cdot },\tau
\right) \right\Vert _{L^{1}}\leq \left\Vert \chi ^{\left( m\right)
}\right\Vert \\
\int_{\mathbb{R}^{d}}\int_{\mathbb{D}_{m}}\left\vert \mathbf{\hat{u}}%
_{1}\left( \mathbf{k}^{\prime }\right) \right\vert \ldots \left\vert \mathbf{%
\hat{u}}_{m}\left( \mathbf{k}^{\left( m\right) }\left( \mathbf{k},\vec{k}%
\right) \right) \right\vert \,\mathrm{d}\mathbf{k}^{\prime }\ldots \,\mathrm{%
d}\mathbf{k}^{\left( m-1\right) }\mathrm{d}\mathbf{k}\,\leq \left\Vert \chi
^{\left( m\right) }\right\Vert \left\Vert \mathbf{\hat{u}}_{1}\right\Vert
_{E}\ldots \left\Vert \mathbf{\hat{u}}_{m}\right\Vert _{E}.
\end{gather*}
\end{proof}

The equation (\ref{varcu}) can be recast as the following abstract equation
in a Banach space 
\begin{equation}
\mathbf{\hat{u}}=\mathcal{F}\left( \mathbf{\hat{u}}\right) +\mathbf{\hat{h}}%
,\ \mathbf{\hat{u}},\mathbf{\hat{h}}\in E,  \label{eqF}
\end{equation}%
and it readily follows from Lemma \ref{Lemma bound} that $\mathcal{F}\left( 
\mathbf{\hat{u}}\right) $ has the following properties.

\begin{lemma}
\label{Lemma Flip} The operator\ $\mathcal{F}\left( \mathbf{\hat{u}}\right) $
defined by (\ref{varcu})-(\ref{Fu}) satisfies the Lipschitz condition 
\begin{equation}
\left\Vert \mathcal{F}\left( \mathbf{\hat{u}}_{1}\right) -\mathcal{F}\left( 
\mathbf{\hat{u}}_{2}\right) \right\Vert _{E}\leq \tau _{\ast
}C_{F}\left\Vert \mathbf{\hat{u}}_{1}-\mathbf{\hat{u}}_{2}\right\Vert _{E}
\label{Flip}
\end{equation}%
where $C_{F}\leq C_{\chi }m_{F}^{2}\left( 4R\right) ^{m_{F}-1}$ if $%
\left\Vert \mathbf{\mathbf{\hat{u}}}_{1}\right\Vert _{E},\left\Vert \mathbf{%
\mathbf{\hat{u}}}_{2}\right\Vert _{E}\leq 2R$, with $C_{\chi }$ as in (\ref%
{chiCR}).
\end{lemma}

We also will use the following form of the contraction principle.

\begin{lemma}[Contraction principle]
\label{Lemma contr}Consider equation \ 
\begin{equation}
\mathbf{x}=\mathcal{F}\left( \mathbf{x}\right) +\mathbf{h},\ \mathbf{x},%
\mathbf{h}\in B,  \label{eqFa}
\end{equation}%
where $B$ is a Banach space, $\mathcal{F}$ is an operator in $B$. Suppose
that for some constants $R_{0}>0$ and $0<q<1$ we have%
\begin{eqnarray}
\left\Vert \mathbf{h}\right\Vert &\leq &R_{0},\ \left\Vert \mathcal{F}\left( 
\mathbf{x}\right) \right\Vert \leq R_{0}\text{ if }\left\Vert \mathbf{x}%
\right\Vert \leq 2R_{0},  \label{Fle} \\
\left\Vert \mathcal{F}\left( \mathbf{x}_{1}\right) -\mathcal{F}\left( 
\mathbf{x}_{2}\right) \right\Vert &\leq &q\left\Vert \mathbf{x}_{1}-\mathbf{x%
}_{2}\right\Vert \text{ if }\left\Vert \mathbf{x}_{1}\right\Vert ,\left\Vert 
\mathbf{x}_{2}\right\Vert \leq 2R_{0}.  \label{Flipa}
\end{eqnarray}%
Then there exists a unique solution $\mathbf{x}$ to the equation (\ref{eqFa}%
) such that $\left\Vert \mathbf{x}\right\Vert \leq 2R_{0}$. Let $\left\Vert 
\mathbf{h}_{1}\right\Vert ,\left\Vert \mathbf{h}_{2}\right\Vert \leq R_{0}$
then the two corresponding solutions $\mathbf{x}_{1},\mathbf{x}_{2}$ satisfy 
\begin{equation}
\left\Vert \mathbf{x}_{1}\right\Vert ,\left\Vert \mathbf{x}_{2}\right\Vert
\leq 2R_{0},\ \left\Vert \mathbf{x}_{1}-\mathbf{x}_{2}\right\Vert \leq
\left( 1-q\right) ^{-1}\left\Vert \mathbf{h}_{1}-\mathbf{h}_{2}\right\Vert .
\label{iminu0}
\end{equation}%
Let $\mathbf{x}_{1},\mathbf{x}_{2}$ be the two solutions of correspondingly
two equations of the form (\ref{eqFa}) with $\mathcal{F}_{1}$, $\mathbf{h}%
_{1}\mathbf{\ }$and\ $\mathcal{F}_{2}$, $\mathbf{h}_{2}$. Assume that that $%
\mathcal{F}_{1}\left( \mathbf{u}\right) $ satisfies (\ref{Fle}), (\ref{Flipa}%
) with a Lipschitz constant $q<1$ and that $\left\Vert \mathcal{F}_{1}\left( 
\mathbf{x}\right) -\mathcal{F}_{2}\left( \mathbf{x}\right) \right\Vert \leq
\delta $ for $\left\Vert \mathbf{x}\right\Vert \leq 2R_{0}$. Then 
\begin{equation}
\left\Vert \mathbf{x}_{1}-\mathbf{x}_{2}\right\Vert \leq \left( 1-q\right)
^{-1}\left( \delta +\left\Vert \mathbf{h}_{1}-\mathbf{h}_{2}\right\Vert
\right) .  \label{uminu}
\end{equation}
\end{lemma}

Lemma \ref{Lemma Flip} and the contraction principle as in Lemma \ref{Lemma
contr} imply the following existence and uniqueness theorem.

\begin{theorem}
\label{Theorem exist}Let $\left\Vert \mathbf{h}\right\Vert _{E}\leq R$, let $%
\tau _{\ast }<1/C_{F}$ where $C_{F}$ is a constant from Lemma \ref{Lemma
Flip}. Then equation (\ref{varcu}) has a solution $\mathbf{\hat{u}}\in
E=C\left( \left[ 0,\tau _{\ast }\right] ,L^{1}\right) $ which satisfies $%
\left\Vert \mathbf{\hat{u}}\right\Vert _{E}\leq 2R$, and such a solution is
unique.
\end{theorem}

The following existence and uniqueness theorem follows from Theorem \ref%
{Theorem exist}.

\begin{theorem}
\label{Theorem Existence1}\ Let (\ref{difeqfou}) satisfy (\ref{chiCR}) and $%
\mathbf{\hat{h}}\in L^{1}\left( \mathbb{R}^{d}\right) $,$\left\Vert \mathbf{%
\hat{h}}\right\Vert _{L^{1}}\leq R$. Then there exists a unique solution to
the modal evolution equation (\ref{difeqfou}) in the functional space $%
C^{1}\left( \left[ 0,\tau _{\ast }\right] ,L^{1}\right) $. The number $\tau
_{\ast }$ depends on $R$ and $C_{\chi }$.
\end{theorem}

Using the inequality (\ref{Linf}) and applying the inverse Fourier transform
we readily obtain the existence of an $F-$solution of (\ref{difeqintr}) in $%
C^{1}\left( \left[ 0,\tau _{\ast }\right] ,L^{\infty }\left( \mathbb{R}%
^{d}\right) \right) $ from the existence of the solution of equation (\ref%
{difeqfou}) in $C^{1}\left( \left[ 0,\tau _{\ast }\right] ,L^{1}\right) $.
The existence of $F$-solutions in spaces of spatially smooth functions can
be derived by replacing Lemma \ref{Lemma bound} with an estimate similar to
the one in Lemma \ref{Lemma fscallip}.

Let us recast now the system (\ref{varcu})-(\ref{Fu}) into modal components
using the projections $\Pi _{n,\zeta }\left( \mathbf{\mathbf{k}}\right) $ as
in (\ref{Pin}). The first step to introduce elementary modal
susceptibilities $\chi _{n,\zeta ,\vec{\xi}}^{\left( m\right) }$ having
one-dimensional range in $\mathbb{C}^{2J}$ and vanishing if one of its
arguments $\mathbf{\hat{u}}_{j}$ belongs to a $\left( 2J-1\right) $%
-dimensional linear subspace in $\mathbb{C}^{2J}$ ($j$-th null-space \ of $%
\chi _{n,\zeta ,\vec{\xi}}^{\left( m\right) }$ ). For example, in the linear
case $m=1$ when $\chi ^{\left( 1\right) }$ acts in $\mathbb{C}^{2J}$ \ and
is presented in the standard orthonormal basis $\left\{ \mathbf{e}_{n,\zeta
}\right\} $ in $\mathbb{C}^{2J}$ by a $2J\times 2J$ matrix with elements $%
a_{\xi ,\xi ^{\prime }}^{\left( 1\right) }=a_{n,\zeta ,n^{\prime },\zeta
^{\prime }}^{\left( 1\right) }$, where index $\xi =n,\zeta $ takes $2J$
values, the action of elementary susceptibility $\chi _{n,\zeta ,n^{\prime
},\zeta ^{\prime }}^{\left( 1\right) }$ on a vector $\mathbf{v}\in \mathbb{C}%
^{2J}$ is given by the formula $\chi _{n,\zeta ,n^{\prime },\zeta ^{\prime
}}^{\left( 1\right) }\mathbf{v}=a_{n,\zeta ,n^{\prime },\zeta ^{\prime
}}^{\left( 1\right) }\left( \mathbf{v\cdot e}_{n^{\prime },\zeta ^{\prime
}}\right) \mathbf{e}_{n,\zeta }$ \ where $\left\{ \mathbf{e}_{n,\zeta
}\right\} $ is the standard orthonormal basis in $\mathbb{C}^{2J}$. \
Obviously $\chi _{n,\zeta ,n^{\prime },\zeta ^{\prime }}^{\left( 1\right) }%
\mathbf{v}=\Pi _{n,\zeta }\chi ^{\left( 1\right) }\Pi _{n^{\prime },\zeta
^{\prime }}\mathbf{v}$ \ and $\chi ^{\left( 1\right) }\mathbf{v=}%
\sum_{n,\zeta ,n^{\prime },\zeta ^{\prime }}\chi _{n,\zeta ,n^{\prime
},\zeta ^{\prime }}^{\left( 1\right) }\mathbf{v}$. The general definition
follows.

\begin{definition}[elementary susceptibilities]
\label{Definition elemsus}Let$\ $%
\begin{equation}
\vec{\xi}=\left( \vec{n},\vec{\zeta}\right) \in \left\{ 1,\ldots ,J\right\}
^{m}\times \left\{ -1,1\right\} ^{m}=\Xi ^{m},\left( n,\zeta \right) \in \Xi
\label{Xi}
\end{equation}%
and $\chi _{\ }^{\left( m\right) }\left( \mathbf{\mathbf{k}},\vec{k}\right) %
\left[ \mathbf{\hat{u}}_{1}\left( \mathbf{k}^{\prime }\right) ,\ldots ,%
\mathbf{\hat{u}}_{m}\left( \mathbf{k}^{\left( m\right) }\right) \right] $ be 
$\ m$-linear symmetric\ tensor (susceptibility) as in (\ref{Fmintr}). We
introduce \emph{elementary} \emph{susceptibilities} $\chi _{n,\zeta ,\vec{\xi%
}}^{\left( m\right) }\left( \mathbf{\mathbf{k}},\vec{k}\right) :\left( 
\mathbb{C}^{2J}\right) ^{m}\rightarrow \mathbb{C}^{2J}$) as$\ m$-linear
tensors defined for almost all $\mathbf{\mathbf{k}},\vec{k}$ by the
following formula 
\begin{gather}
\chi _{n,\zeta ,\vec{\xi}}^{\left( m\right) }\left( \mathbf{\mathbf{k}},\vec{%
k}\right) \left[ \mathbf{\hat{u}}_{1}\left( \mathbf{k}^{\prime }\right)
,\ldots ,\mathbf{\hat{u}}_{m}\left( \mathbf{k}^{\left( m\right) }\right) %
\right] =\chi _{n,\zeta ,\vec{n},\vec{\zeta}}^{\left( m\right) }\left( 
\mathbf{\mathbf{k}},\vec{k}\right) \left[ \mathbf{\hat{u}}_{1}\left( \mathbf{%
k}^{\prime }\right) ,\ldots ,\mathbf{\hat{u}}_{m}\left( \mathbf{k}^{\left(
m\right) }\right) \right] =  \label{chim} \\
\Pi _{n,\zeta }\left( \mathbf{\mathbf{k}}\right) \chi ^{\left( m\right)
}\left( \mathbf{\mathbf{k}},\vec{k}\right) \left[ \left( \Pi _{n_{1},\zeta
^{\prime }}\left( \mathbf{k}^{\prime }\right) \mathbf{\hat{u}}_{1}\left( 
\mathbf{k}^{\prime }\right) ,\ldots ,\Pi _{n_{m},\zeta ^{\left( m\right)
}}\left( \mathbf{k}^{\left( m\right) }\left( \mathbf{k},\vec{k}\right)
\right) \mathbf{\hat{u}}_{m}\left( \mathbf{k}^{\left( m\right) }\left( 
\mathbf{k},\vec{k}\right) \right) \right) \right] .  \notag
\end{gather}
\end{definition}

Then using (\ref{sumPi}) and the \emph{elementary} \emph{susceptibilities} (%
\ref{chim}) we get%
\begin{equation}
\chi ^{\left( m\right) }\left( \mathbf{\mathbf{k}},\vec{k}\right) \left[ 
\mathbf{\hat{u}}_{1}\left( \mathbf{k}^{\prime }\right) ,\ldots ,\mathbf{\hat{%
u}}_{m}\left( \mathbf{k}^{\left( m\right) }\right) \right]
=\sum\nolimits_{n,\zeta }\sum\nolimits_{\vec{\xi}}\chi _{n,\zeta ,\vec{\xi}%
}^{\left( m\right) }\left( \mathbf{\mathbf{k}},\vec{k}\right) \left[ \mathbf{%
\hat{u}}_{1}\left( \mathbf{k}^{\prime }\right) ,\ldots ,\mathbf{\hat{u}}%
_{m}\left( \mathbf{k}^{\left( m\right) }\right) \right] .  \label{chisumn}
\end{equation}%
Consequently the modal components $\mathcal{F}_{n,\zeta ,\vec{\xi}}^{\left(
m\right) }$\ of the operators $\mathcal{F}^{\left( m\right) }$ in (\ref{Fu})
are $m$-linear \emph{oscillatory integral operators} defined in terms of the
elementary susceptibilities (\ref{chisumn}) as follows.

\begin{definition}[interaction phase]
Using notations from (\ref{Fmintr}) we introduce for $\vec{\xi}=\left( \vec{n%
},\vec{\zeta}\right) \in \Xi ^{m}$ operator 
\begin{gather}
\mathcal{F}_{n,\zeta ,\vec{\xi}}^{\left( m\right) }\left( \mathbf{\tilde{u}}%
_{1}\ldots \mathbf{\tilde{u}}_{m}\right) \left( \mathbf{k},\tau \right)
=\int_{0}^{\tau }\int_{\mathbb{D}_{m}}\exp \left\{ \mathrm{i}\phi _{n,\zeta ,%
\vec{\xi}}\left( \mathbf{\mathbf{k}},\vec{k}\right) \frac{\tau _{1}}{\varrho 
}\right\}  \label{Fm} \\
\chi _{n,\zeta ,\vec{\xi}}^{\left( m\right) }\left( \mathbf{\mathbf{k}},\vec{%
k}\right) \left[ \mathbf{\tilde{u}}_{1}\left( \mathbf{k}^{\prime },\tau
_{1}\right) ,\ldots ,\mathbf{\tilde{u}}_{m}\left( \mathbf{k}^{\left(
m\right) }\left( \mathbf{\mathbf{k}},\vec{k}\right) ,\tau _{1}\right) \right]
\mathrm{\tilde{d}}^{\left( m-1\right) d}\vec{k}\mathrm{d}\tau _{1},  \notag
\end{gather}%
with the \emph{interaction phase function }$\phi $ defined by 
\begin{gather}
\phi _{n,\zeta ,\vec{\xi}}\left( \mathbf{\mathbf{k}},\vec{k}\right) =\phi
_{n,\zeta ,\vec{n},\vec{\zeta}}\left( \mathbf{\mathbf{k}},\vec{k}\right)
\label{phim} \\
=\zeta \omega _{n}\left( \zeta \mathbf{k}\right) -\zeta ^{\prime }\omega
_{n_{1}}\left( \zeta ^{\prime }\mathbf{k}^{\prime }\right) -\ldots -\zeta
^{\left( m\right) }\omega _{n_{m}}\left( \zeta ^{\left( m\right) }\mathbf{k}%
^{\left( m\right) }\right) ,\ \mathbf{k}^{\left( m\right) }=\mathbf{k}%
^{\left( m\right) }\left( \mathbf{\mathbf{k}},\vec{k}\right) .  \notag
\end{gather}
\end{definition}

Using $\mathcal{F}_{n,\zeta ,\vec{\xi}}^{\left( m\right) }$ in (\ref{Fm}) we
recast $\mathcal{F}^{\left( m\right) }\left( \mathbf{u}^{m}\right) $ in the
system (\ref{varcu})-(\ref{Fu}) as 
\begin{equation}
\mathcal{F}^{\left( m\right) }\left[ \mathbf{\hat{u}}_{1}\ldots ,\mathbf{%
\hat{u}}_{m}\right] \left( \mathbf{k},\tau \right) =\sum\nolimits_{n,\zeta ,%
\vec{\xi}}\mathcal{F}_{n,\zeta ,\vec{\xi}}^{\left( m\right) }\left[ \mathbf{%
\hat{u}}_{1}\ldots \mathbf{\hat{u}}_{m}\right] \left( \mathbf{k},\tau
\right) ,  \label{Tplusmin}
\end{equation}%
yielding the following system for the modal components $\mathbf{\hat{u}}%
_{n,\zeta }\left( \mathbf{k},\tau \right) $ as in (\ref{Pin}) 
\begin{equation}
\mathbf{\hat{u}}_{n,\zeta }\left( \mathbf{k},\tau \right)
=\sum\nolimits_{m\in \mathfrak{M}_{F}}\sum\nolimits_{\vec{\xi}\in \Xi ^{m}}%
\mathcal{F}_{n,\zeta ,\vec{\xi}}^{\left( m\right) }\left( \mathbf{\hat{u}}%
^{m}\right) \left( \mathbf{k},\tau \right) +\mathbf{\hat{h}}_{n,\zeta
}\left( \mathbf{k}\right) ,\ \left( n,\zeta \right) \in \Xi .  \label{equfa}
\end{equation}

\section{Wavepacket interaction system}

The wavepacket preservation property of the nonlinear evolutionary system in
any of its forms (\ref{difeqintr}), (\ref{difeqfou}), (\ref{varcu}), (\ref%
{eqF}), (\ref{equfa}) is not easy to see directly. It turns out though that
dynamics of wavepackets is well described by a system in a larger space $%
E^{2N}$ based on the original equation (\ref{varcu}) in the space $E$. We
call it \emph{wavepacket interaction system}, which is useful in three ways:
(i) the wavepacket preservation is quite easy to see and verify; (ii) it can
be used to prove the wavepacket preservation for the original nonlinear
problem; (iii) it can be used to study more subtle properties of the
original problem, such as NLS approximation. We start with the system (\ref%
{varcu}) where $\mathbf{\hat{h}}\left( \mathbf{k}\right) $ is a
multiwavepacket with a given $nk$-spectrum $S=\left\{ \left( \mathbf{k}%
_{\ast l},n_{l}\right) ,\ l=1,...,N\right\} $ as in (\ref{P0}) and $k$%
-spectrum $K_{S}=\left\{ \mathbf{k}_{\ast i},\ i=1,...,\left\vert
K_{S}\right\vert \right\} $ as in (\ref{K0}).

When constructing the wavepacket interaction system it is convenient to have
relevant functions to be explicitly localized about the $k$-spectrum $K_{S}$
of the initial data. We implement that by making up the following cutoff
functions based on (\ref{j0}), (\ref{Psik}) 
\begin{equation}
\Psi _{i,\vartheta }\left( \mathbf{k}\right) =\Psi \left( \mathbf{k}%
,\vartheta \mathbf{k}_{\ast i}\right) =\Psi \left( \beta ^{-\left(
1-\epsilon \right) }\left( \mathbf{k}-\vartheta \mathbf{k}_{\ast i}\right)
\right) ,\ \mathbf{k}_{\ast i}\in K_{S},\ i=1,\ldots ,\left\vert
K_{S}\right\vert ,\ \vartheta =\pm  \label{Psilz}
\end{equation}%
with $\epsilon $ as in Definition \ref{dwavepack} and $\beta >0$ small
enough to satisfy 
\begin{equation}
\beta ^{1/2}\leq \pi _{0},\text{ where }\pi _{0}=\pi _{0}\left( S\right) <%
\frac{1}{2}\min_{\mathbf{k}_{\ast i}\in K_{S}}\limfunc{dist}\left\{ \mathbf{k%
}_{\ast i},\sigma _{bc}\right\} .  \label{kpi0}
\end{equation}%
In what follows we use notations from (\ref{setlam}) and 
\begin{equation}
\vec{l}=\left( l_{1},...,l_{m}\right) \in \left\{ 1,...,N\right\} ^{m},\ 
\vec{\vartheta}=\left( \vartheta ^{\prime },...,\vartheta ^{\left( m\right)
}\right) \in \left\{ -1,1\right\} ^{m},\ \vec{\lambda}=\left( \vec{l},\vec{%
\vartheta}\right) \in \Lambda ^{m},  \label{narrow}
\end{equation}%
\begin{gather}
\vec{n}=\left( n_{1},\ldots ,n_{m}\right) \in \left\{ 1,...,J\right\} ^{m},\ 
\vec{\zeta}\in \left\{ -1,1\right\} ^{m},  \label{zetaar} \\
\vec{\xi}=\left( \vec{n},\vec{\zeta}\right) \in \Xi ^{m}\ ,\vec{k}=\left( 
\mathbf{\mathbf{k}}^{\prime },\ldots ,\mathbf{\mathbf{k}}^{\left( m\right)
}\right) \in \mathbb{R}^{m},\text{ where }\Xi ^{m}\text{ as in (\ref{Xi})}. 
\notag
\end{gather}%
Based on the above we introduce now the \emph{wavepacket interaction system} 
\begin{gather}
\mathbf{\hat{w}}_{l,\vartheta }\left( \mathbf{\cdot }\right) =\Psi \left( 
\mathbf{\cdot },\vartheta \mathbf{k}_{\ast i_{l}}\right) \Pi
_{n_{l},\vartheta }\left( \mathbf{\cdot }\right) \mathcal{F}\left(
\sum\nolimits_{\left( l^{\prime },\vartheta ^{\prime }\right) \in \Lambda }%
\mathbf{\hat{w}}_{l^{\prime },\vartheta ^{\prime }}\right) +\Psi \left( 
\mathbf{\cdot },\vartheta \mathbf{k}_{\ast i_{l}}\right) \Pi
_{n_{l},\vartheta }\left( \mathbf{\cdot }\right) \mathbf{\hat{h}},\left(
l,\vartheta \right) \in \Lambda ,  \label{sysloc} \\
\mathbf{\vec{w}}=\left( \mathbf{\hat{w}}_{1,+},\mathbf{\hat{w}}_{1,-},...,%
\mathbf{\hat{w}}_{N,+},\mathbf{\hat{w}}_{N,-}\right) \in E^{2N},\ \mathbf{%
\hat{w}}_{l,\vartheta }\in E,\left( l,\vartheta \right) \in \Lambda ,  \notag
\end{gather}%
with $\Psi \left( \mathbf{\cdot },\vartheta \mathbf{k}_{\ast i}\right) ,\Pi
_{n,\vartheta }$ being as in (\ref{Psilz}), (\ref{Pin}), $\mathcal{F}$
defined by (\ref{varcu}), and the norm in \ $E^{2N}$ defined based on (\ref%
{Elat}) by the formula 
\begin{equation*}
\left\Vert \mathbf{\vec{w}}\right\Vert _{E^{2N}}=\sum\nolimits_{l,\vartheta
}\left\Vert \mathbf{\hat{w}}_{l,\vartheta }\right\Vert _{E},\ E=C\left( %
\left[ 0,\tau _{\ast }\right] ,L^{1}\right) .
\end{equation*}%
The index $\left( l,\vartheta \right) $ which takes $2N$ values labels
equations and variables, the right-hand side of (\ref{sysloc}) is
well-defined for all $\mathbf{\vec{w}}\in E^{2N}$ and the equality (\ref%
{sysloc}) is understood as equality of elements of $E^{2N}$. We also use the
following concise form of the wave interaction system (\ref{sysloc}) 
\begin{gather}
\mathbf{\vec{w}}=\mathcal{F}_{_{\Psi }}\left( \mathbf{\vec{w}}\right) +%
\mathbf{\vec{h}}_{_{\Psi }},\text{ where}  \label{syslocF} \\
\mathbf{\vec{h}}_{_{\Psi }}=\left( \Psi _{i_{1},+}\Pi _{n_{1},+}\mathbf{\hat{%
h}},\Psi _{i_{1},-}\Pi _{n_{1},-}\mathbf{\hat{h}},...,\Psi _{i_{N},+}\Pi
_{n_{N},+}\mathbf{\hat{h}},\Psi _{i_{N},-}\Pi _{n_{N},-}\mathbf{\hat{h}}%
\right) \in E^{2N}.  \notag
\end{gather}%
The following lemma is analogous to Lemmas \ref{Lemma bound}, \ref{Lemma
Flip}.

\begin{lemma}
\label{Lemma Flippr} Polynomial operator $\mathcal{F}_{\Psi }\left( \mathbf{%
\vec{w}}\right) $ is bounded in $E^{2N}$, $\mathcal{F}_{\Psi }\left( \mathbf{%
0}\right) =\mathbf{0}$, and it satisfies Lipschitz condition 
\begin{equation}
\left\Vert \mathcal{F}_{\Psi }\left( \mathbf{\vec{w}}_{1}\right) -\mathcal{F}%
_{\Psi }\left( \mathbf{\vec{w}}_{2}\right) \right\Vert _{E^{2N}}\leq C\tau
_{\ast }\left\Vert \mathbf{\vec{w}}_{1}-\mathbf{\vec{w}}_{2}\right\Vert
_{E^{2N}},  \label{lipFpsi}
\end{equation}%
where $C$ depends only on $C_{\chi }$ as in (\ref{chiCR}), on the degree of $%
\mathcal{F}$ and on $\left\Vert \mathbf{\vec{w}}_{1}\right\Vert
_{E^{2N}}+\left\Vert \mathbf{\vec{w}}_{2}\right\Vert _{E^{2N}}$, and it does
not depend on $\beta $ and $\varrho $.
\end{lemma}

\begin{proof}
We consider every operator $\mathcal{F}_{n,\zeta ,\vec{\xi}}^{\left(
m\right) }\left( \mathbf{\vec{w}}\right) $ defined by (\ref{Fm}) and prove
its boundedness and the Lipschitz property as in Lemma \ref{Lemma bound}\
using the inequality $\left\vert \exp \left\{ \mathrm{i}\phi _{n,\zeta ,\vec{%
\xi}}\frac{\tau _{1}}{\varrho }\right\} \right\vert \leq 1$ and estimates (%
\ref{j0}), (\ref{chiCR}). Note that the integration in $\tau _{1}$ yields
the factor $\tau _{\ast }$ and consequent summation \ with respect to $%
n,\zeta ,\vec{\xi}$ yields (\ref{lipFpsi}).
\end{proof}

Lemma \ref{Lemma Flippr} \ and the contraction principle as in Lemma \ref%
{Lemma contr} yield the following statement.

\begin{theorem}
\label{Theorem existpr}Let $\left\Vert \mathbf{\vec{h}}_{_{\Psi
}}\right\Vert _{E^{2N}}\leq R.$ Then there exists $R_{1}>0$ and $\tau _{\ast
}>0$ such that equation (\ref{sysloc}) has a solution $\mathbf{\vec{w}}\in
E^{2N}$ which satisfies $\left\Vert \mathbf{\vec{w}}\right\Vert
_{E^{2N}}\leq R_{1}$ and such a solution is unique.
\end{theorem}

\begin{lemma}
\label{Lemma wiswave}Every function $\mathbf{\hat{w}}_{l,\zeta }\left( 
\mathbf{k},\tau \right) $ corresponding to the solution of \ (\ref{syslocF}) 
$\ $from $\ E^{2N}$ is a wavepacket with $nk$-pair $\left( \mathbf{k}_{\ast
l},n_{l}\right) $ with the degree of regularity which can be any $s>0$.
\end{lemma}

\begin{proof}
Note that according to (\ref{Psilz}) and (\ref{syslocF}) the function 
\begin{equation*}
\mathbf{\hat{w}}_{l,\vartheta }\left( \mathbf{k},\tau \right) =\Psi \left( 
\mathbf{k},\vartheta \mathbf{k}_{\ast i_{l}}\right) \Pi _{n_{l},\vartheta }%
\mathcal{F}\left( \mathbf{k},\tau \right) ,\ \left\Vert \mathcal{F}\left(
\tau \right) \right\Vert _{L^{1}}\leq C,\ 0\leq \tau \leq \tau _{\ast }
\end{equation*}%
\ involves the factor $\Psi _{l,\vartheta }\left( \mathbf{k}\right) =\Psi
\left( \beta ^{-\left( 1-\epsilon \right) }\left( \mathbf{k}-\vartheta 
\mathbf{k}_{\ast l}\right) \right) $ where $\epsilon $ is as in Definition %
\ref{dwavepack}. Hence,%
\begin{gather}
\Pi _{n,\vartheta ^{\prime }}\mathbf{\hat{w}}_{l,\vartheta }\left( \mathbf{k}%
,\tau \right) =0\text{ if }n\neq n_{l}\text{ or }\vartheta ^{\prime }\neq
\vartheta ,  \label{piweq0} \\
\mathbf{\hat{w}}_{l,\vartheta }\left( \mathbf{k},\tau \right) =\Psi \left( 
\mathbf{k},\vartheta \mathbf{k}_{\ast i_{l}}\right) \mathbf{\hat{w}}%
_{l,\vartheta }\left( \mathbf{k},\tau \right) ,\mathbf{\hat{w}}_{l,\vartheta
}\left( \mathbf{k},\tau \right) =0\text{ if\ }\left\vert \mathbf{k}%
-\vartheta \mathbf{k}_{\ast l}\right\vert \geq \beta ^{1-\epsilon },
\label{weq0}
\end{gather}%
and, consequently, Definition \ref{dwavepack} for $\mathbf{\hat{w}}%
_{l,\vartheta }$ is satisfied with $\hat{D}_{h}=0$ for any $s>0$ and $%
C^{\prime }=0$\ in (\ref{sourloc}).
\end{proof}

Now we would like to show that if $\mathbf{\hat{h}}$ is a multiwavepacket,
then the function 
\begin{equation}
\mathbf{\hat{w}}\left( \mathbf{k},\tau \right) =\sum\nolimits_{\left(
l,\vartheta \right) \in \Lambda }\mathbf{\hat{w}}_{l,\vartheta }\left( 
\mathbf{k},\tau \right) =\sum\nolimits_{\lambda \in \Lambda }\mathbf{\hat{w}}%
_{\lambda }\left( \mathbf{k},\tau \right)  \label{wkt}
\end{equation}%
is an approximate solution of equation (\ref{eqF}) (see notation (\ref%
{setlam})). To do that we introduce 
\begin{equation}
\Psi _{\infty }\left( \mathbf{k}\right) =1-\sum\nolimits_{\vartheta =\pm
}\sum\nolimits_{i=1}^{\left\vert K_{S}\right\vert }\Psi \left( \mathbf{k}%
,\vartheta \mathbf{k}_{\ast i}\right) =1-\sum\nolimits_{\vartheta =\pm
}\sum\nolimits_{\mathbf{k}_{\ast i}\in K_{S}}\Psi \left( \frac{\mathbf{k}%
-\vartheta \mathbf{k}_{\ast i}}{\beta ^{1-\epsilon }}\right) .
\label{Psiinf}
\end{equation}%
\ Expanding $m$-linear operator $\mathcal{F}^{\left( m\right) }\left( \left(
\sum_{l,\vartheta }\mathbf{\hat{w}}_{l,\vartheta }\right) ^{m}\right) $ and
using notations (\ref{setlam}), (\ref{lamprop}) we get%
\begin{gather}
\mathcal{F}^{\left( m\right) }\left( \left( \sum\nolimits_{l,\vartheta }%
\mathbf{\hat{w}}_{l,\vartheta }\right) ^{m}\right) =\sum\nolimits_{\vec{%
\lambda}\in \Lambda ^{m}}\mathcal{F}^{\left( m\right) }\left( \mathbf{\vec{w}%
}_{\vec{\lambda}}\right) ,\ \text{where}  \label{binom} \\
\mathbf{\vec{w}}_{\vec{\lambda}}=\mathbf{\hat{w}}_{\lambda _{1}}...\mathbf{%
\hat{w}}_{\lambda _{m}},\ \vec{\lambda}=\left( \lambda _{1},...,\lambda
_{m}\right) \in \Lambda ^{m}.  \label{wlam}
\end{gather}%
The next statement shows that (\ref{wkt}) defines an approximate solution to
integrated evolution equation (\ref{varcu}).

\begin{theorem}
\label{Theorem Dsmall}Let $\mathbf{\hat{h}}$ be a multi-wavepacket with
resonance invariant $nk$-spectrum $S$ with regularity degree $s$, $\mathbf{%
\vec{w}}$ be a solution of (\ref{syslocF}) and $\mathbf{\hat{w}}\left( 
\mathbf{k},\tau \right) $ be defined by (\ref{wkt}). Let \ 
\begin{equation}
\mathbf{\hat{D}}\left( \mathbf{\hat{w}}\right) =\mathbf{\hat{w}}-\mathcal{F}%
\left( \mathbf{\hat{w}}\right) -\mathbf{\hat{h}}.  \label{Dw}
\end{equation}%
Then there exists $\beta _{0}>0$ such that we have the estimate 
\begin{equation}
\left\Vert \mathbf{\hat{D}}\left( \mathbf{\hat{w}}\right) \right\Vert
_{E}\leq C\varrho +C\beta ^{s},\text{ if }0<\varrho \leq 1,\ \beta \leq
\beta _{0}.  \label{Dwb}
\end{equation}
\end{theorem}

\begin{proof}
Let 
\begin{equation}
\mathcal{F}^{-}\left( \mathbf{\hat{w}}\right) =\left(
1-\sum\nolimits_{l,\vartheta }\Psi _{i_{l},\vartheta }\Pi _{n_{l},\vartheta
}\right) \mathcal{F}\left( \mathbf{\hat{w}}\right) ,\ \mathbf{\hat{h}}^{-}=%
\mathbf{\hat{h}}-\sum\nolimits_{l,\vartheta }\Psi _{i_{l},\vartheta }\Pi
_{n_{l},\vartheta }\mathbf{\hat{h}}.  \label{Fhmin}
\end{equation}%
Summation of (\ref{sysloc}) with respect to $l,\vartheta $ yields 
\begin{equation*}
\mathbf{\hat{w}}=\sum\nolimits_{l,\vartheta }\Psi _{i_{l},\vartheta }\Pi
_{n_{l},\vartheta }\mathcal{F}\left( \mathbf{\hat{w}}\right)
+\sum\nolimits_{l,\vartheta }\Psi _{i_{l},\vartheta }\Pi _{n_{l},\vartheta }%
\mathbf{\hat{h}}.
\end{equation*}%
Hence, from (\ref{sysloc}) and (\ref{Dw}) we obtain 
\begin{equation}
\mathbf{\hat{D}}\left( \mathbf{\hat{w}}\right) =\mathbf{\hat{h}}^{-}-%
\mathcal{F}^{-}\left( \mathbf{\hat{w}}\right) .  \label{Dw1}
\end{equation}%
Using (\ref{hbold}) and (\ref{sourloc}) we consequently obtain 
\begin{equation*}
\left\Vert \Pi _{n_{l},\vartheta }\mathbf{\hat{h}}_{i}\right\Vert
_{L^{1}}\leq C\beta ^{s}\text{ if }n_{l}\neq n_{i};\ \left\Vert \Psi
_{i_{l},\vartheta }\mathbf{\hat{h}}_{i}\right\Vert _{L^{1}}\leq C\beta ^{s}%
\text{ if }\mathbf{k}_{\ast i_{l}}\neq \mathbf{k}_{\ast i},
\end{equation*}%
\begin{equation}
\left\Vert \mathbf{\hat{h}}^{-}\right\Vert _{E}\leq C_{1}\beta ^{s}.
\label{hbet}
\end{equation}%
Now, to show (\ref{Dwb}) it is sufficient to prove that 
\begin{equation}
\left\Vert \mathcal{F}^{-}\left( \mathbf{\hat{w}}\right) \right\Vert
_{E}\leq C_{2}\varrho .  \label{Fbet}
\end{equation}%
Obviously, 
\begin{equation}
\mathcal{F}^{-}\left( \mathbf{\hat{w}}\right) =\left(
1-\sum\nolimits_{l,\vartheta }\Psi _{i_{l},\vartheta }\Pi _{n_{l},\vartheta
}\right) \sum_{m}\mathcal{F}^{\left( m\right) }\left( \mathbf{\hat{w}}%
^{m}\right) .  \label{Fminw}
\end{equation}%
Note that 
\begin{equation}
\sum_{l,\vartheta }\Psi _{i_{l},\vartheta }\Pi _{n_{l},\vartheta
}=\sum\nolimits_{\vartheta =\pm }\sum\nolimits_{\left( n,k_{\ast }\right)
\in S}\Psi \left( \mathbf{\cdot },\vartheta \mathbf{k}_{\ast }\right) \Pi
_{n,\vartheta }.  \label{pspm}
\end{equation}%
Using (\ref{sumPi}) and ( \ref{Psiinf}) we consequently obtain 
\begin{gather}
\sum\nolimits_{\vartheta =\pm }\sum\nolimits_{\left( n,k_{\ast }\right) \in
\Sigma }\Psi \left( \mathbf{\cdot },\vartheta \mathbf{k}_{\ast }\right) \Pi
_{n,\vartheta }+\Psi _{\infty }=1,  \label{sum1} \\
\left( 1-\sum\nolimits_{l,\vartheta }\Psi _{i_{l},\vartheta }\Pi
_{n_{l},\vartheta }\right) =\Psi _{\infty }+\sum\nolimits_{\vartheta =\pm
}\sum\nolimits_{\left( n,k_{\ast }\right) \in \Sigma \setminus S}\Psi \left( 
\mathbf{\cdot },\vartheta \mathbf{k}_{\ast }\right) \Pi _{n,\vartheta }.
\label{psipi}
\end{gather}%
with $\Sigma $ defined in (\ref{ssp1}). Let us expand now $\mathcal{F}%
^{\left( m\right) }\left( \mathbf{\hat{w}}^{m}\right) $ using (\ref{binom}).
According to (\ref{Fminw}) and (\ref{psipi}) to prove (\ref{Fbet}) it is
sufficient to prove that \ for every string $\vec{\lambda}\in \Lambda ^{m}$ $%
\ $the following inequalities hold%
\begin{eqnarray}
\left\Vert \Psi _{\infty }\Pi _{n,\vartheta }\mathcal{F}^{\left( m\right)
}\left( \mathbf{\vec{w}}_{\vec{\lambda}}\right) \right\Vert &\leq
&C_{3}\varrho \text{ for }\left( n,\vartheta \right) \in \Lambda \text{, \
and }  \label{psiinfw} \\
\left\Vert \Psi \left( \mathbf{\cdot },\vartheta \mathbf{k}_{\ast }\right)
\Pi _{n,\vartheta }\mathcal{F}^{\left( m\right) }\left( \mathbf{\vec{w}}_{%
\vec{\lambda}}\right) \right\Vert &\leq &C_{3}\varrho ,\text{ \ if }\left( n,%
\mathbf{k}_{\ast }\right) \in \Sigma \setminus S\text{.}  \label{psiltw}
\end{eqnarray}%
We will use (\ref{piweq0}) and (\ref{weq0}) to obtain the above estimates. \
According to (\ref{Tplusmin})%
\begin{equation}
\mathcal{F}^{\left( m\right) }\left[ \mathbf{\vec{w}}_{\vec{\lambda}}\right]
\left( \mathbf{k},\tau \right) =\sum\nolimits_{n,\zeta }\sum\nolimits_{\vec{%
\xi}}\mathcal{F}_{n,\zeta ,\vec{\xi}}^{\left( m\right) }\left[ \mathbf{\hat{w%
}}_{\lambda _{1}}\ldots \mathbf{\hat{w}}_{\lambda _{m}}\right] \left( 
\mathbf{k},\tau \right) .  \label{Fbin}
\end{equation}%
Note that according to (\ref{piweq0}) if $\lambda _{i}=\left( l,\vartheta
^{\prime }\right) $ 
\begin{equation}
\mathbf{\hat{w}}_{\lambda _{i}}=\Pi _{n,\vartheta }\mathbf{\hat{w}}_{\lambda
_{i}},\ \text{if }n=n_{l}\text{ and }\vartheta ^{\prime }=\vartheta .
\label{wpiw}
\end{equation}%
Let us introduce notation 
\begin{equation}
\vec{n}\left( \vec{l}\right) =\left( n_{l_{1}},...,n_{l_{m}}\right) ,\ \vec{%
\xi}\left( \vec{\lambda}\right) =\left( \vec{n}\left( \vec{l}\right) ,\vec{%
\vartheta}\right) ,\text{\ for\ }\vec{\lambda}=\left( \vec{l},\vec{\vartheta}%
\right) \in \Lambda ^{m}.  \label{nl}
\end{equation}%
Since 
\begin{equation}
\Pi _{n^{\prime },\vartheta }\Pi _{n,\vartheta ^{\prime }}=0,\text{ if }%
n\neq n^{\prime }\text{ or }\vartheta ^{\prime }\neq \vartheta  \label{PnP'}
\end{equation}%
then (\ref{wpiw}) implies 
\begin{eqnarray}
\mathcal{F}_{n,\zeta ,\vec{\xi}}^{\left( m\right) }\left[ \mathbf{\hat{w}}%
_{\lambda _{1}}\ldots \mathbf{\hat{w}}_{\lambda _{m}}\right] &=&0\text{\ if }%
\vec{\xi}=\left( \vec{n},\vec{\zeta}\right) \neq \vec{\xi}\left( \vec{\lambda%
}\right) ,\text{ and, hence,}  \notag \\
\mathcal{F}^{\left( m\right) }\left[ \mathbf{\vec{w}}_{\vec{\lambda}}\right]
\left( \mathbf{k},\tau \right) &=&\sum\nolimits_{n,\zeta }\mathcal{F}%
_{n,\zeta ,\vec{\xi}\left( \vec{\lambda}\right) }^{\left( m\right) }\left[ 
\mathbf{\hat{w}}_{\lambda _{1}}\ldots \mathbf{\hat{w}}_{\lambda _{m}}\right]
\left( \mathbf{k},\tau \right) ,  \label{nzi}
\end{eqnarray}%
where we use notation (\ref{lamprop}), (\ref{nl}). Note also that 
\begin{equation}
\Pi _{n^{\prime },\vartheta }\mathcal{F}_{n,\zeta ,\vec{\xi}}^{\left(
m\right) }=0\text{\ if }n^{\prime }\neq n\text{\ or\ }\vartheta \neq \zeta ,
\label{Pinth0}
\end{equation}%
and, hence, we have nonzero $\Pi _{n^{\prime },\vartheta }\mathcal{F}%
_{n,\zeta ,\vec{\xi}}^{\left( m\right) }\left( \mathbf{\vec{w}}_{\vec{\lambda%
}}\right) $ only if 
\begin{equation}
\vec{\xi}=\vec{\xi}\left( \vec{\lambda}\right) ,\ n^{\prime }=n\text{, }%
\vartheta =\zeta .  \label{nonzero}
\end{equation}%
By (\ref{Fm})%
\begin{gather}
\mathcal{F}_{n,\zeta ,\vec{\xi}\left( \vec{\lambda}\right) }^{\left(
m\right) }\left( \mathbf{\vec{w}}_{\vec{\lambda}}\right) \left( \mathbf{k}%
,\tau \right) =\int_{0}^{\tau }\int_{\mathbb{D}_{m}}\exp \left\{ \mathrm{i}%
\phi _{n,\zeta ,\vec{\xi}\left( \vec{\lambda}\right) }\left( \mathbf{\mathbf{%
k}},\vec{k}\right) \frac{\tau _{1}}{\varrho }\right\}  \label{Fmw} \\
\chi _{n,\zeta ,\vec{\xi}\left( \vec{\lambda}\right) }^{\left( m\right)
}\left( \mathbf{\mathbf{k}},\vec{k}\right) \left[ \mathbf{\hat{w}}_{\lambda
_{1}}\left( \mathbf{k}^{\prime },\tau _{1}\right) ,\ldots ,\mathbf{\hat{w}}%
_{\lambda _{m}}\left( \mathbf{k}^{\left( m\right) }\left( \mathbf{\mathbf{k}}%
,\vec{k}\right) ,\tau _{1}\right) \right] \mathrm{\tilde{d}}^{\left(
m-1\right) d}\vec{k}\mathrm{d}\tau _{1},  \notag
\end{gather}%
Now we use (\ref{weq0}) and notice that according to the convolution
identity in (\ref{Fmintr}) 
\begin{equation}
\left\vert \mathbf{\hat{w}}_{\lambda _{1}}\left( \mathbf{k}^{\prime },\tau
_{1}\right) \right\vert \cdot ...\cdot \left\vert \mathbf{\hat{w}}_{\lambda
_{m}}\left( \mathbf{k}^{\left( m\right) }\left( \mathbf{\mathbf{k}},\vec{k}%
\right) ,\tau _{1}\right) \right\vert =0\text{ if }\left\vert \mathbf{k}%
-\sum\nolimits_{i}\vartheta _{i}\mathbf{k}_{\ast l_{i}}\right\vert \geq
m\beta ^{1-\epsilon }.  \label{prweq0}
\end{equation}%
Hence the integral (\ref{Fmw}) is nonzero only if $\left( \mathbf{k},\vec{k}%
\right) $ belongs to the set 
\begin{equation}
B_{\beta }=\left\{ \left( \mathbf{k},\vec{k}\right) :\left\vert \mathbf{k}%
^{\left( i\right) }-\vartheta _{i}\mathbf{k}_{\ast l_{i}}\right\vert \leq
\beta ^{1-\epsilon },\ i=1,...,m,\ \left\vert \mathbf{k}-\sum\nolimits_{i}%
\vartheta _{i}\mathbf{k}_{\ast l_{i}}\right\vert \leq m\beta ^{1-\epsilon
}\right\} .  \label{konly}
\end{equation}%
We will prove now that if $\left( n,\mathbf{k}_{\ast i}\right) \notin S$
then for small $\beta $ one of the following alternatives holds: \ 
\begin{gather}
\text{either }\Psi \left( \mathbf{\cdot },\vartheta \mathbf{k}_{\ast
i}\right) \Pi _{n^{\prime },\vartheta }\mathcal{F}_{n,\zeta ,\vec{\xi}%
}^{\left( m\right) }\left( \mathbf{\vec{w}}_{\vec{\lambda}}\right) =0
\label{alt1} \\
\text{or (\ref{nonzero}) holds and}\left\vert \phi _{n,\zeta ,\vec{\xi}%
}\left( \mathbf{\mathbf{k}},\vec{k}\right) \right\vert \geq c>0\text{ for }%
\left( \mathbf{k},\vec{k}\right) \in B_{\beta }.  \label{figrc}
\end{gather}%
\ Note then since $\phi _{n,\zeta ,\vec{\xi}}\left( \mathbf{\mathbf{k}},\vec{%
k}\right) $ is smooth then using notation (\ref{kapzel}) we get 
\begin{gather}
\left\vert \phi _{n,\zeta ,\vec{\xi}}\left( \mathbf{\mathbf{k}},\vec{k}%
\right) -\phi _{n^{\prime },\zeta ,\vec{\xi}}\left( \mathbf{\mathbf{k}}%
_{\ast \ast },\vec{k}_{\ast }\right) \right\vert \leq C\beta ^{1-\epsilon }%
\text{ for }\left( \mathbf{k},\vec{k}\right) \in B_{\beta },  \label{fiminfi}
\\
\vec{\vartheta}=\left( \vartheta _{1},...,\vartheta _{m}\right) ,\ \mathbf{%
\mathbf{k}}_{\ast \ast }=\zeta \sum\nolimits_{i}\vartheta _{i}\mathbf{k}%
_{\ast l_{i}}=\zeta \varkappa _{m}\left( \vec{\vartheta},\vec{l}\right) , 
\notag
\end{gather}%
Hence the alternative (\ref{figrc}) holds if%
\begin{equation}
\phi _{n,\zeta ,\vec{\xi}}\left( \mathbf{\mathbf{k}}_{\ast \ast },\vec{k}%
_{\ast }\right) \neq 0,  \label{fineq0}
\end{equation}%
\ and, consequently, it suffices to prove that either (\ref{alt1}) or (\ref%
{fineq0}) holds. Combining (\ref{konly}) with $\Psi \left( \mathbf{k}%
,\vartheta \mathbf{k}_{\ast i}\right) =0$ for $\left\vert \mathbf{k}%
-\vartheta \mathbf{k}_{\ast i}\right\vert \geq \beta ^{1-\epsilon }$ we find
that $\Psi _{i,\vartheta }\mathcal{F}^{\left( m\right) }\left[ \mathbf{\vec{w%
}}_{\vec{\lambda}}\right] $ can be non-zero for small $\beta $ only in a
small neighborhood of a point $\zeta \varkappa _{m}\left( \vec{\vartheta},%
\vec{l}\right) \in \left[ S\right] _{K,\text{out}}$, and that is possible
only if 
\begin{equation}
\mathbf{\mathbf{k}}_{\ast \ast }=\zeta \varkappa _{m}\left( \vec{\vartheta},%
\vec{l}\right) =\vartheta \mathbf{k}_{\ast i},\ \mathbf{k}_{\ast i}\in K_{S}.
\label{kstar}
\end{equation}%
Let us show that the equality 
\begin{equation}
\phi _{n,\zeta ,\vec{\xi}}\left( \mathbf{\mathbf{k}}_{\ast \ast },\vec{k}%
_{\ast }\right) =0  \label{fiopp}
\end{equation}%
is impossible for $\mathbf{\mathbf{k}}_{\ast \ast }$ as in (\ref{kstar}) and 
$n^{\prime }=n$ as in (\ref{Pinth0}), keeping in mind that $\left( n,\mathbf{%
\mathbf{k}}_{\ast i}\right) \notin S$. It follows from (\ref{Omzet}) and (%
\ref{phim}) that the equation (\ref{fiopp}) has the form of the resonance
equation (\ref{Omeq0}). Since $nk$-spectrum $S$ is resonance invariant, in
view of Definition \ref{Definition omclos} the resonance equation (\ref%
{fiopp}) may have a solution only if $\mathbf{\mathbf{k}}_{\ast \ast }=%
\mathbf{\mathbf{k}}_{\ast i}$, $i=i_{l}$, $n=n_{l}$, with $\left( n_{l},%
\mathbf{\mathbf{k}}_{\ast i_{l}}\right) \in S$. Since $\left( n,\mathbf{%
\mathbf{k}}_{\ast i}\right) \notin S$ that implies (\ref{fiopp}) does not
have a solution and, hence, (\ref{fineq0}) holds when $\left( n,\mathbf{%
\mathbf{k}}_{\ast i}\right) \notin S$. Notice that Theorem \ref{Theorem
existpr} and (\ref{dtf1}) yield bounds 
\begin{equation*}
\left\Vert \mathbf{\hat{w}}_{\lambda _{i}}\right\Vert _{E}\leq R_{1},\
\left\Vert \partial _{\tau }\mathbf{\hat{w}}_{\lambda _{i}}\right\Vert
_{E}\leq C.
\end{equation*}%
These bounds combined with Lemma \ref{Lemma intbyparts}, proven below, imply
that if (\ref{fineq0}) holds then (\ref{psiltw}) holds. Now let us turn to (%
\ref{psiinfw}). According to ( \ref{Psiinf}) and (\ref{prweq0}) the term $%
\Psi _{\infty }\Pi _{n^{\prime },\vartheta }\mathcal{F}^{\left( m\right)
}\left( \mathbf{\vec{w}}_{\vec{\lambda}}\right) $ can be non-zero only if $%
\zeta \varkappa _{m}\left( \vec{\lambda}\right) =\mathbf{\mathbf{k}}_{\ast
\ast }\notin K_{S}$. Since $nk$-spectrum $S$ is resonance invariant we
conclude as above that inequality (\ref{fineq0}) holds in this case as well.
The fact that the set of all $\varkappa _{m}\left( \vec{\lambda}\right) $ is
finite, combined with inequality (\ref{fineq0}), imply (\ref{figrc}) for
sufficiently small $\beta $. Using Lemma \ref{Lemma intbyparts} as above we
derive (\ref{psiinfw}). Hence, all terms in the expansion (\ref{Fminw}) are
either zero or satisfy (\ref{psiinfw}) or (\ref{psiltw}) implying
consequently (\ref{Fbet}) and (\ref{Dwb}).
\end{proof}

Here is the lemma used in the above proof.

\begin{lemma}
\label{Lemma intbyparts}Assume that%
\begin{gather}
\left\vert \Psi _{i,\vartheta ^{\prime }}\Pi _{n^{\prime },\zeta }\chi
_{n,\zeta ,\vec{\xi}}^{\left( m\right) }\left( \mathbf{\mathbf{k}},\vec{k}%
\right) \left[ \mathbf{\hat{w}}_{\lambda _{1}}\left( \mathbf{k}^{\prime
},\tau _{1}\right) ,\ldots ,\mathbf{\hat{w}}_{\lambda _{m}}\left( \mathbf{k}%
^{\left( m\right) }\left( \mathbf{\mathbf{k}},\vec{k}\right) ,\tau
_{1}\right) \right] \right\vert =0\text{ for }\left( \mathbf{\mathbf{k}},%
\vec{k}\right) \in B_{\beta },  \notag \\
\text{and }\left\vert \phi _{n,\zeta ,\vec{\xi}}\left( \mathbf{\mathbf{k}},%
\vec{k}\right) \right\vert \geq \omega _{\ast }>0\text{ for }\left( \mathbf{%
\mathbf{k}},\vec{k}\right) \notin B_{\beta },\text{ with }B_{\beta }\text{
as in (\ref{konly}).}  \label{figrom}
\end{gather}%
Then%
\begin{gather}
\left\Vert \Psi \left( \mathbf{\cdot },\vartheta ^{\prime }\mathbf{k}_{\ast
i}\right) \Pi _{n^{\prime },\zeta }\mathcal{F}_{n,\zeta ,\vec{\xi}}^{\left(
m\right) }\left( \mathbf{\vec{w}}_{\vec{\lambda}}\right) \right\Vert _{E}\leq
\label{estrho} \\
\frac{4\varrho }{\omega _{\ast }}\left\Vert \chi ^{\left( m\right)
}\right\Vert \dprod\nolimits_{j}\left\Vert \mathbf{\hat{w}}_{\lambda
_{j}}\right\Vert _{E}+\frac{2\varrho \tau _{\ast }}{\omega _{\ast }}%
\left\Vert \chi ^{\left( m\right) }\right\Vert \sum\nolimits_{i}\left\Vert
\partial _{\tau }\mathbf{\hat{w}}_{\lambda _{i}}\right\Vert
_{E}\dprod\nolimits_{j\neq i}\left\Vert \mathbf{\hat{w}}_{\lambda
_{j}}\right\Vert _{E}.  \notag
\end{gather}
\end{lemma}

\begin{proof}
Notice that the oscillatory factor in (\ref{Fm}) equals to 
\begin{equation*}
\exp \left\{ \mathrm{i}\phi \left( \mathbf{\mathbf{k}},\vec{k}\right) \frac{%
\tau _{1}}{\varrho }\right\} =\frac{\varrho }{\mathrm{i}\phi \left( \mathbf{%
\mathbf{k}},\vec{k}\right) }\partial _{\tau _{1}}\exp \left\{ \mathrm{i}\phi
\left( \mathbf{\mathbf{k}},\vec{k}\right) \frac{\tau _{1}}{\varrho }\right\}
.
\end{equation*}%
Denoting $\phi _{n,\zeta ,\vec{\xi}}=\phi $, $\Psi _{i,\vartheta ^{\prime
}}\Pi _{n^{\prime },\zeta }\chi _{n,\zeta ,\vec{\xi}}^{\left( m\right)
}=\chi _{\vec{\eta}}^{\left( m\right) }$ and integrating (\ref{Fm})\ by
parts with respect to $\tau _{1}$ we obtain%
\begin{gather}
\Psi \left( \mathbf{k},\vartheta ^{\prime }\mathbf{k}_{\ast i}\right) \Pi
_{n^{\prime },\zeta }\mathcal{F}_{n,\zeta ,\vec{\xi}}^{\left( m\right)
}\left( \mathbf{\vec{w}}_{\vec{\lambda}}\right) \left( \mathbf{k},\tau
\right) =  \label{NFMintbp} \\
\int_{B}\Psi \left( \mathbf{k},\vartheta ^{\prime }\mathbf{k}_{\ast
i}\right) \frac{\varrho \mathrm{e}^{\mathrm{i}\phi \left( \mathbf{\mathbf{k}}%
,\vec{k}\right) \frac{\tau }{\varrho }}}{\mathrm{i}\phi \left( \mathbf{%
\mathbf{k}},\vec{k}\right) }\chi _{\vec{\eta}}^{\left( m\right) }\left( 
\mathbf{\mathbf{k}},\vec{k}\right) \mathbf{\hat{w}}_{\lambda _{1}}\left( 
\mathbf{k}^{\prime },\tau \right) \ldots \mathbf{\hat{w}}_{\lambda
_{m}}\left( \mathbf{k}^{\left( m\right) }\left( \mathbf{k},\vec{k}\right)
,\tau \right) \,\mathrm{\tilde{d}}^{\left( m-1\right) d}\vec{k}  \notag \\
-\int_{B}\Psi \left( \mathbf{k},\vartheta ^{\prime }\mathbf{k}_{\ast
i}\right) \frac{\varrho }{\mathrm{i}\phi \left( \mathbf{\mathbf{k}},\vec{k}%
\right) }\chi _{\vec{\eta}}^{\left( m\right) }\left( \mathbf{\mathbf{k}},%
\vec{k}\right) \mathbf{\hat{w}}_{\lambda _{1}}\left( \mathbf{k}^{\prime
},0\right) \ldots \mathbf{\hat{w}}_{\lambda _{m}}\left( \mathbf{k}^{\left(
m\right) }\left( \mathbf{k},\vec{k}\right) ,0\right) \,\mathrm{\tilde{d}}%
^{\left( m-1\right) d}\vec{k}  \notag \\
-\int_{0}^{\tau }\int_{B}\Psi \left( \mathbf{k},\vartheta ^{\prime }\mathbf{k%
}_{\ast i}\right) \frac{\varrho \mathrm{e}^{\mathrm{i}\phi \left( \mathbf{%
\mathbf{k}},\vec{k}\right) \frac{\tau _{1}}{\varrho }}}{\mathrm{i}\phi
\left( \mathbf{\mathbf{k}},\vec{k}\right) }\chi _{\vec{\eta}}^{\left(
m\right) }\left( \mathbf{\mathbf{k}},\vec{k}\right) \partial _{\tau _{1}}%
\left[ \mathbf{\hat{w}}_{\lambda _{1}}\left( \mathbf{k}^{\prime }\right)
\ldots \mathbf{\hat{w}}_{\lambda _{m}}\left( \mathbf{k}^{\left( m\right)
}\left( \mathbf{k},\vec{k}\right) \right) \right] \,\mathrm{\tilde{d}}%
^{\left( m-1\right) d}\vec{k}d\tau _{1},  \notag
\end{gather}%
where $B$ is the set of $\mathbf{k}^{\left( i\right) }$ for which (\ref%
{konly}) holds. The relations (\ref{chiCR}) and (\ref{j0}) imply $\left\vert
\chi _{\vec{\eta}}^{\left( m\right) }\left( \mathbf{\mathbf{k}},\vec{k}%
\right) \right\vert \leq \left\Vert \chi ^{\left( m\right) }\right\Vert $.
Using then (\ref{figrom}),the Leibnitz formula and (\ref{Yconv}) we obtain (%
\ref{estrho}).
\end{proof}

The main result of this subsection is the next theorem which, when combined
with Lemma \ref{Lemma wiswave}, implies the wavepacket preservation, namely
that the solution $\mathbf{\hat{u}}_{n,\vartheta }\left( \mathbf{k},\tau
\right) $ of (\ref{equfa}) is a multi-wavepacket for all $\tau \in \left[
0,\tau _{\ast }\right] $.

\begin{theorem}
\label{Theorem uminw} Assume that conditions of Theorem \ref{Theorem Dsmall}
are fulfilled. Let $\mathbf{\hat{u}}_{n,\vartheta }\left( \mathbf{k},\tau
\right) $ for $n=n_{l}$ and $\mathbf{\hat{w}}_{l,\vartheta }\left( \mathbf{k}%
,\tau \right) \ $be the solutions to respective systems (\ref{equfa}) (\ref%
{sysloc}), $\mathbf{\hat{w}}$ be defined by (\ref{wkt}). Then there exists $%
\beta _{0}>0$ such that 
\begin{equation}
\left\Vert \mathbf{\hat{u}}_{n_{l},\vartheta }-\Pi _{n_{l},\vartheta }%
\mathbf{\hat{w}}\right\Vert _{E}\leq C\varrho +C^{\prime }\beta ^{s}\text{ \
for}\ 0<\beta \leq \beta _{0}.  \label{uminw}
\end{equation}
\end{theorem}

\begin{proof}
Note that $\mathbf{\hat{u}}_{n,\vartheta }=\Pi _{n,\vartheta }\mathbf{\hat{u}%
}$ where $\mathbf{\hat{u}}$ is a solution of (\ref{varcu}) and, according to
Theorem \ref{Theorem exist}, $\left\Vert \mathbf{\hat{u}}\right\Vert
_{E}\leq 2R$. Comparing the equations (\ref{varcu}) and (\ref{Dw}) , which
are $\mathbf{\hat{u}}=\mathcal{F}\left( \mathbf{\hat{u}}\right) +\mathbf{%
\hat{h}}$ and $\mathbf{\hat{w}}=\mathcal{F}\left( \mathbf{\hat{w}}\right) +%
\mathbf{\hat{h}+\hat{D}}\left( \mathbf{\hat{w}}\right) $, we find that Lemma %
\ref{Lemma contr} can be applied. Then we notice that by Lemma \ref{Lemma
Flip} $\mathcal{F}$ has the Lipschitz constant $C_{F}\tau _{\ast }$ for such 
$\mathbf{\hat{u}}$. Taking $C_{F}\tau _{\ast }<1$ as in Theorem \ref{Theorem
exist} we obtain (\ref{uminw}) from (\ref{iminu0}).
\end{proof}

Notice that Theorem \ref{Theorem sumwave} \ is a direct corollary of Theorem %
\ref{Theorem uminw} and Lemma \ref{Lemma wiswave}. The following corollary
shows that inequality (\ref{uminw}) and, therefore, Theorems \ref{Theorem
sumwave} and \ref{Theorem invarwave} on preservation of wavepackets hold in
the case when the coefficients of operator $\mathbf{\hat{F}}\left( \mathbf{%
\hat{U}}\right) $ in (\ref{difeqfou}), (\ref{Fmintr}) regularly depends on
small $\varrho $, $\mathbf{\hat{F}}\left( \mathbf{\hat{U}}\right) =\mathbf{%
\hat{F}}\left( \mathbf{\hat{U}},\varrho \right) $.

\begin{corollary}[parameter dependent nonlinearity]
\label{Corollary Fofrho}Assume that conditions of \ Theorem \ref{Theorem
Dsmall} are fulfilled. Consider a perturbed equation (\ref{varcu}) $\mathbf{%
\hat{u}}\left( \mathbf{k},\tau \right) =\mathcal{F}\left( \mathbf{\hat{u}}%
\right) \left( \mathbf{k},\tau \right) +\mathcal{F}_{1}\left( \mathbf{\hat{u}%
},\varrho \right) \left( \mathbf{k},\tau \right) +\mathbf{\hat{h}}\left( 
\mathbf{k}\right) $ where operator $\mathcal{F}_{1}\left( \mathbf{\hat{u}}%
,\varrho \right) $ satisfies the inequality $\left\Vert \mathcal{F}%
_{1}\left( \mathbf{\hat{u}},\varrho \right) \right\Vert _{E}\leq C\varrho
^{q}$ for $\left\Vert \mathbf{\hat{u}}\right\Vert _{E}\leq 2R$ with some $q,$
$0<q\leq $ $1$. Let $\mathbf{\hat{w}}_{l,\vartheta }\left( \mathbf{k},\tau
\right) \ $be the solution of \ (\ref{sysloc}). Then \ $\left\Vert \Pi
_{n,\vartheta }\mathbf{\hat{u}\ }-\mathbf{\hat{w}}_{l,\vartheta }\right\Vert
_{E}\leq C\varrho ^{q}+C^{\prime }\beta ^{s}$.
\end{corollary}

\begin{proof}
The statement follows from (\ref{uminw}) and Lemma \ref{Lemma contr}.
\end{proof}

The following theorem shows that any multi-wavepacket solution to (\ref%
{varcu}) yields a solution to the wavepacket interaction system (\ref{sysloc}%
).

\begin{theorem}
\label{Theorem necessary}Let $\mathbf{\hat{u}}\left( \mathbf{k},\tau \right) 
$ be a solution of (\ref{varcu}) and assume that $\mathbf{\hat{u}}\left( 
\mathbf{k},\tau \right) $ and $\mathbf{\hat{h}}\left( \mathbf{k}\right) $
are multiwavepackets with $nk$-spectrum $S=\left\{ \left( n_{l},\mathbf{k}%
_{\ast l}\right) \text{, }l=1,...,N\right\} $ and the regularity degree $s$.
Let also $\Psi _{i_{l},\vartheta }=\Psi _{i_{l},\vartheta }$ be defined by (%
\ref{Psilz}). Then functions $\mathbf{\hat{w}}_{l,\vartheta }^{\prime
}\left( \mathbf{k},\tau \right) =\Psi _{i_{l},\vartheta }\Pi
_{n_{l},\vartheta }\mathbf{\hat{u}}\left( \mathbf{k},\tau \right) $ are a
solution to the system (\ref{sysloc}) with $\mathbf{\hat{h}}\left( \mathbf{k}%
\right) $ replaced by $\mathbf{\hat{h}}^{\prime }\left( \mathbf{k},\tau
\right) $ satisfying 
\begin{equation}
\left\Vert \mathbf{\hat{h}}\left( \mathbf{k}\right) -\mathbf{\hat{h}}%
^{\prime }\left( \mathbf{k},\tau \right) \right\Vert _{L^{1}}\leq C\beta
^{s},\;0\leq \tau \leq \tau _{\ast }.  \label{hprime}
\end{equation}
\end{theorem}

\begin{proof}
Multiplying (\ref{varcu}) by $\Psi _{i_{l},\vartheta }\Pi _{n_{l},\vartheta
} $ we get%
\begin{equation}
\mathbf{\hat{w}}_{l,\vartheta }^{\prime }=\Psi \left( \mathbf{\cdot }%
,\vartheta \mathbf{k}_{\ast i_{l}}\right) \Pi _{n_{l},\vartheta }\mathcal{F}%
\left( \mathbf{\hat{u}}\right) \left( \mathbf{k},\tau \right) +\Psi \left( 
\mathbf{\cdot },\vartheta \mathbf{k}_{\ast i_{l}}\right) \Pi
_{n_{l},\vartheta }\mathbf{\hat{h}}\left( \mathbf{k}\right) ,\ \mathbf{\hat{w%
}}_{l,\vartheta }^{\prime }=\Psi \left( \mathbf{\cdot },\vartheta \mathbf{k}%
_{\ast i_{l}}\right) \Pi _{n_{l},\vartheta }\mathbf{\hat{u}}.  \label{wprime}
\end{equation}%
Since $\mathbf{\hat{u}}\left( \mathbf{k},\tau \right) $ is a multiwavepacket
with regularity $s$ we have 
\begin{equation}
\left\Vert \mathbf{\hat{u}}\left( \mathbf{\cdot },\tau \right) -\mathbf{\hat{%
w}}^{\prime }\left( \mathbf{\cdot },\tau \right) \right\Vert _{L^{1}}\leq
C_{\epsilon }\beta ^{s}\text{ \ where \ }\mathbf{\hat{w}}^{\prime }\left( 
\mathbf{\cdot },\tau \right) =\sum\nolimits_{l,\vartheta }\Psi \left( 
\mathbf{\cdot },\vartheta \mathbf{k}_{\ast i_{l}}\right) \mathbf{\hat{u}}%
\left( \mathbf{\cdot },\tau \right) .  \label{uwprime}
\end{equation}%
Let us recast (\ref{wprime}) in the form 
\begin{gather}
\mathbf{\hat{w}}_{l,\vartheta }^{\prime }=\Psi \left( \mathbf{\cdot }%
,\vartheta \mathbf{k}_{\ast i_{l}}\right) \Pi _{n_{l},\vartheta }\mathcal{F}%
\left( \mathbf{\hat{w}}^{\prime }\right) \left( \mathbf{k},\tau \right)
+\Psi \left( \mathbf{\cdot },\vartheta \mathbf{k}_{\ast i_{l}}\right) \Pi
_{n_{l},\vartheta }\left[ \mathbf{\hat{h}}\left( \mathbf{k}\right) +\mathbf{%
\hat{h}}^{\prime \prime }\left( \mathbf{k},\tau \right) \right] ,
\label{wprime1} \\
\mathbf{\hat{h}}^{\prime \prime }\left( \mathbf{k},\tau \right) =\left[ 
\mathcal{F}\left( \mathbf{\hat{u}}\right) -\mathcal{F}\left( \mathbf{\hat{w}}%
^{\prime }\right) \right] \left( \mathbf{k},\tau \right) .  \notag
\end{gather}%
Denoting $\mathbf{\hat{h}}\left( \mathbf{k}\right) +\mathbf{\hat{h}}^{\prime
\prime }\left( \mathbf{k},\tau \right) =$ $\mathbf{\hat{h}}^{\prime }\left( 
\mathbf{k},\tau \right) $ we observe that (\ref{wprime1}) has the form of (%
\ref{sysloc}) with $\mathbf{\hat{h}}\left( \mathbf{k}\right) $ replaced by $%
\mathbf{\hat{h}}^{\prime }\left( \mathbf{k},\tau \right) $. Inequality (\ref%
{hprime}) follows then from (\ref{uwprime}) and (\ref{Flip}).
\end{proof}

\section{Reduction of wavepacket interaction system to a minimal interaction
system}

\ Our goal in this section is to substitute the wavepacket interaction
system (\ref{sysloc}) with a simpler (minimal) \emph{interaction system}
which describes the evolution of wavepackets with the same accuracy. We fix
the $nk$-spectrum $S=\left\{ \left( n_{l},\mathbf{k}_{\ast l}\right) ,\
l=1,...,N\right\} $ of the initial multiwavepacket and assume everywhere
below that it is resonance invariant. The minimal interaction system is
built based on operators $\mathbf{L}$ and $\mathbf{\hat{F}}\left( \mathbf{%
\hat{U}}\right) $ and on $S$. We want the minimal interaction system to
satisfy the following requirements. Firstly, the approximation of solutions
of (\ref{sysloc}) \ by solutions of the minimal interaction system of the
order $\left( \mu ,\nu \right) $ \ has to be of the order $\varrho $ in
suitable region of parameters $\left( \varrho ,\beta \right) $ (which is
larger for larger $\mu ,\nu $). Secondly, the minimal interaction system of
the order $\left( \mu ,\nu \right) $ should be defined by $S$ and by the
values of $\mathbf{L}\left( \mathbf{k}\right) $ and its derivatives of the
order up to $\mu $ and by the values $\chi ^{\left( m\right) }\left( \mathbf{%
\mathbf{k}},\vec{k}\right) $ and its derivatives of order up to $\nu $ at $\ 
\mathbf{k}_{\ast l}\in S_{K}$.

The construction of the minimal interaction system consists of the following
consecutive steps: (i) introduction of a time averaged wavepacket
interaction system obtained by discarding non-resonant terms in the
nonlinearity; (ii) reduction of the system for vector components $\mathbf{%
\hat{v}}_{l,\vartheta }$ to an equivalent one for scalar amplitudes $\hat{v}%
_{l,\vartheta }$; (iii) change of variables $\mathbf{k}=\vartheta \mathbf{k}%
_{\ast l}+\beta \mathbf{\eta }$ in the equation for $\hat{v}_{l,\vartheta }$
resulting in a regular dependence of coefficients on small $\beta \mathbf{%
\eta }$; (iv) substitution of the general dependence on $\beta \mathbf{\eta }
$ in the linear part with a certain polynomial one of the order $\mu $, and
the general dependence on $\beta \mathbf{\eta }$ of coefficients of the
nonlinearity with a certain trigonometric polynomial of the order $\nu $;
(v) substitution of the cutoff functions $\Psi \left( \mathbf{\cdot }%
,\vartheta \mathbf{k}_{\ast i_{l}}\right) $ from (\ref{sysloc}), which were
preserved up to this step, with $1$.

As a result we obtain a minimal interaction system with weakly universal
nonlinearity, which in the simplest case, where $S$ is just a single element 
$\left( \mathbf{k}_{\ast },n\right) $, is equivalent to the classical NLS
equation, and in the case when $S$ consists of only two elements $\left( 
\mathbf{k}_{\ast },n\right) $, $\left( -\mathbf{k}_{\ast },n\right) $, is
equivalent to the classical coupled modes system.

\subsection{Time averaged wavepacket interaction system}

Here we modify the wavepacket interaction system (\ref{sysloc}),
substituting its nonlinearity with a certain universal or conditionally
universal one obtained by the time averaging, and prove that this
substitution produces a small error of order $\varrho $. As the first step
we recast (\ref{sysloc}) in a slightly different form by using expansions (%
\ref{binom}), (\ref{Fbin}) together with (\ref{nzi})\ and (\ref{Pinth0}) and
writing the nonlinearity in the equation (\ref{sysloc}) in the form 
\begin{eqnarray*}
\Psi \left( \mathbf{\cdot },\vartheta \mathbf{k}_{\ast i_{l}}\right) \Pi
_{n_{l},\vartheta }\mathcal{F}\left( \mathbf{\cdot },\tau \right)
&=&\sum\nolimits_{m\in \mathfrak{M}_{F}}\sum\nolimits_{\vec{\lambda}\in
\Lambda ^{m}}\Psi \left( \mathbf{\cdot },\vartheta \mathbf{k}_{\ast
i_{l}}\right) \mathcal{F}_{n_{l},\vartheta ,\vec{\xi}\left( \vec{\lambda}%
\right) }^{\left( m\right) }\left( \mathbf{\vec{w}}_{\vec{\lambda}}\right) ,%
\text{\ }\vec{\lambda}=\left( \vec{l},\vec{\zeta}\right) , \\
\mathcal{F}_{n_{l},\vartheta ,\vec{\xi}\left( \vec{\lambda}\right) }^{\left(
m\right) }\left( \mathbf{\vec{w}}_{\vec{\lambda}}\right) \left( \mathbf{k}%
,\tau \right) &=&\mathcal{F}_{n,\zeta ,\vec{n},\vec{\zeta}}^{\left( m\right)
}\left[ \mathbf{\hat{w}}_{\lambda _{1}}\ldots \mathbf{\hat{w}}_{\lambda _{m}}%
\right] \left( \mathbf{k},\tau \right) ,\vec{n}=\vec{n}\left( \vec{l}\right)
,\text{\ }\left( n,\zeta \right) =\left( n_{l},\vartheta \right) ,
\end{eqnarray*}%
with $\mathcal{F}_{n,\zeta ,\vec{n},\vec{\zeta}}^{\left( m\right) }$ as in (%
\ref{Fm}) and $\vec{n}\left( \vec{l}\right) \ $as in (\ref{nl}).
Consequently, the wavepacket interaction system (\ref{sysloc}) can be
written in an equivalent form 
\begin{equation}
\mathbf{\hat{w}}_{l,\vartheta }=\sum_{m\in \mathfrak{M}_{F}}\sum_{\vec{%
\lambda}\in \Lambda ^{m}}\Psi \left( \mathbf{\cdot },\vartheta \mathbf{k}%
_{\ast i_{l}}\right) \mathcal{F}_{n_{l},\vartheta ,\vec{\xi}\left( \vec{%
\lambda}\right) }^{\left( m\right) }\left( \mathbf{\vec{w}}_{\vec{\lambda}%
}\right) +\Psi \left( \mathbf{\cdot },\vartheta \mathbf{k}_{\ast
i_{l}}\right) \Pi _{n_{l},\vartheta }\mathbf{\hat{h}},\ l=1,...N,\ \vartheta
=\pm .  \label{sysloc1}
\end{equation}%
The construction of the above mentioned time averaged equation reduces to
discarding certain terms in the original system (\ref{sysloc1}). First we
introduce the following sets of indices related to the resonance equation (%
\ref{Omeq0}) and $\Omega _{m}$ defined by (\ref{Omzet}): 
\begin{equation}
\Lambda _{n_{l},\vartheta }^{m}=\left\{ \vec{\lambda}=\left( \vec{l},\vec{%
\zeta}\right) \in \Lambda ^{m}:\Omega _{m}\left( \vartheta ,n_{l},\vec{%
\lambda}\right) =0\right\} ,  \label{resset}
\end{equation}%
and then the \emph{time-averaged nonlinearity by} 
\begin{equation}
\mathcal{F}_{\limfunc{av},n_{l},\vartheta }\left( \mathbf{\vec{w}}\right)
=\sum\nolimits_{m\in \mathfrak{M}_{F}}\mathcal{F}_{n_{l},\vartheta }^{\left(
m\right) },\;\mathcal{F}_{n_{l},\vartheta }^{\left( m\right)
}=\sum\nolimits_{\vec{\lambda}\in \Lambda _{n_{l},\vartheta }^{m}}\mathcal{F}%
_{n_{l},\vartheta ,\vec{\xi}\left( \vec{\lambda}\right) }^{\left( m\right)
}\left( \mathbf{\vec{w}}_{\vec{\lambda}}\right) .  \label{Fav}
\end{equation}%
Note that the nonlinearity $\mathcal{F}_{\limfunc{av},n_{l},\vartheta
}^{\left( m\right) }\left( \mathbf{\vec{w}}\right) $ can be obtained from $%
\mathcal{F}_{n_{l},\vartheta }^{\left( m\right) }$ by the averaging \
formula (\ref{Gav}) where $A_{T}\mathbf{\ }$is defined by formula (\ref%
{cana1}) with frequencies $\phi _{j}=\omega _{n_{j}}\left( \mathbf{k}_{\ast
i_{j}}\right) $. Consequently, the desired equation with time-averaged
nonlinearity is 
\begin{equation}
\mathbf{\hat{v}}_{l,\vartheta }=\Psi \left( \mathbf{\cdot },\vartheta 
\mathbf{k}_{\ast i_{l}}\right) \mathcal{F}_{\limfunc{av},n_{l},\vartheta
}\left( \mathbf{\vec{v}}\right) +\Psi \left( \mathbf{\cdot },\vartheta 
\mathbf{k}_{\ast i_{l}}\right) \Pi _{n_{l},\vartheta }\mathbf{\hat{h}},\
l=1,...N,\vartheta =\pm ,  \label{eqav}
\end{equation}%
which similarly to (\ref{syslocF}) we recast concisely as 
\begin{equation}
\mathbf{\vec{v}}=\mathcal{F}_{\text{$\limfunc{av}$},\Psi }\left( \mathbf{%
\vec{v}}\right) +\mathbf{\vec{h}}_{_{\Psi }}.  \label{eqavF}
\end{equation}%
The following lemma is analogous to Lemmas \ref{Lemma Flippr}, \ref{Lemma
Flip}.

\begin{lemma}
\label{Lemma Flipprav} Operator\ $\mathcal{F}_{\limfunc{av},\Psi }\left( 
\mathbf{\vec{v}}\right) $ is bounded \ for bounded $\mathbf{\vec{v}}\in
E^{2N}$, $\mathcal{F}_{\limfunc{av},\Psi }\left( \mathbf{0}\right) =\mathbf{0%
}$. Polynomial operator $\mathcal{F}_{\limfunc{av},\Psi }\left( \mathbf{\vec{%
v}}\right) $ satisfies the Lipschitz condition 
\begin{equation}
\left\Vert \mathcal{F}_{\text{$\limfunc{av}$},\Psi }\left( \mathbf{\vec{v}}%
_{1}\right) -\mathcal{F}_{\text{$\limfunc{av}$},\Psi }\left( \mathbf{\vec{v}}%
_{2}\right) \right\Vert _{E^{2N}}\leq C\tau _{\ast }\left\Vert \mathbf{\vec{v%
}}_{1}-\mathbf{\vec{v}}_{2}\right\Vert _{E^{2N}}
\end{equation}%
where $C$ depends only on $C_{\chi }$ a in (\ref{chiCR}), on the power of $%
\mathcal{F}$ and on $\left\Vert \mathbf{\vec{v}}_{1}\right\Vert
_{E^{2N}}+\left\Vert \mathbf{\vec{v}}_{2}\right\Vert _{E^{2N}}$, and, in
particular, it does not depend on $\beta $.
\end{lemma}

From Lemma \ref{Lemma Flipprav} and the contraction principle we obtain the
following Theorem similarly to Theorem \ref{Theorem existpr}.

\begin{theorem}
\label{Theorem existprav}Let $\left\Vert \mathbf{\vec{h}}_{\Psi }\right\Vert
_{E^{2N}}\leq R$. Then there exists $R_{1}>0$ and $\tau _{\ast }>0$ such
that equation (\ref{eqavF}) has a solution $\mathbf{\vec{v}}\in E^{2N}$
satisfying $\left\Vert \mathbf{\vec{v}}\right\Vert _{E^{2N}}\leq R_{1}$, and
such a solution is unique.
\end{theorem}

\begin{theorem}
\label{Theorem uminwav} Let $\mathbf{\hat{v}}_{l,\vartheta }\left( \mathbf{k}%
,\tau \right) \mathbf{\ }$be solution of (\ref{eqav}) and $\mathbf{\hat{w}}%
_{l,\vartheta }\left( \mathbf{k},\tau \right) \ $be the solution of \ (\ref%
{sysloc}). Then the $\mathbf{\hat{v}}_{l,\vartheta }\left( \mathbf{k},\tau
\right) $ is a wavepacket satisfying (\ref{piweq0}), (\ref{weq0}) with $%
\mathbf{\hat{w}}$ replaced by $\mathbf{\hat{v}}$. In addition to that, there
exists $\beta _{0}>0$ such that 
\begin{equation}
\left\Vert \mathbf{\hat{v}}_{l,\vartheta }\mathbf{\ }-\mathbf{\hat{w}}%
_{l,\vartheta }\right\Vert _{E}\leq C\varrho ,\ l=1,...,N;\ \vartheta =\pm ,%
\text{ for }0<\varrho \leq 1,\ 0<\beta \leq \beta _{0}.  \label{vminw}
\end{equation}
\end{theorem}

\begin{proof}
Formula (\ref{piweq0}), (\ref{weq0}) for $\mathbf{\hat{v}}_{l,\vartheta
}\left( \mathbf{k},\tau \right) $ \ follow from (\ref{eqav}). We note that $%
\mathbf{\vec{w}}$ is an approximate solution of (\ref{eqav}), namely we have
an estimate for $\mathbf{\hat{D}}_{\text{av}}\left( \mathbf{\hat{w}}\right) =%
\mathbf{\hat{w}}-\mathcal{F}_{\limfunc{av},\Psi }-\mathbf{\hat{h}}_{\Psi }$
which is similar to (\ref{Dw}), (\ref{Dwb}):%
\begin{equation}
\left\Vert \mathbf{\hat{D}}_{\text{av}}\left( \mathbf{\hat{w}}\right)
\right\Vert =\left\Vert \mathbf{\hat{w}}-\mathcal{F}_{\limfunc{av},\Psi }-%
\mathbf{\hat{h}}\right\Vert _{E}\leq C\varrho ,\text{ \ if }0<\varrho \leq
1,\beta \leq \beta _{0}.  \label{wav}
\end{equation}%
The proof of (\ref{wav}) is similar to the proof of (\ref{Fbet}) with minor
simplifications thanks to the absence of terms with $\Psi _{\infty }$. Using
(\ref{wav}) we apply Lemma \ref{Lemma contr} \ and obtain (\ref{vminw}).
\end{proof}

\subsection{Averaged system for scalar amplitudes}

Now we recast (\ref{eqav}) in the form of an equivalent system of scalar
equations for \emph{amplitudes} $\hat{v}_{l,\vartheta }=\hat{v}_{\lambda }$
of solutions $\mathbf{\hat{v}}_{\lambda _{l}}$ defined based on (\ref{Pin}),
namely%
\begin{equation}
\mathbf{\hat{v}}_{\lambda _{l}}\left( \mathbf{\mathbf{k}}\right) =\Psi
\left( \mathbf{k},\zeta ^{\left( l\right) }\mathbf{k}_{\ast i_{l}}\right)
\Pi _{n_{l},\zeta ^{\left( l\right) }}\left( \mathbf{\mathbf{k}}\right) 
\mathbf{\hat{v}}_{\lambda _{l}}\left( \mathbf{\mathbf{k}}\right) =\hat{v}%
_{l,\zeta ^{\left( l\right) }}\left( \mathbf{k}\right) \mathbf{g}%
_{n_{l},\zeta ^{\left( l\right) }}\left( \mathbf{k}\right) .  \label{wPi}
\end{equation}%
Note that according to (\ref{weq0}) support of $\hat{v}_{l,\zeta ^{\left(
l\right) }}$ is localized near $\zeta \mathbf{k}_{\ast i_{l}},$and we can
assume that $\mathbf{g}_{n_{l},\zeta ^{\left( l\right) }}\left( \mathbf{k}%
\right) $ depend smoothly on $\mathbf{k}$ near this point. Multiplying (\ref%
{eqav}) by $\mathbf{g}_{n_{l},\zeta _{l}}\left( \mathbf{k}\right) $ (with
the standard scalar product in $\mathbb{C}^{2j}$) and using (\ref{wPi}) we
obtain the following system of scalar amplitude equations 
\begin{gather}
\hat{v}_{l,\vartheta }=\Psi \left( \mathbf{\cdot },\vartheta \mathbf{k}%
_{\ast i_{l}}\right) f_{\limfunc{av},n_{l},\vartheta }\left( \vec{v}\right)
+\Psi \left( \mathbf{\cdot },\vartheta \mathbf{k}_{\ast i_{l}}\right) \hat{h}%
_{n_{l},\vartheta },\ l=1,...,N,\;\vartheta =\pm ,\text{ where}
\label{wscal} \\
\hat{h}_{n_{l},\vartheta }=\mathbf{g}_{n_{l},\vartheta }\cdot \Pi
_{n_{l},\vartheta }\mathbf{\hat{h}},\ f_{\limfunc{av},n_{l},\vartheta
}\left( \vec{v}\right) =\sum\nolimits_{m\in \mathfrak{M}_{F}}\sum\nolimits_{%
\vec{\lambda}\in \Lambda _{n_{l},\vartheta }^{m}}f_{n_{l},\vartheta ,\vec{\xi%
}\left( \vec{\lambda}\right) }^{\left( m\right) }\left( \vec{v}_{\vec{\lambda%
}}\right) .  \label{hhat}
\end{gather}%
According to (\ref{Fmw}) the $m$-linear operators in the above equation are
given by 
\begin{gather}
f_{n,\vartheta ,\vec{\xi}}^{\left( m\right) }\left( \vec{v}_{\vec{\lambda}%
}\right) \left( \mathbf{k},\tau \right) =\int_{0}^{\tau }\int_{\mathbb{D}%
_{m}}\mathrm{e}^{\mathrm{i}\phi _{n,\vartheta ,\vec{\xi}}\left( \mathbf{%
\mathbf{k}},\vec{k}\right) \frac{\tau _{1}}{\varrho }}Q_{n,\vartheta ,\vec{%
\xi}}^{\left( m\right) }\left( \mathbf{\mathbf{k}},\vec{k}\right)
\dprod\limits_{i=1}^{m}\hat{v}_{\lambda _{i}}\mathrm{\tilde{d}}^{\left(
m-1\right) d}\vec{k}\mathrm{d}\tau _{1},  \label{opfw} \\
Q_{n,\vartheta ,\vec{\xi}}^{\left( m\right) }\left( \mathbf{\mathbf{k}},\vec{%
k}\right) =\mathbf{g}_{n,\vartheta }\left( \mathbf{k}\right) \cdot \chi
_{n,\vartheta ,\vec{\xi}}^{\left( m\right) }\left( \mathbf{\mathbf{k}},\vec{k%
}\right) \left[ \mathbf{g}_{\lambda _{1}}\left( \mathbf{k}^{\prime }\right)
,\ldots ,\mathbf{g}_{\lambda _{m}}\left( \mathbf{k}^{\left( m\right) }\left( 
\mathbf{\mathbf{k}},\vec{k}\right) \right) \right] .  \label{Qg}
\end{gather}%
The concise form for the system (\ref{wscal}) of scalar equations for
amplitudes is 
\begin{equation}
\vec{v}=f_{\Psi }\left( \vec{v}\right) +\hat{h}_{\Psi },\text{ }\vec{v}\in
E_{\mathrm{sc}}^{2N},  \label{vpsih}
\end{equation}%
where the components $\hat{v}_{l,\vartheta }$ of $\vec{v}$ belong to the
space $E_{\mathrm{sc}}$ of scalar functions with the norm defined by (\ref%
{Elat}), (\ref{L1}) applied to scalar functions. Note that $Q_{n,\vartheta ,%
\vec{\xi}}^{\left( m\right) }\left( \mathbf{\mathbf{k}},\vec{k}\right) $ can
be extended in an arbitrary way \ as bounded functions for arguments $%
\mathbf{\mathbf{k}},\vec{k}$ where (\ref{konly}) is not satisfied, for
example the extension can be zero, the extension does not affect solutions
of (\ref{wscal}) because this equation involves factors $\Psi \left( \mathbf{%
\cdot },\vartheta \mathbf{k}_{\ast i_{l}}\right) $ and (\ref{weq0}) holds.

\begin{lemma}
\label{Lemma fscallip1} Operator\ $f_{\Psi }$ is bounded \ for bounded $\vec{%
v}\in E_{\mathrm{sc}}^{2N}$ and$\mathcal{\ }f_{\Psi }\left( \mathbf{0}%
\right) =\mathbf{0}$. Polynomial operator $f_{\Psi }\left( \vec{v}\right) $
satisfies the Lipschitz condition 
\begin{equation*}
\left\Vert f_{\Psi }\left( \vec{v}_{1}\right) -f_{\Psi }\left( \vec{v}%
_{2}\right) \right\Vert _{E_{\mathrm{sc}}^{2N}}\leq C\tau _{\ast }\left\Vert 
\vec{v}_{1}-\vec{v}_{2}\right\Vert _{E_{\mathrm{sc}}^{2N}}
\end{equation*}%
where $C$ depends only on $C_{\chi }$ as in (\ref{chiCR}), on the order of $%
\mathcal{F}$ as a polynomial and on $\left\Vert \vec{v}_{1}\right\Vert
_{E^{2N}}+\left\Vert \vec{v}_{2}\right\Vert _{E^{2N}}$, and it does not
depend on $\beta $.
\end{lemma}

From Lemma \ref{Lemma Flipprav} and the contraction principle we obtain the
following Theorem similarly to Theorem \ref{Theorem existpr}.

\begin{theorem}
\label{Theorem existprav1}Let $\left\Vert \hat{h}_{\Psi }\right\Vert _{E_{%
\mathrm{sc}}^{2N}}\leq R$. Then there exists $R_{1}>0$ and $\tau _{\ast }>0$
such that the (\ref{vpsih}) has a solution $\vec{v}\in E_{\mathrm{sc}}^{2N}$
\ satisfying $\left\Vert \vec{v}\right\Vert _{E_{\mathrm{sc}}^{2N}}\leq
R_{1} $, and such a solution is unique.
\end{theorem}

\subsection{Rescaled amplitude equations}

According to (\ref{weq0}) amplitudes $\hat{v}_{l,\vartheta }\left( \zeta 
\mathbf{k}_{\ast l}+\mathbf{\eta }\right) $ are localized about the point $%
\mathbf{\eta }=\mathbf{0}$, and to study its behavior in a vicinity of $%
\mathbf{\eta }=\mathbf{0}$ we introduce a group of dilation operators 
\begin{equation}
\left( B_{\beta }\hat{v}\right) \left( \mathbf{\eta }\right) =\beta ^{d}\hat{%
v}\left( \beta \mathbf{\eta }\right) ,\ \beta >0,  \label{Bbet}
\end{equation}%
which preserve the $L^{1}$-norm and commute with the convolution, i.e.%
\begin{equation}
\left\Vert B_{\beta }\hat{v}\right\Vert _{L^{1}}=\left\Vert \hat{v}%
\right\Vert _{L^{1}},\ B_{\beta }\hat{v}\ast B_{\beta }\hat{w}=B_{\beta
}\left( \hat{v}\ast \hat{w}\right) .  \label{BbetL}
\end{equation}%
We introduce then a rescaled and shifted version of initial data $\hat{h}%
_{n_{l},\vartheta }\ \ $in (\ref{hhat}) by the formula 
\begin{equation}
\hat{H}_{n_{l},\vartheta }\left( \mathbf{k}\right) =B_{\beta }\hat{h}%
_{n_{l},\vartheta }\left( \mathbf{k}+\vartheta \mathbf{k}_{\ast l}\right) ,\;%
\hat{h}_{n_{l},\vartheta }\left( \mathbf{k}\right) =\beta ^{-d}\hat{H}%
_{n_{l},\vartheta }\left( \beta ^{-1}\left( \mathbf{k}-\vartheta \mathbf{k}%
_{\ast l}\right) \right) ,  \label{hH}
\end{equation}%
where $B_{\beta }$ is defined by (\ref{Bbet}), $\left\vert \mathbf{k}%
-\vartheta \mathbf{k}_{\ast l}\right\vert \leq \beta ^{1-\epsilon }$, and
new variables 
\begin{equation}
\mathbf{\eta }_{l}=\beta ^{-1}\left( \mathbf{\mathbf{k}}-\vartheta \mathbf{%
\mathbf{k}}_{\ast l}\right) ,\ l=1,...,N,\ \vec{\eta}=\left( \mathbf{\eta }%
_{1},\ldots ,\mathbf{\eta }_{N}\right) .  \label{etak}
\end{equation}%
In this and the following sections we assume that $\hat{H}_{n_{l},\vartheta
}\left( \beta ,\mathbf{\eta }\right) $ are defined for all $\mathbf{\eta }%
\in \mathbb{R}^{d}$, including $\left\vert \mathbf{\eta }\right\vert \geq
\beta ^{-\epsilon }$. Though (\ref{wscal}) involves $\hat{h}%
_{n_{l},\vartheta }$ with a cutoff factor, namely \ $\Psi \left( \mathbf{k}%
,\vartheta \mathbf{k}_{\ast i_{l}}\right) \hat{h}_{n_{l},\vartheta }\left( 
\mathbf{k}\right) =\Psi \left( \mathbf{k},\vartheta \mathbf{k}_{\ast
i_{l}},\beta ^{1-\epsilon }\right) \hat{h}_{n_{l},\vartheta }\left( \mathbf{k%
}\right) $ \ as in (\ref{Psik}), we will later use $\hat{H}_{n_{l},\vartheta
}\left( \beta ,\mathbf{\eta }\right) $ defined for all $\mathbf{\eta }$, and
assume that%
\begin{equation}
\left\Vert \left( 1-\Psi \left( \beta ^{\epsilon }\mathbf{\eta }\right)
\right) \hat{H}_{n_{l},\vartheta }\left( \beta ,\mathbf{\eta }\right)
\right\Vert _{L^{1}}\leq C\beta ^{s},  \label{halfH}
\end{equation}%
where (i) $\Psi \left( \beta ^{\epsilon }\mathbf{\eta }\right) =\Psi \left( 
\mathbf{\eta },0,\beta ^{-\epsilon }\right) $ is as in (\ref{j0}), (\ref%
{Psik}); (ii) $\epsilon $ and $s$ are the same as in Definition \ref%
{dwavepack};\ (iii) condition (\ref{halfH}) is consistent with (\ref%
{halfwave}) and (\ref{sourloc}).

For a solution $\ \hat{v}_{l,\vartheta }\left( \mathbf{\mathbf{k}},\tau
\right) $ of (\ref{wscal}) using (\ref{weq0}) we introduce the following
functions%
\begin{equation}
\hat{z}_{l,\vartheta }\left( \mathbf{\eta },\tau \right) =\beta ^{d}\hat{v}%
_{l,\vartheta }\left( \vartheta \mathbf{k}_{\ast l}+\beta \mathbf{\eta }%
,\tau \right) ,\;\hat{z}_{l,\vartheta }\left( \mathbf{\eta },\tau \right)
=\Psi \left( \beta ^{\epsilon }\mathbf{\eta }\right) \hat{z}_{l,\vartheta
}\left( \mathbf{\eta },\tau \right) ,\ \mathbf{\eta }\in \mathbb{R}^{d},
\label{vz}
\end{equation}%
which satisfy a rescaled version of (\ref{wscal}) provided below. Note that
since $\left( \vec{n},\vec{\zeta}\right) =\vec{\lambda}\in \Lambda
_{n_{l},\vartheta }^{m}$ and the $nk$-spectrum $S$ is resonance invariant we
have $\varkappa _{m}\left( \vec{\lambda}\right) =\sum\nolimits_{i}\zeta
^{\left( i\right) }\mathbf{k}_{\ast l_{i}}=\zeta \mathbf{\mathbf{k}}_{\ast
l}=\vartheta \mathbf{\mathbf{k}}_{\ast l}$. Since $\mathbf{\mathbf{k}},\vec{k%
}$\ satisfy the convolution identity (\ref{conv}) the variables $\mathbf{%
\eta },\vec{\eta}$ defined by (\ref{etak}) satisfy similar identity as well,
namely 
\begin{equation}
\mathbf{\eta }=\sum\nolimits_{i=1}^{m}\mathbf{\eta }^{\left( i\right) },\ 
\mathbf{\eta }^{\left( m\right) }\left( \mathbf{\mathbf{k}},\vec{\eta}%
\right) =\mathbf{\eta }-\sum\nolimits_{i=1}^{m-1}\mathbf{\eta }^{\left(
i\right) }.  \label{conveta}
\end{equation}%
Change of variables (\ref{etak}) in the integral operator $f_{\limfunc{av}%
,n_{l},\vartheta }$ defined by (\ref{opfw}) yields the following \emph{%
amplitude system} for $z_{l,\vartheta }$ which is equivalent to (\ref{wscal}%
):\emph{\ }%
\begin{equation}
\hat{z}_{l,\vartheta }\left( \mathbf{\eta }\right) =\Psi \left( \beta
^{\epsilon }\mathbf{\eta }\right) f_{\limfunc{av},n_{l},\vartheta ,\beta
}\left( \vec{z}\right) \left( \mathbf{\eta }\right) +\Psi \left( \beta
^{\epsilon }\mathbf{\eta }\right) \hat{H}_{n_{l},\vartheta }\left( \mathbf{%
\eta }\right) ,\ l=1,...N,\ \vartheta =\pm .  \label{lavsys}
\end{equation}%
According to (\ref{Psilz}), (\ref{hhat}) \ and (\ref{opfw}) 
\begin{gather}
\Psi \left( \mathbf{k},\vartheta \mathbf{k}_{\ast i_{l}},\beta ^{1-\epsilon
}\right) =\Psi \left( \beta ^{\epsilon }\mathbf{\eta }\right) ,\ f_{\limfunc{%
av},n_{l},\vartheta ,\beta }\left( \vec{z}\right) =\sum\nolimits_{m\in 
\mathfrak{M}_{F}}f_{\limfunc{av},n_{l},\vartheta ,\beta }^{\left( m\right)
}\left( \vec{z}\right) ,  \label{psibet} \\
f_{\limfunc{av},n_{l},\vartheta ,\beta }^{\left( m\right) }\left( \vec{z}%
\right) =\sum\nolimits_{\vec{\lambda}\in \Lambda _{n_{l},\vartheta
}^{m}}f_{n_{l},\vartheta ,\vec{\xi}\left( \vec{\lambda}\right) ,\beta
}^{\left( m\right) }\left( \vec{z}_{\vec{\lambda}}\right) ,  \notag
\end{gather}%
\begin{gather}
f_{n,\vartheta ,\vec{\xi}\left( \vec{\lambda}\right) ,\beta }^{\left(
m\right) }\left( \vec{z}_{\vec{\lambda}}\right) \left( \mathbf{\eta },\tau
\right) =\int_{0}^{\tau }\int_{\mathbf{\eta }^{\prime }+...+\mathbf{\eta }%
^{\left( m\right) }=\mathbf{\eta }}\exp \left\{ \mathrm{i}\phi _{n,\vartheta
,\vec{\xi}\left( \vec{\lambda}\right) }\left( \vartheta \mathbf{\mathbf{k}}%
_{\ast l}+\beta \mathbf{\eta },\vec{k}_{\ast }+\beta \vec{\eta}\right) \frac{%
\tau _{1}}{\varrho }\right\}  \label{fbet} \\
Q_{n,\vartheta ,\vec{\xi}\left( \vec{\lambda}\right) }^{\left( m\right)
}\left( \vartheta \mathbf{\mathbf{k}}_{\ast l}+\beta \mathbf{\eta },\vec{k}%
_{\ast }+\beta \vec{\eta}\right) \dprod\nolimits_{i=1}^{m}\hat{z}_{\lambda
_{i}}\left( \mathbf{\eta }^{\left( i\right) }\right) \mathrm{\tilde{d}}%
^{\left( m-1\right) d}\vec{\eta}\mathrm{d}\tau _{1}.  \notag
\end{gather}%
Note that the condition (\ref{konly}) on the domain of integration takes in
the new variables the form 
\begin{equation}
\left\vert \mathbf{\eta }^{\left( i\right) }\right\vert \leq \beta
^{-\epsilon },\ i=1,...,m\text{ and }\left\vert \mathbf{\eta }\right\vert
\leq m\beta ^{-\epsilon }.  \label{konly1}
\end{equation}%
Finally, we rewrite \emph{amplitude system} (\ref{lavsys}) in the concise
form 
\begin{equation}
\vec{z}=\Psi \left( \beta ^{\epsilon }\mathbf{\cdot }\right) f_{\limfunc{av}%
,\beta }\left( \vec{z}\right) +\Psi \left( \beta ^{\epsilon }\mathbf{\cdot }%
\right) \hat{H}_{\beta },\;\vec{z}\in E_{\mathrm{sc}}^{2N}.  \label{lavsys1}
\end{equation}%
Let us show now that (\ref{lavsys1}) is of the form of (\ref{varcu}) with $%
2J $-component vector $\mathbf{\hat{u}}$ substituted with $2N$-component
vector $\vec{z}$, the matrix $\mathbf{L}\left( \mathbf{k}\right) $
substituted with a diagonal matrix $\vec{L}$ with entries $\vartheta \omega
_{n_{l}}\left( \vartheta \mathbf{\mathbf{k}}_{\ast l}+\beta \mathbf{\eta }%
\right) $. For that we introduce $S$-\emph{averaged tensor }$Q_{\text{av}%
}^{\left( m\right) }$ defined on $\vec{z}\in \mathbb{C}^{2Nm}$ by the
formula 
\begin{equation}
Q_{\text{av},n,\vartheta }^{\left( m\right) }\left( \beta \mathbf{\eta }%
,\beta \vec{\eta},\vec{z}\right) =\sum\nolimits_{\vec{\lambda}\in \Lambda
_{n,\vartheta }^{m}}Q_{n,\vartheta ,\vec{\xi}\left( \vec{\lambda}\right)
}^{\left( m\right) }\left( \vartheta \mathbf{\mathbf{k}}_{\ast l}+\beta 
\mathbf{\eta },\vec{k}_{\ast }+\beta \vec{\eta}\right)
\dprod\nolimits_{i=1}^{m}\hat{z}_{\lambda _{i}}  \label{resinvten}
\end{equation}%
which depends on $S$ through $\Lambda _{n,\vartheta }^{m}$ and acts from $%
\mathbb{C}^{2Nm}$ into $\mathbb{C}^{2N}$. Note that $\hat{z}_{\lambda _{i}}$
and $Q_{n,\vartheta ,\vec{\xi}}^{\left( m\right) }$ are scalar factors, $%
\hat{z}_{\lambda _{i}}$is a scalar projection in $\mathbb{C}^{2N}$ onto a
line along $\ \lambda _{i}$-th eigenvector of $\vec{L}$. Hence,\ the
right-hand side of (\ref{resinvten}) is a sum of elementary susceptibilities
obtained from $Q_{\text{av}}^{\left( m\right) }$ as in (\ref{chisumn}) and (%
\ref{lavsys}) has the form of (\ref{equfa}). \ Note that non-zero terms in (%
\ref{resinvten}) contain products $\hat{z}_{\lambda _{i}}$ which satisfy (%
\ref{Omeq0}). Therefore, if $\beta =0$ and $S$ is resonance invariant $Q_{%
\text{av}}^{\left( m\right) }$ has the form of weakly universal
nonlinearity; if $S$ is universally resonance invariant then $Q_{\text{av}%
}^{\left( m\right) }$ has the form of a universal nonlinearity as in (\ref%
{Ppm}).

\subsection{Amplitude\emph{\ }system with polynomial dispersion relations}

Now we introduce\emph{\ amplitude system} \emph{with polynomial dispersion }%
which is\emph{\ }similar to (\ref{lavsys}) and provides (i) sufficiently
accurate approximation to (\ref{lavsys}); (ii) standard polynomial
dependence of coefficients on $\mathbf{\eta },\vec{\eta}$ in the sense
clarified below. The amplitude system has the form 
\begin{gather}
\hat{u}_{l,\vartheta }=\Psi \left( \beta ^{\epsilon }\mathbf{\eta }\right)
f_{n_{l},\vartheta }^{\left( \mu ,\nu \right) }\left( \vec{u}\right) +\Psi
\left( \beta ^{\epsilon }\mathbf{\eta }\right) \hat{H}_{n_{l},\vartheta },\
l=1,...N,\;\vartheta =\pm ,  \label{sysu} \\
f_{n_{l},\vartheta }^{\left( \mu ,\nu \right) }\left( \vec{u}\right)
=\sum\nolimits_{m\in \mathfrak{M}_{F}}\sum\nolimits_{\vec{\lambda}\in
\Lambda _{n_{l},\vartheta }^{m}}f_{n_{l},\vartheta ,\vec{\xi}\left( \vec{%
\lambda}\right) }^{\left( m,\mu ,\nu \right) }\left( \vec{u}_{\vec{\lambda}%
}\right) ,  \label{fsnu}
\end{gather}%
where $\Psi \left( \beta ^{\epsilon }\mathbf{\eta }\right) $ are
cutoff-factors defined in (\ref{psibet}), (\ref{Psilz}) and approximations $%
f_{n_{l},\vartheta ,\vec{\xi}\left( \vec{\lambda}\right) }^{\left( m,\mu
,\nu \right) }$ for $f_{n_{l},\vartheta ,\vec{\xi}\left( \vec{\lambda}%
\right) }^{\left( m\right) }$ are defined below. The indices $\mu =1,2$, $%
\nu =0,1$ determine the order of approximation: (i) $\mu $ determines the
order of approximation of the dispersion relation by a polynomial of the
degree $\mu $; (ii) $\nu $ determines the order of approximation of\ the
susceptibility coefficients (\ref{Qg}) by a trigonometric polynomial of the
degree $\nu $. As before, we recast (\ref{sysu}) in a concise form 
\begin{equation}
\vec{u}=\Psi _{\beta }f^{\left( \mu ,\nu \right) }\left( \vec{u}\right)
+\Psi _{\beta }\hat{H}.  \label{eq}
\end{equation}%
where $\Psi _{\beta }\left( \mathbf{\eta }\right) =\Psi \left( \beta
^{\epsilon }\mathbf{\eta }\right) $. Finally, we eliminate \ in (\ref{sysu})
the cutoff factor $\Psi \left( \beta ^{\epsilon }\mathbf{\eta }\right) $ by
setting $\ \Psi \left( \beta ^{\epsilon }\mathbf{\eta }\right) =\Psi \left( 
\mathbf{0}\right) =1$, and introduce the amplitude system with \emph{\
weakly universal nonlinearity and polynomial dispersion} without cutoff 
\begin{equation}
\hat{u}_{l,\vartheta }\left( \mathbf{\eta }\right) =f_{n_{l},\vartheta
}^{\left( \mu ,\nu \right) }\left( \vec{u}\right) \left( \mathbf{\eta }%
\right) +\hat{H}_{n_{l},\vartheta }\left( \mathbf{\eta }\right) ,\
l=1,...N,\ \vartheta =\pm ,  \label{eqcan}
\end{equation}%
which can be written in the form of (\ref{eq}) with $\Psi _{\beta }=1$.

Let us turn now to the construction of the approximations. For every $nk$%
-pair $\left( \mathbf{k}_{\ast l},n_{l}\right) $ we introduce the Taylor
polynomials of order $\mu $ of the dispersion relation $\omega
_{n_{l}}\left( \mathbf{k}_{\ast l}+\beta \mathbf{\eta }\right) $: 
\begin{gather*}
\gamma _{1}\left( \mathbf{k}_{\ast l},n_{l},\beta \mathbf{\eta }\right)
=\omega _{n_{l}}\left( \mathbf{k}_{\ast l}\right) +\beta \omega
_{n_{l}}^{\prime }\left( \mathbf{k}_{\ast l}\right) \mathbf{\eta },\text{ \ }
\\
\gamma _{2}\left( \mathbf{k}_{\ast l},n_{l},\beta \mathbf{\eta }\right)
=\gamma _{1}\left( \mathbf{k}_{\ast l},n_{l},\beta \mathbf{\eta }\right) +%
\frac{\beta ^{2}}{2}\left( \mathbf{\eta },\omega _{n_{l}}^{\prime \prime
}\left( \mathbf{k}_{\ast l}\right) \mathbf{\eta }\right) ,
\end{gather*}%
and similarly $\gamma _{3}$ for $\mu =3$. Obviously we have the inequality
(see (\ref{konly}))%
\begin{equation}
\left\vert \omega _{n_{l}}\left( \mathbf{k}_{\ast l}+\beta \mathbf{\eta }%
\right) -\gamma _{\mu }\left( \mathbf{k}_{\ast l},n_{l},\beta \mathbf{\eta }%
\right) \right\vert \leq C\beta ^{\left( \mu +1\right) \left( 1-\epsilon
_{1}\right) },\ \left( \mathbf{\mathbf{k}},\vec{k}\right) \in B_{\beta }.
\label{omming}
\end{equation}%
The phase function $\phi _{n,\zeta ,\vec{\xi}}\left( \mathbf{\mathbf{k}},%
\vec{k}\right) $, $\vec{\xi}=\left( \vec{n},\vec{\zeta}\right) $, defined by
(\ref{phim}), is approximated then by a polynomial phase function \ 
\begin{gather}
\phi _{n_{l},\zeta ,\vec{\xi}}^{\left( \mu \right) }\left( \zeta \mathbf{%
\mathbf{k}}_{\ast l},\vec{k}_{\ast },\beta \mathbf{\eta },\beta \vec{\eta}%
\right) =  \label{fisn} \\
\zeta \gamma _{\mu }\left( \mathbf{k}_{\ast l},n_{l},\beta \mathbf{\eta }%
\right) -\zeta ^{\prime }\gamma _{\mu }\left( \mathbf{k}_{\ast
l_{1}},n^{\prime },\beta \mathbf{\eta }^{\prime }\right) -\ldots -\zeta
^{\left( m\right) }\gamma _{\mu }\left( \mathbf{k}_{\ast l_{m}},n^{\left(
m\right) },\beta \mathbf{\eta }^{\left( m\right) }\right) .  \notag
\end{gather}%
Note that since $\vec{\xi}=\vec{\xi}\left( \vec{\lambda}\right) $ with $\vec{%
\lambda}\in \Lambda _{n_{l},\vartheta }^{m}$ defined by (\ref{resset}),\ the
equation (\ref{Omeq0}) is fulfilled. Hence, $\phi _{n_{l},\vartheta ,\vec{\xi%
}}^{\left( \mu \right) }\left( \vartheta \mathbf{\mathbf{k}}_{\ast l},\vec{k}%
_{\ast },\mathbf{0},\mathbf{0}\right) =0$ and the function $\phi
_{n_{l},\vartheta ,\vec{\xi}}^{1}$ depends linearly on $\mathbf{\eta },\vec{%
\eta}$ and $\phi _{n_{l},\vartheta ,\vec{\xi}}^{2}$ is quadratic, namely%
\begin{equation}
\phi _{n_{l},\vartheta ,\vec{\xi}}^{1}\left( \vartheta \mathbf{\mathbf{k}}%
_{\ast l},\vec{k}_{\ast },\beta \mathbf{\eta },\beta \vec{\eta}\right)
=\beta \phi _{n_{l},\vartheta ,\vec{\xi}}^{1}\left( \vartheta \mathbf{%
\mathbf{k}}_{\ast l},\vec{k}_{\ast },\mathbf{\eta },\vec{\eta}\right) ,
\label{fihom1}
\end{equation}%
\begin{equation}
\phi _{n_{l},\vartheta ,\vec{\xi}}^{2}\left( \vartheta \mathbf{\mathbf{k}}%
_{\ast l},\vec{k}_{\ast },\beta \mathbf{\eta },\beta \vec{\eta}\right)
=\beta \phi _{n_{l},\vartheta ,\vec{\xi}}^{1}\left( \vartheta \mathbf{%
\mathbf{k}}_{\ast l},\vec{k}_{\ast },\mathbf{\eta },\vec{\eta}\right) +\beta
^{2}\phi _{n_{l},\vartheta ,\vec{\xi}}^{2,2}\left( \vartheta \mathbf{\mathbf{%
k}}_{\ast l},\vec{k}_{\ast },\mathbf{\eta },\vec{\eta}\right) .
\label{fihom}
\end{equation}%
In the case $\mu =2$ the polynomial phase function involves two parameters $%
\varrho _{1},\varrho _{2}$: 
\begin{gather}
\phi _{n_{l},\vartheta ,\vec{\xi}}^{2}\left( \vartheta \mathbf{\mathbf{k}}%
_{\ast l},\vec{k}_{\ast },\beta \mathbf{\eta },\beta \vec{\eta}\right) \frac{%
\tau _{1}}{\varrho }=\mathrm{i}\phi _{n_{l},\vartheta ,\vec{\xi}}^{1}\left(
\vartheta \mathbf{\mathbf{k}}_{\ast l},\vec{k}_{\ast },\mathbf{\eta },\vec{%
\eta}\right) \frac{\tau _{1}}{\varrho _{1}}+\mathrm{i}\phi _{n_{l},\vartheta
,\vec{\xi}}^{2,2}\left( \vartheta \mathbf{\mathbf{k}}_{\ast l},\vec{k}_{\ast
},\mathbf{\eta },\vec{\eta}\right) \frac{\tau _{1}}{\varrho _{2}},
\label{fir1r2} \\
\varrho _{1}=\frac{\varrho }{\beta },\ \varrho _{2}=\frac{\varrho }{\beta
^{2}};\ 0<\varrho _{1}<\infty ,\ 0<\varrho _{2}\leq \infty ,  \label{r1r2}
\end{gather}%
\ where $\varrho _{1}$ and $\varrho _{2}$ may be large or small depending on
relation between $\varrho $ and $\beta $. Sometimes it is convenient to
consider $\varrho _{1}$ and $\varrho _{2}$ as independent parameters. If $%
\mu =1$ we formally set $\varrho _{2}=\infty $, $\frac{\tau _{1}}{\varrho
_{2}}=0$. If (\ref{konly}) holds \ we have the estimate%
\begin{equation}
\left\vert \mathrm{e}^{\left\{ \mathrm{i}\phi _{n_{l},\vartheta ,\vec{\xi}%
}^{\mu }\left( \vartheta \mathbf{\mathbf{k}}_{\ast l},\vec{k}_{\ast },\beta 
\mathbf{\eta },\beta \vec{\eta}\right) \frac{\tau _{1}}{\varrho }\right\} }-%
\mathrm{e}^{\left\{ \mathrm{i}\phi _{n_{l},\vartheta ,\vec{\xi}}\left(
\vartheta \mathbf{\mathbf{k}}_{\ast l}+\beta \mathbf{\eta },\vec{k}_{\ast
}+\beta \vec{\eta}\right) \frac{\tau _{1}}{\varrho }\right\} }\right\vert
\leq C\tau _{\ast }\frac{\beta ^{\left( \mu +1\right) \left( 1-\epsilon
\right) }}{\varrho },\;\mu =1,2.  \label{fifib}
\end{equation}%
To ensure that the approximation error is small for given $\mu $ we assume
that $\varrho $ and $\beta $ satisfy 
\begin{equation}
\varrho \rightarrow 0,\,\beta \rightarrow 0,\ \frac{\beta ^{\left( \mu
+1\right) \left( 1-\epsilon \right) }}{\varrho }\rightarrow 0.  \label{kap13}
\end{equation}%
Now we approximate the dependence of $Q_{n,\zeta ,\vec{\xi}}^{\left(
m\right) }\left( \vartheta \mathbf{\mathbf{k}}_{\ast l}+\beta \mathbf{\eta },%
\vec{k}_{\ast }+\beta \vec{\eta}\right) $ on $\mathbf{\eta }$,$\ \vec{\eta}$
given by (\ref{Qg}) by trigonometric polynomials. Zero order approximation
with $\nu =0$ is given by 
\begin{equation}
Q_{n,\zeta ,\vec{\xi}}^{\left( m,0\right) }\left( \vartheta \mathbf{\mathbf{k%
}}_{\ast l}+\beta \mathbf{\eta },\vec{k}_{\ast }+\beta \vec{\eta}\right)
=Q_{n,\zeta ,\vec{\xi}}^{\left( m\right) }\left( \vartheta \mathbf{\mathbf{k}%
}_{\ast l},\vec{k}_{\ast }\right) .  \label{Qm0}
\end{equation}%
To define the first order approximation we modify the standard Taylor
expansion \ using trigonometric polynomials instead of algebraic ones.
Taking the first derivative with respect to $\beta $ at $\beta =0$ 
\begin{equation*}
Q_{n,\zeta ,\vec{\xi}}^{\left( m\right) \prime }\left( \vartheta \mathbf{%
\mathbf{k}}_{\ast l},\mathbf{\eta },\vec{k}_{\ast },\vec{\eta}\right)
=\left. \frac{d}{d\beta }\right\vert _{\beta =0}Q_{n,\zeta ,\vec{\xi}%
}^{\left( m\right) }\left( \vartheta \mathbf{\mathbf{k}}_{\ast l}+\beta 
\mathbf{\eta },\vec{k}_{\ast }+\beta \vec{\eta}\right) ,
\end{equation*}%
which obviously is a linear function with respect to $\mathbf{\eta }$, $\vec{%
\eta}$, we express then $\mathbf{\eta }$ in terms of $\vec{\eta}$ using (\ref%
{conveta}): 
\begin{equation*}
Q_{n,\zeta ,\vec{\xi}}^{\left( m\right) \prime }\left( \vartheta \mathbf{%
\mathbf{k}}_{\ast l},\mathbf{\eta },\vec{k}_{\ast },\vec{\eta}\right)
=\sum\nolimits_{j=1}^{m}q_{n,\zeta ,\vec{\xi}}^{\left( m\right) ,j}\left(
\vartheta \mathbf{\mathbf{k}}_{\ast l},\vec{k}_{\ast }\right) \cdot \mathbf{%
\eta }^{\left( j\right) },\ \mathbf{\eta }^{\left( j\right) }=\left( \eta
_{1}^{\left( j\right) },...,\eta _{d}^{\left( j\right) }\right) .
\end{equation*}%
Then the first order approximation is 
\begin{equation}
Q_{n,\zeta ,\vec{\xi}}^{\left( m,1\right) }\left( \vartheta \mathbf{\mathbf{k%
}}_{\ast l}+\beta \mathbf{\eta },\vec{k}_{\ast }+\beta \vec{\eta}\right)
=Q_{n,\zeta ,\vec{\xi}}^{\left( m\right) }\left( \vartheta \mathbf{\mathbf{k}%
}_{\ast l},\vec{k}_{\ast }\right) +\sum_{j=1}^{m}q_{n,\zeta ,\vec{\xi}%
}^{\left( m\right) ,j}\left( \vartheta \mathbf{\mathbf{k}}_{\ast l},\vec{k}%
_{\ast }\right) \cdot \sin \beta \mathbf{\eta }^{\left( j\right) },  \notag
\end{equation}%
where $\sin \mathbf{\eta }^{\left( j\right) }=\left( \sin \eta _{1}^{\left(
j\right) },...,\sin \eta _{d}^{\left( j\right) }\right) $. An advantage of
this approximation is that the multiplication by $\sin \eta _{1}^{\left(
j\right) }$ is a bounded operator which equals the Fourier transform of a
finite-difference operator whereas the multiplication by $\eta _{1}^{\left(
j\right) }$ corresponds to the partial derivative and is unbounded. Since
the original nonlinearity does not involve unbounded operators the use of
bounded operators is natural and convenient. In fact, it is well known that
the presence of the derivatives in the nonlinearity of NLS-type equations
causes well known technical difficulties, see \cite{Bourgain}. In our
approach the approximating equation provides the same accuracy and its
nonlinearity involves only bounded finite-difference operators bypassing
those difficulties altogether.

According to Condition \ref{cnonbound} the susceptibility is smooth and if (%
\ref{konly1}) holds we have the following inequality 
\begin{equation}
\left\vert Q_{n,\zeta ,\vec{\xi}}^{\left( m\right) }\left( \vartheta \mathbf{%
\mathbf{k}}_{\ast l}+\beta \mathbf{\eta },\vec{k}_{\ast }+\beta \vec{\eta}%
\right) -Q_{n,\zeta ,\vec{\xi}}^{\left( m,\nu \right) }\left( \vartheta 
\mathbf{\mathbf{k}}_{\ast l},\vec{k}_{\ast },\beta \mathbf{\eta },\beta \vec{%
\eta}\right) \right\vert \leq C\beta ^{\left( \nu +1\right) \left(
1-\epsilon _{1}\right) }.  \label{chistar0}
\end{equation}%
We introduce components $f_{n_{l},\vartheta ,\vec{\lambda}}^{\left( m,\mu
,\nu \right) }$ of the weakly universal nonlinearity $f^{\left( \mu ,\nu
\right) }$ by the formula\ \ 
\begin{gather}
f_{n_{l},\vartheta ,\vec{\lambda}}^{\left( m,\mu ,\nu \right) }\left( \vec{z}%
_{\vec{\lambda}}\right) \left( \mathbf{\eta },\tau \right) =\int_{0}^{\tau
}\int_{\mathbf{\eta }^{\prime }+...+\mathbf{\eta }^{\left( m\right) }=%
\mathbf{\eta }}\mathrm{e}^{\mathrm{i}\phi _{n_{l},\vartheta ,\vec{\xi}%
}^{1}\left( \vartheta \mathbf{\mathbf{k}}_{\ast l},\vec{k}_{\ast },\mathbf{%
\eta },\vec{\eta}\right) \frac{\tau _{1}}{\varrho _{1}}+\mathrm{i}\phi
_{n_{l},\vartheta ,\vec{\xi}}^{2,2}\left( \vartheta \mathbf{\mathbf{k}}%
_{\ast l},\vec{k}_{\ast },\mathbf{\eta },\vec{\eta}\right) \frac{\tau _{1}}{%
\varrho _{2}}}  \label{fbet0} \\
Q_{n_{l},\vartheta ,\vec{\xi}}^{\left( m,\nu \right) }\left( \vartheta 
\mathbf{\mathbf{k}}_{\ast l},\vec{k}_{\ast }\right) \dprod\limits_{i=1}^{m}%
\hat{z}_{\lambda _{i}}\left( \mathbf{\eta }^{\left( i\right) }\right) 
\mathrm{\tilde{d}}^{\left( m-1\right) d}\vec{k}\,\mathrm{d}\tau _{1}.  \notag
\end{gather}%
As before, we establish standard properties of the operator $f^{\left( \mu
,\nu \right) }$ defined by the above formula.

\begin{lemma}
\label{Lemma lfscallip} Operator\ \ $\Psi _{\beta }f^{\left( \mu ,\nu
\right) }$ is bounded \ for bounded $\vec{u}\in E_{\mathrm{sc}}^{2N}$, $%
f_{\Psi }\left( \mathbf{0}\right) =\mathbf{0}$. Polynomial operator $\Psi
_{\beta }f^{\left( \mu ,\nu \right) }$ satisfies the Lipschitz condition 
\begin{equation}
\left\Vert \Psi _{\beta }f^{\left( \mu ,\nu \right) }\left( \vec{u}%
_{1}\right) -\Psi _{\beta }f^{\left( \mu ,\nu \right) }\left( \vec{u}%
_{2}\right) \right\Vert _{E_{\mathrm{sc}}^{2N}}\leq C\tau _{\ast }\left\Vert 
\vec{u}_{1}-\vec{u}_{2}\right\Vert _{E_{\mathrm{sc}}^{2N}},
\end{equation}%
where $C$ depends only on $C_{\chi }$ a in (\ref{chiCR}), on the power of $%
\mathcal{F}$ and on $\left\Vert \vec{u}_{1}\right\Vert _{E_{\mathrm{sc}%
}^{2N}}+\left\Vert \vec{u}_{2}\right\Vert _{E_{\mathrm{sc}}^{2N}}$. In
particular, it does not depend on $\beta \geq 0$ and on $0<\varrho
_{1}<\infty $,\ $0<\varrho _{2}\leq \infty $.
\end{lemma}

From Lemma \ref{Lemma Flipprav} \ and the contraction principle we obtain
the following Theorem completely similar to Theorem \ref{Theorem existpr}.

\begin{theorem}
\label{Theorem existpravl}Let $\left\Vert \hat{h}_{\Psi }\right\Vert _{E_{%
\mathrm{sc}}^{2N}}\leq R.$ Then there exists $R_{1}>0$ and $\tau _{\ast }>0$
such that equation (\ref{eqavF}) has a solution $\vec{z}\in E_{\mathrm{sc}%
}^{2N}$ satisfying $\left\Vert \vec{z}\right\Vert _{E_{\mathrm{sc}%
}^{2N}}\leq R_{1}$. Such a solution is unique and $\hat{z}_{l,\vartheta
}\left( \mathbf{k},\tau \right) =0$\ if $\left\vert \mathbf{k}\right\vert
\geq \beta ^{-\epsilon }$.
\end{theorem}

\begin{theorem}
\label{Theorem uminz} Let $\hat{u}_{l,\vartheta }\left( \mathbf{k},\tau
\right) \mathbf{\ }$be a solution to\ (\ref{sysu}) and $\hat{z}_{l,\vartheta
}\left( \mathbf{k},\tau \right) \ $be the solution of \ (\ref{lavsys1}).
Then the following inequality holds: 
\begin{equation}
\left\Vert \hat{u}_{l,\vartheta }-\hat{z}_{l,\vartheta }\right\Vert _{E_{%
\mathrm{sc}}}\leq C\beta ^{\left( \mu +1\right) \left( 1-\epsilon \right)
}+C\varrho ^{-1}\beta ^{\left( \mu +1\right) \left( 1-\epsilon \right) },\
l=1,...,N;\ \vartheta =\pm ,  \label{vminz}
\end{equation}%
for all $0<\varrho \leq $ $1$ and $0<\beta \leq \beta _{0}$, where $\epsilon 
$ is the same as in Definition \ref{dwavepack}, $\beta _{0}$ is sufficiently
small.
\end{theorem}

\begin{proof}
To obtain (\ref{vminz}) we note that $u_{l,\vartheta }$ is an approximate
solution of (\ref{lavsys1}), namely 
\begin{equation*}
\vec{u}-\Psi _{\beta }f^{\left( \mu ,\nu \right) }\left( \vec{u}\right) -%
\hat{h}_{\Psi }=\hat{D}\text{ where }\hat{D}\text{ is small.}
\end{equation*}%
To estimate $\left\Vert \hat{D}\right\Vert $ observe that integrals
involving $\vec{u}$ have the integration domain as in (\ref{konly}). Hence,
using (\ref{chistar0}) and (\ref{fifib}) we obtain 
\begin{equation*}
\left\Vert \hat{D}\right\Vert _{E_{\mathrm{sc}}^{2N}}\leq C\beta ^{\left(
\mu +1\right) \left( 1-\epsilon \right) }+C\varrho ^{-1}\beta ^{\left( \mu
+1\right) \left( 1-\epsilon \right) },
\end{equation*}%
and applying Lemma \ref{Lemma contr} we get (\ref{vminz}).
\end{proof}

\subsection{Decay of solutions and elimination of cutoff factors}

In this subsection we show how to remove the cutoff function in (\ref{sysu})
and to obtain the averaged interaction system with a weakly universal
nonlinearity. If $\mu =1$, $\nu =0$ and the $nk$-spectrum $S$ is
resonance-invariant, the amplitude system coincides with the system (\ref%
{cansystext}) with a weakly universal nonlinearity. For $\mu >1$ or $\nu >0$
the amplitude system involves additional terms. In particular, if $\mu =2$, $%
\nu =0$ and $S=\left\{ \left( \mathbf{k}_{\ast },n\right) \right\} $ is just
a single element then the linear part has the second order and the
nonlinearity is universal, and amplitude system turns into the classical NLS
system: 
\begin{equation}
\partial _{\tau }u_{\zeta }=\zeta \frac{1}{\varrho }\gamma _{2}\left( 
\mathbf{k}_{\ast },n,-i\zeta \beta \nabla _{r}\mathbf{\eta }\right)
+b_{\zeta }u_{\zeta }^{2}u_{-\zeta },\ u_{\zeta }\left( 0\right) =\hat{H}%
_{\zeta },\ \zeta =\pm .  \notag
\end{equation}%
This system is equivalent to (\ref{NLS0}) when $\hat{H}_{-}=\hat{H}%
_{+}^{\ast }$, $b_{-}=b_{+}^{\ast }$, $u_{-}=u_{+}^{\ast }$. When $\nu >0$
the nonlinearity involves additional terms with finite difference operators.

The possibility to remove cutoff functions is based on the fast decay of $%
\hat{u}\left( \mathbf{k}\right) $ as $\left\vert \mathbf{k}\right\vert
\rightarrow \infty $ which is equivalent to high smoothness of $u\left( 
\mathbf{r}\right) $. The factor $\Psi _{\beta }$ can be replaced by $1$ with
a small error when data $\hat{H}\left( \mathbf{k}\right) $ decay
sufficiently fast. To to describe the decay we introduce weighted Banach
spaces of scalar functions $\hat{H}\left( \mathbf{k}\right) $ described as
follows.

\begin{definition}[weight function]
\label{Definition weight}For $a\geq 0$ we call a positive function $\psi
\left( r\right) $, $r\geq 0$, a weight function from class $W\left( a\right) 
$ if it satisfies the following conditions: (i) $\psi \left( 0\right) >0$, $%
\psi \left( r_{1}\right) \geq \psi \left( r_{2}\right) $ for $r_{1}\geq
r_{2}\geq 0$; (ii) $\psi \left( r_{1}+r_{2}\right) \leq \psi \left(
r_{1}\right) +\psi \left( r_{2}\right) +C$ where $C$ does not depend on $%
r_{1},r_{2}$ ($\psi $ is sublinear); (iii) $\psi \left( r\right) -a\ln r\geq
C^{\prime }>0$ \ for all $r>0$ ($\psi \left( r\right) $ is superlogarithmic).
\end{definition}

We introduce $L^{1}\left( \psi \right) $\ as a space of scalar functions $%
\hat{H}\left( \mathbf{k}\right) ,$ $\mathbf{k}\in \mathbb{R}^{d}$ with the
norm 
\begin{equation}
\left\Vert \hat{H}\right\Vert _{L^{1}\left( \psi \right) }=\int_{\mathbb{R}%
^{d}}\mathrm{e}^{\psi \left( \left\vert \mathbf{k}\right\vert \right)
}\left\vert \hat{H}\left( \mathbf{k}\right) \right\vert \,\mathrm{d}\mathbf{k%
}.  \label{L1psi}
\end{equation}%
For vector-functions we use the same formula with Euclidean norm $\left\vert
\cdot \right\vert $. In the simplest case of $\psi \left( r\right) =a\ln
\left( 1+r\right) $ we have $\psi \in W\left( a\right) $ and obtain $%
L^{1}\left( \psi \right) =L^{1,a}$ with the norm (\ref{L1a}). If the weight
function belongs to $W\left( a\right) $ for all $a$ the space $L^{1}\left(
\psi \right) $ consists of the Fourier transforms of infinitely smooth
functions. The following Lemma shows that $L^{1}\left( \psi \right) $ is
closed with respect to the convolution.

\begin{lemma}
\label{Lemma convpsi} Let $\hat{H}_{1},\hat{H}_{2}\in L^{1}\left( \psi
\right) $ and 
\begin{equation*}
\hat{H}_{3}\left( \mathbf{k}\right) =\int_{\mathbb{R}^{d}}\hat{H}_{1}\left( 
\mathbf{k}-\mathbf{k}^{\prime }\right) \hat{H}_{2}\left( \mathbf{k}-\mathbf{k%
}^{\prime }\right) \,\mathrm{d}\mathbf{k}^{\prime }.
\end{equation*}%
\begin{equation}
\text{Then }\left\Vert \hat{H}_{3}\left( \mathbf{k}\right) \right\Vert
_{L^{1}\left( \psi \right) }\leq C\left\Vert \hat{H}_{1}\left( \mathbf{k}%
\right) \right\Vert _{L^{1}\left( \psi \right) }\left\Vert \hat{H}_{1}\left( 
\mathbf{k}\right) \right\Vert _{L^{1}\left( \psi \right) }.  \label{L1conpsi}
\end{equation}
\end{lemma}

\begin{proof}
Using Definition \ref{Definition weight} (ii) we obtain 
\begin{gather*}
\mathrm{e}^{\psi \left( \left\vert \mathbf{k}\right\vert \right) }\left\vert 
\hat{H}_{3}\left( \mathbf{k}\right) \right\vert \leq \int_{\mathbb{R}^{d}}%
\mathrm{e}^{\psi \left( \left\vert \mathbf{k}\right\vert \right) }\left\vert 
\hat{H}_{1}\left( \mathbf{k}-\mathbf{k}^{\prime }\right) \right\vert
\left\vert \hat{H}_{2}\left( \mathbf{k}^{\prime }\right) \right\vert \,%
\mathrm{d}\mathbf{k}^{\prime }\leq \\
\mathrm{e}^{C}\int_{\mathbb{R}^{d}}\mathrm{e}^{\psi \left( \left\vert 
\mathbf{k}^{\prime }\right\vert \right) }\mathrm{e}^{\psi \left( \left\vert 
\mathbf{k-k}^{\prime }\right\vert \right) }\left\vert \hat{H}_{1}\left( 
\mathbf{k}-\mathbf{k}^{\prime }\right) \right\vert \left\vert \hat{H}%
_{2}\left( \mathbf{k}^{\prime }\right) \right\vert \,\mathrm{d}\mathbf{k}%
^{\prime }.
\end{gather*}%
Applying Young's inequality (\ref{Yconv}) we obtain 
\begin{equation*}
\int_{\mathbb{R}^{d}}\mathrm{e}^{\psi \left( \left\vert \mathbf{k}%
\right\vert \right) }\left\vert \hat{H}_{3}\left( \mathbf{k}\right)
\right\vert \mathrm{d}\mathbf{k}\leq \mathrm{e}^{C}\int_{\mathbb{R}^{d}}%
\mathrm{e}^{\psi \left( \left\vert \mathbf{k}\right\vert \right) }\left\vert 
\hat{H}_{1}\left( \mathbf{k}\right) \right\vert \mathrm{d}\mathbf{k}^{\prime
}\int_{\mathbb{R}^{d}}\mathrm{e}^{\psi \left( \left\vert \mathbf{k}%
\right\vert \right) }\left\vert \hat{H}_{2}\left( \mathbf{k}\right)
\right\vert \,\mathrm{d}\mathbf{k}^{\prime }
\end{equation*}%
implying (\ref{L1conpsi}).
\end{proof}

Let us introduce the norm in the space $E_{\mathrm{sc}}\left( \psi \right) $
by the formula (\ref{Elat}) 
\begin{equation}
\left\Vert \hat{H}\left( \cdot ,\cdot \right) \right\Vert _{E\left( \psi
\right) }=\left\Vert \hat{H}\left( \cdot ,\cdot \right) \right\Vert
_{C\left( \left[ 0,\tau _{\ast }\right] ,L^{1}\left( \psi \right) \right)
}=\sup_{0\leq \tau \leq \tau _{\ast }}\int_{\mathbb{R}^{d}}\mathrm{e}^{\psi
\left( \left\vert \mathbf{k}\right\vert \right) }\left\vert \hat{H}\left( 
\mathbf{k},\tau \right) \right\vert \,\mathrm{d}\mathbf{k}.  \label{Epsi}
\end{equation}%
Using (\ref{L1conpsi}) instead of (\ref{L1}) we obtain as in Lemma \ref%
{Lemma bound} the following statement.

\begin{lemma}
\label{Lemma fscallip} Operator\ \ $\Psi _{\beta }f^{\left( s,\nu \right) }$
in (\ref{eq}) is bounded \ for bounded $\vec{u}\in E_{\mathrm{sc}%
}^{2N}\left( \psi \right) $, $f\left( \mathbf{0}\right) =\mathbf{0}$ \ and
satisfies Lipschitz condition 
\begin{equation}
\left\Vert \Psi _{\beta }f^{\left( s,\nu \right) }\left( \vec{u}_{1}\right)
-\Psi _{\beta }f^{\left( s,\nu \right) }\left( \vec{u}_{2}\right)
\right\Vert _{E_{\mathrm{sc}}^{2N}\left( \psi \right) }\leq C\tau _{\ast
}\left\Vert \vec{u}_{1}-\vec{u}_{2}\right\Vert _{E_{\mathrm{sc}}^{2N}\left(
\psi \right) }
\end{equation}%
where $C$ depends only on $C_{\chi }$ a in (\ref{chiCR}), on the power of \
polynomial $\ f^{\left( s,\nu \right) }$ and on $\left\Vert \vec{u}%
_{1}\right\Vert _{E_{\mathrm{sc}}^{2N}\left( \psi \right) }+\left\Vert \vec{u%
}_{1}\right\Vert _{E_{\mathrm{sc}}^{2N}\left( \psi \right) }$ \ and does not
depend on $\beta \geq 0$ and on $0<\varrho _{1}<\infty $,\ $0<\varrho
_{2}\leq \infty $ .
\end{lemma}

From Lemma \ref{Lemma Flipprav} \ and the contraction principle we obtain
the following Theorem completely similar to Theorem \ref{Theorem existpr}.

\begin{theorem}
\label{Theorem existcan}Let $\left\Vert \hat{H}\right\Vert _{E_{\mathrm{sc}%
}^{2N}\left( \psi \right) }\leq R$. Then there exists $R_{1}>0$ and $\tau
_{\ast }>0$ such that equation (\ref{eq}) has a solution $\vec{u}\in E_{%
\mathrm{sc}}^{2N}\left( \psi \right) $ \ which satisfies $\left\Vert \vec{u}%
\right\Vert _{E_{\mathrm{sc}}^{2N}\left( \psi \right) }\leq R_{1}$, and such
a solution is unique.
\end{theorem}

The following lemma shows that $\Psi $ can be replaced by one with a small
error.

\begin{lemma}
\label{Lemma Lpsiwp}Let $\left\Vert \hat{H}\right\Vert _{L^{1}\left( \psi
\right) }\leq C$, $\psi \in W\left( a\right) ,\ \Psi $ as in (\ref{j0}). If $%
s>0,\;\epsilon >0$ and $\frac{s}{\epsilon }<a$ \ then (\ref{halfH}) holds.
\end{lemma}

\begin{proof}
We have 
\begin{gather}
\int \left( 1-\Psi \left( \beta ^{\epsilon }\mathbf{\eta }\right) \right)
\left\vert \hat{H}\left( \mathbf{\eta }\right) \right\vert \mathrm{d}\mathbf{%
\eta }\leq \int_{\left\vert \mathbf{\eta }\right\vert \geq \beta ^{-\epsilon
}}\left\vert \hat{H}\left( \mathbf{\eta }\right) \right\vert \mathrm{d}%
\mathbf{\eta }=\int_{\left\vert \mathbf{\eta }\right\vert \geq \beta
^{-\epsilon }}\mathrm{e}^{-\psi \left( \left\vert \mathbf{\eta }\right\vert
\right) }\left\vert \mathrm{e}^{\psi \left( \left\vert \mathbf{\eta }%
\right\vert \right) }\hat{H}\left( \mathbf{\eta }\right) \right\vert \mathrm{%
d}\mathbf{\eta }  \label{ineq1} \\
\leq \int_{\left\vert \mathbf{\eta }\right\vert \geq \beta ^{-\epsilon }}%
\mathrm{e}^{-\psi \left( \beta ^{-\epsilon }\right) }\left\vert \mathrm{e}%
^{\psi \left( \left\vert \mathbf{k}\right\vert \right) }\hat{H}\left( 
\mathbf{\eta }\right) \right\vert \mathrm{d}\mathbf{\eta }\leq \beta ^{s}%
\mathrm{e}^{\ln \left( \beta ^{-\epsilon }\right) s/\epsilon -\psi \left(
\beta ^{-\epsilon }\right) }\left\Vert \hat{H}\right\Vert _{L^{1}\left( \psi
\right) }.  \notag
\end{gather}%
According to Definition \ref{Definition weight} (iii) 
\begin{equation*}
\ln \left( \beta ^{-\epsilon }\right) s/\epsilon -\psi \left( \beta
^{-\epsilon }\right) \leq a\ln \left( \beta ^{-\epsilon }\right) -\psi
\left( \beta ^{-\epsilon }\right) \leq C
\end{equation*}%
and we obtain (\ref{halfH}) from (\ref{ineq1}).
\end{proof}

\begin{theorem}
\label{Theorem moll}Let $\left\Vert \hat{H}\right\Vert _{E_{\mathrm{sc}%
}^{2N}\left( \psi \right) }\leq R$ where weight function $\psi $ belongs to $%
W\left( a\right) $ and let $\frac{s}{\epsilon }<a$. Let $\vec{u}$ \ and $%
\vec{u}_{0}$ be solutions to respectively the minimal equation with cutoff
factor \ and without cutoff factor respectively. Then\ there exists $C_{s}$
and $\beta _{0}$\ such that 
\begin{equation}
\left\Vert \vec{u}-\vec{u}_{0}\right\Vert _{E_{\mathrm{sc}}^{2N}\left( \psi
\right) }\leq C_{s}\beta ^{s},\ 0<\beta \leq \beta _{0}.  \label{uu0}
\end{equation}
\end{theorem}

\begin{proof}
We show that $\vec{u}$ is an approximate solution to $\vec{u}_{0}=f^{\left(
\mu ,\nu \right) }\left( \vec{u}_{0}\right) +\hat{H}$. Namely, 
\begin{equation*}
\vec{u}=\Psi _{\beta }f^{\left( \mu ,\nu \right) }\left( \vec{u}\right)
+\Psi _{\beta }\hat{H}=f^{\left( \mu ,\nu \right) }\left( \vec{u}\right) +%
\hat{H}+\hat{D},\ \hat{D}=\left( \Psi _{\beta }-1\right) f^{\left( \mu ,\nu
\right) }\left( \vec{u}\right) +\left( \Psi _{\beta }-1\right) \hat{H}.
\end{equation*}%
According to Lemma \ref{Lemma convpsi} if $\vec{u}\in E_{\mathrm{sc}%
}^{2N}\left( \psi \right) $ then $f^{\left( \mu ,\nu \right) }\left( \vec{u}%
\right) \in E_{\mathrm{sc}}^{2N}\left( \psi \right) $. \ Applying Lemma \ref%
{Lemma Lpsiwp} we obtain 
\begin{equation}
\left\Vert \hat{D}\right\Vert _{E_{\mathrm{sc}}^{2N}\left( \psi \right)
}\leq C\beta ^{s},\ 0<\beta \leq \beta _{0}.  \label{Dbp}
\end{equation}%
Lemma \ref{Lemma contr} combined with (\ref{Dbp}) yields (\ref{uu0}).
\end{proof}

Now we give the theorem on approximation by solutions of a minimal system
without cutoff.

\begin{theorem}
\label{Theorem canappr}Let $\hat{H}_{l,\zeta }\left( \mathbf{k}\right) $, $%
l=1,...,N$ \ \ be functions bounded in $L^{1}\left( \psi \right) $ where $%
\psi $ belongs to $W\left( a\right) $, let $\frac{s}{\epsilon }<a$. Let $\ 
\hat{h}_{l,\zeta }\left( \mathbf{k}\right) $ be defined by (\ref{hH}) and $%
\Psi \mathbf{\hat{h}}_{l,\zeta }\left( \mathbf{k}\right) =\Psi \hat{h}%
_{l,\zeta }\left( \mathbf{k}\right) \mathbf{g}_{n_{l},\zeta }\left( \mathbf{k%
}\right) $. Let $\mathbf{\hat{u}}\left( \mathbf{k},\tau \right) \ $be a
solution of equation (\ref{varcu}) with multiwavepacket initial data of the
form (\ref{JJ1}). Let $u_{l,\vartheta }\left( \mathbf{k},\tau \right) $ be a
solution to the system with a weakly universal nonlinearity (\ref{eqcan})
with initial data $u_{l,\vartheta }\left( \mathbf{k},0\right) =\hat{H}%
_{l,\vartheta }\left( \mathbf{k}\right) $ and 
\begin{equation*}
\mathbf{\hat{u}}_{\min }\left( \mathbf{k},\tau \right)
=\sum\nolimits_{\vartheta }\sum\nolimits_{l=1}^{N}\beta ^{-d}u_{l,\vartheta
}\left( \beta ^{-1}\left( \mathbf{k}-\zeta \mathbf{k}_{\ast i_{l}}\right)
,\tau \right) \mathbf{g}_{n_{l},\vartheta }\left( \mathbf{k}\right) .
\end{equation*}%
Then 
\begin{equation}
\left\Vert \mathbf{\hat{u}}-\mathbf{\hat{u}}_{\min }\right\Vert _{E}\leq
C_{\epsilon ,s}\beta ^{s}+C\beta ^{\left( \nu +1\right) \left( 1-\epsilon
\right) }+C\varrho ^{-1}\beta ^{\left( \mu +1\right) \left( 1-\epsilon
\right) }+C\varrho .  \label{DD}
\end{equation}
\end{theorem}

\begin{proof}
We take $\mathbf{\hat{u}}=\sum\nolimits_{\vartheta
}\sum\nolimits_{l=1}^{N}u_{l,\vartheta }$ and estimate $\left\Vert \mathbf{%
\hat{u}}\left( \mathbf{k},\tau \right) -\mathbf{\hat{u}}_{\min }\left( 
\mathbf{k},\tau \right) \right\Vert _{E}$ applying subsequently Theorems \ref%
{Theorem uminw}, \ref{Theorem uminwav}, formulas (\ref{wPi}) and (\ref{vz}),
Theorem \ref{Theorem uminz} and finally Theorem \ref{Theorem moll} to obtain
inequality (\ref{DD}).
\end{proof}

Note that Theorem \ref{Theorem sumcanwave} is a direct corollary of Theorem %
\ref{Theorem canappr}.

\begin{remark}
Note that (\ref{eqcan}) is the Fourier integral version of the following
system of equations based on weakly universal \ nonlinearity and slightly
more general than (\ref{cansystext}) 
\begin{gather}
\partial _{\tau }u_{l,\vartheta }=\frac{1}{\varrho _{1}}\omega
_{n_{l}}^{\prime }\left( \mathbf{k}_{\ast l}\right) \cdot \nabla
_{x}u_{l,\vartheta }+\frac{i}{2\varrho _{2}}\nabla _{r}\cdot \omega
_{n_{l}}^{\prime \prime }\left( \mathbf{k}_{\ast l}\right) \nabla
_{r}u_{l,\vartheta }+f_{n_{l},\vartheta }^{\left( \mu ,\nu \right) }\left( 
\vec{u},\delta \vec{u}\right) ,  \label{canphys} \\
\left. u_{l\vartheta }\right\vert _{\tau =0}=\hat{H}_{l\vartheta },\text{
where }\delta _{i}u_{l}\left( \mathbf{r}\right) =u_{j}\left( \mathbf{r}%
+e_{i}\right) -u_{j}\left( \mathbf{r}-e_{i}\right)  \notag
\end{gather}%
where $\varrho _{1},\varrho _{1}$ are as in (\ref{r1r2}) and $e_{i}$ is $i$%
-th standard ort in $\mathbb{R}^{d}$. In the case when \ (\ref{wdisp}) holds 
$1/\varrho _{2}$ is bounded or small and the dependence on the coefficient $%
1/\varrho _{2}\ $is regular for small $\varrho $ and $\beta $ \ and $%
u_{\vartheta ,j}\left( \mathbf{k},\tau \right) $ may be looked at as a shape
function. When $\varrho _{1}=\varrho $ and $1/\varrho _{2}$ is substituted
by zero we obtain an equation exactly of the form (\ref{cansystext}).
\end{remark}

When $\nu =0$, $\mu =1$ and the $nk$-spectrum $S$ is universally resonance
invariant as in Definition \ref{Definition omclos}, the nonlinearities $%
f_{n_{l},\vartheta ,0}^{\left( 1,0\right) }$ are universal of the form (\ref%
{Ppm}). When the $nk$-spectrum $S$ is resonance invariant but not
universally resonance invariant, the nonlinearities are weakly universal,
but may be not universal, that allows, in particular, for the second and the
third harmonic generation.

\textbf{Acknowledgment:} Effort of A. Babin and A. Figotin is sponsored by
the Air Force Office of Scientific Research, Air Force Materials Command,
USAF, under grant number FA9550-04-1-0359.

\end{document}